\documentclass[11pt]{article}
\input amssymb.sty
\input amssym.def
\input amssym
\input epsf
\usepackage{graphics}

\newtheorem{theorem}{Theorem}[section]
\newtheorem{hypothesis}[theorem]{Hypothesis}
\newtheorem{theorem-definition}[theorem]{Theorem-Definition}
\newtheorem{theorem-construction}[theorem]{Theorem-Construction}
\newtheorem{lemma-definition}[theorem]{Lemma--Definition}
\newtheorem{proposition-definition}[theorem]{Proposition--Definition}
\newtheorem{lemma}[theorem]{Lemma}
\newtheorem{proposition}[theorem]{Proposition}
\newtheorem{corollary}[theorem]{Corollary}
\newtheorem{conjecture}[theorem]{Conjecture}
\newtheorem{definition}[theorem]{Definition}

\begin{document}
\newcommand{\be}{\begin{equation}}
\newcommand{\ee}{\end{equation}}
\newcommand{\bt}{\begin{theorem}}
\newcommand{\et}{\end{theorem}}
\newcommand{\bp}{\begin{proposition}}
\newcommand{\ep}{\end{proposition}}
\newcommand{\bl}{\begin{lemma}}
\newcommand{\el}{\end{lemma}}
\newcommand{\bc}{\begin{corollary}}
\newcommand{\ec}{\end{corollary}}
\newcommand{\bd}{\begin{definition}}
\newcommand{\ed}{\end{definition}}
\newcommand{\la}{\label}
\newcommand{\Z}{{\mathbb Z}}
\newcommand{\R}{{\mathbb R}}
\newcommand{\Q}{{\mathbb Q}}
\newcommand{\C}{{\mathbb C}}
\newcommand{\lms}{\longmapsto}
\newcommand{\lra}{\longrightarrow}
\newcommand{\hra}{\hookrightarrow}
\newcommand{\ra}{\rightarrow}
\newcommand{\sgn}{\rm sgn}
\begin{titlepage}
\title{Cluster ensembles, quantization  and  the dilogarithm}  
\author{V.V. Fock, A. B. Goncharov}
\end{titlepage}
\date{}
%\stepcounter{page}
\maketitle
%\addtocounter{page}{+1}
%\doublespace

\tableofcontents

 \vskip 6mm \noindent

\section{Introduction and main definitions with simplest examples}

Cluster algebras are a remarkable discovery of 
S. Fomin and A. Zelevinsky 
\cite{FZI}. They are  certain commutative algebras defined 
by a very simple and general data. 
%Cluster algebras of geometric origin are quite 
%general and probably the most important examples of cluster algebras. 

We show that a cluster algebra  %of geometric origin 
is part of a richer structure, which we call a 
{\it cluster ensemble}. A cluster ensemble 
is a pair $({\cal X}, {\cal A})$ of {\it positive spaces} 
(which are 
varieties equipped with positive atlases), coming 
with an action of a certain  discrete symmetry group $\Gamma$. 
These two spaces are 
related by a morphism $p: {\cal A} \lra {\cal X}$, which in general, 
as well as in many interesting examples, is neither 
injective nor surjective. 
The space ${\cal A}$ has a degenerate symplectic structure, and 
the space  ${\cal X}$ has a Poisson structure. The map $p$ relates the 
Poisson and degenerate symplectic  structures in a natural way. 
Amazingly, the dilogarithm together with 
 its motivic and quantum avatars plays a central role 
in the cluster ensemble structure.  The space  
${\cal A}$ is closely related to 
the spectrum 
of a cluster algebra. 
On the other hand, in many situations
  the most interesting part 
of the structure is the space ${\cal X}$.

We define a canonical
 non-commutative $q$-deformation of the ${\cal X}$-space. 
We show that when $q$ is a root of unity the algebra of functions on the 
$q$-deformed ${\cal X}$-space has a large center, 
identified with the algebra of functions on the original ${\cal X}$-space. 
Cluster ensembles admit canonical  quantization.  

\vskip 2mm
The main example, as well as the main application of this theory so far,  
is provided by the $({\cal X}, {\cal A})$-pair 
of moduli spaces assigned in \cite{FG1} to a  
topological surface $S$ with a finite set of points at the boundary 
 and a semisimple 
algebraic group $G$. In particular, the ${\cal X}$-space 
in the simplest 
case when $G=PGL_2$ and $S$ is a disc with $n$ points 
at the boundary is the moduli space ${\cal M}_{0, n}$. 

This pair of moduli spaces 
is an algebraic-geometric avatar of higher Teichm{\"u}ller theory on $S$ related
to $G$. In the case $G=SL_2$ we get the classical Teichm{\"u}ller theory, 
as well as its generalization 
to surfaces with a finite set of points on the boundary. 
A survey of the Teichm{\"u}ller theory emphasizing
 the cluster point of view can be found in \cite{FG}. 

\vskip 2mm
We suggest  that there exists a  
duality between the  ${\cal A}$ and ${\cal X}$ spaces.  
One of its manifestations is our package of duality conjectures 
 in Section 4. These conjectures assert 
that the {\it tropical points} of the ${\cal A}/{\cal X}$-space 
parametrise a basis in a certain class of 
functions on the {\it Langlands dual} ${\cal X}/{\cal A}$-space. 
It can be viewed as  a canonical function (the {\it universal kernel}) on the product of the 
set of tropical points of one space and the 
Langlands dual space.

To support these conjectures, we define in Section 5.1 
the tropical limit of such a universal kernel in the finite type case. 
Another piece of evidence is provided by Chapter 12 in \cite{FG1}. 

\vskip 2mm
In the rest of the Introduction we define cluster ${\cal X}$- and 
${\cal A}$-varieties and describe  
their key 
features. Section 1.1 provides background on positive
spaces, borrowed from Chapter 4 of \cite{FG1}. 
Cluster varieties are defined in Section 1.2. 
In Section 1.3 we discuss one of the simplest examples: cluster 
${\cal X}$-variety structures of the moduli space ${\cal M}_{0, n+3}$. 
In Section 1.4 we summarize their main structures. In Section 1.4 we discuss how they 
appear in our main example - higher Teichm\"uller theory.

\subsection{Positive schemes and positive spaces} 
A {\it semifield}  is a   
 set $P$  equipped with the  operations of
 addition and multiplication, so that 
addition is commutative and associative, multiplication makes $P$ an abelian group,
and they are compatible in a usual way: $(a+b)c = ac + bc$ for $a,b,c \in P$. 
A standard example is given by the set $\R_{>0}$ of positive real numbers. 
Here are more exotic examples.  
 Let ${\Bbb A}$ be one of the sets  
$\Z$, $\Q$ or $\R$. The {\it tropical semifield} ${\Bbb A}^t$ associated with ${\Bbb A}$ 
is the set ${\Bbb A}$  with the multiplication $\otimes $ 
and addition $\oplus$ given by 
$$
a \otimes b := a+b, \quad a \oplus b := {\rm max}(a,b).
$$
One more example is given by the semifield $\R_{>0}((\varepsilon))$ of Laurent series in $\varepsilon$ with 
real coefficients and a positive leading coefficient, 
equipped with the usual addition and multiplication. There is a homomorphism of semifields
$-{\rm deg}: \R_{>0}((\varepsilon))\to \Z^t$, given by $f \lms -{\rm deg}(f)$. 
It explains the origin of the tropical semifield $\Z^t$. 

\vskip 3mm
Recall the standard notation ${\Bbb G}_m$ for the multiplicative group. 
It is an affine algebraic group. 
The ring of regular functions on ${\Bbb G}_m$ is $\Z[X,X^{-1}]$, and for any field $F$ one has 
${\Bbb G}_m(F)= F^*$. A product of multiplicative groups 
is known as a {\it  split algebraic torus} over $\Z$, or simply a split algebraic torus. 

Let $H$ be a split algebraic torus. 
 A rational function $f$ on $H$ is called {\it positive}  
if it belongs to the semifield generated, in the field of rational functions on $H$, 
 by the characters of 
$H$. 
So it  can be written 
as $f = f_1/f_2$ where  $f_1, f_2$ are linear combinations of characters 
with positive integral coefficients. 
A {\it positive rational map} between two split tori $H_1, H_2$ 
is a rational map $f: H_1 \to H_2$ such that $f^*$ induces a homomorphism of the semifields  
of positive rational functions. Equivalently, for any character $\chi$ 
of $H_2$ the composition $\chi \circ f$ is a positive rational function 
on $H_1$. A composition of positive rational functions is positive. 
Let ${\rm Pos}$ be the category whose objects are split algebraic tori
and morphisms are positive rational maps. A {\it positive divisor} in a torus $H$ is a divisor given by 
an equation $f=0$, where $f$ is a positive rational function on $H$.

\begin{definition} \label{7.28.03.1} 
A { positive atlas} on an irreducible scheme/stack $X$ over $\Q$ 
is a family of  birational isomorphisms
\begin{equation} \label{7.29.03.2}
\psi_{\alpha}: H_{\alpha} \lra X, \quad \alpha \in {\cal C}_X,
\end{equation}
between split algebraic tori $H_{\alpha}$ and $X$, 
parametrised by a non empty set ${\cal C}_X$, 
 such that: 

i) each $\psi_{\alpha}$ is regular on the complement of a positive divisor 
in $H_{\alpha}$;

ii)  for any  ${\alpha}, {\beta} \in {\cal C}_X$ the map 
$
\psi_{\alpha, \beta}:= \psi_{\beta}^{-1} \circ \psi_{\alpha}: H_{\alpha} \lra H_{\beta}
$
is a positive rational map\footnote{A positive atlas covers a non-empty Zariski 
open subset of $X$, but not necessarily the whole space $X$}. 

\noindent
A positive atlas is called regular if each $\psi_{\alpha}$ is regular. 
\end{definition}

Birational isomorphisms (\ref{7.29.03.2}) are 
called positive coordinate systems on $X$.  
A {\it positive scheme} is 
a scheme equipped with a positive 
atlas. 
We will need an equivariant version of this definition. 
\begin{definition} \label{7.28.03.1s} 
Let $\Gamma$ be a group of automorphisms of  $X$. A 
positive atlas (\ref{7.29.03.2}) on $X$ is {\it $\Gamma$-equivariant} if 
 $\Gamma$ acts on the set ${\cal C}_X$, and 
for every $\gamma \in \Gamma$ there is an isomorphism of algebraic tori 
$i_{\gamma}: H_{\alpha} \stackrel{\sim}{\lra} H_{\gamma(\alpha)}$ making the 
following diagram commutative:
\begin{equation} \label{7.20.03.1}
\begin{array}{ccc}
H_{\alpha} & \stackrel{\psi_{\alpha}}{\lra}& X\\
\downarrow i_{\gamma}& &\downarrow \gamma\\
H_{\gamma(\alpha)}&\stackrel{\psi_{\gamma(\alpha)}}{\lra}  & X
\end{array}
\end{equation}
\end{definition}

Quite often a collection of positive coordinate systems is the only data 
we need when working with a positive  scheme.  
Axiomatizing this observation, we arrive at the {\it category of 
positive spaces} defined below.   

A groupoid is a category where all morphisms are isomorphisms. 
We assume that the set of morphisms between any two objects  is non-empty.  
The 
{\it fundamental group  of a groupoid} is 
the automorphism group of an object of  the groupoid. 
 It is well defined up to an inner automorphism.

 \begin{definition} \label{10.31.03.1} Let ${\cal G}_{\cal X}$ be a groupoid. 
A {\em positive space} is a functor 
\begin{equation} \label{11.6.03.99}
\psi_{\cal X}: {\cal G}_{\cal X}\lra {\rm Pos}.
\end{equation}
\end{definition}
The groupoid ${\cal G}_{\cal X}$ is 
called the {\it coordinate groupoid} of a positive
space. 
Thus for every object $\alpha$ of   ${\cal G}_{\cal X}$ 
there is an algebraic torus $H_{\alpha}$, called a {\it coordinate torus} of the
positive space ${\cal X}$, 
and for every morphism $f: \alpha \lra \beta$ in the groupoid 
there is a positive birational isomorphism 
$\psi_{f}: H_{\alpha} \lra  H_{\beta}$.

Let $\psi_1$ and $\psi_2$ be functors from coordinate groupoids ${\cal G}_1$
and ${\cal G}_2$ 
to the category ${\rm Pos}$. A morphism 
from  $\psi_1$ to $\psi_2$ is a pair consisting of a functor $\mu: {\cal G}_2
\to {\cal G}_1$ and 
a natural transformation $F: \psi_2 \to \psi_1\circ \mu$. 
A morphism  is called a {\it monomial morphism}  if 
for every object $\alpha \in {\cal G}_2$ the map $F_{\alpha}: \psi_2(\alpha) \to 
\psi_1(\mu(\alpha))$ is a homomorphism of algebraic tori. 
\vskip 3mm

{\bf Example 1}. A positive variety $X$  provides a functor 
(\ref{11.6.03.99}) as follows.  The fundamental group of the coordinate
groupoid ${\cal G}_{\cal X}$ is trivial, so it is just  a
set. 
Precisely, the objects of ${\cal G}_{\cal X}$ form the set ${\cal C}_X$
 of 
coordinate charts of the positive atlas on  $X$. The morphisms are given by 
the subset of ${\cal C}_X \times {\cal C}_X$ consisting of 
pairs of charts with nontrivial intersection, with the obvious source and target maps. In particular, the morphisms form  the 
set ${\cal C}_X \times {\cal C}_X$ if $X$ is irreducible. 
The functor $\psi_{\cal X}$ is given by  $\psi_{\cal X}(\alpha):= H_{\alpha}$  and 
 $\psi_{\cal X}(\alpha\to \beta):= \psi_{\alpha, \beta}$.

{\bf Example 2}. A $\Gamma$-equivariant positive scheme $X$ 
 provides a positive space ${\cal X}$ given by a functor 
(\ref{11.6.03.99}). The fundamental group of its coordinate groupoid  
 is isomorphic to $\Gamma$.   

\vskip 3mm
Given a split torus $H$ and a  semifield $P$ we define the set of $P$-valued points of $H$ as 
$$
H(P):= X_*(H) \otimes_{\Z} P,
$$ 
where $X_*(H)$ is the group of cocharacters of
$H$, and the tensor product is with the abelian group defined by the 
semifield $P$.  
A positive birational isomorphism $\psi: H \to H'$ induces 
a map $\psi_*: H(P) \to H'(P)$. 

\vskip 2mm
{\bf Example 3.} If ${\Bbb A}^t$ is a tropical semifield, then 
the map $\psi_*$ is given by a piece-wise 
linear map, the {\it tropicalization} of the map $\psi$. 
\vskip 2mm

An inverse to a positive map may not be positive -- the inverse 
of the map $x'=x+y, y'=y$ is $x=x'-y', y=y'$. 
If $\psi^{-1}$ is also positive, 
the map $\psi_*$ is an isomorphism. 

Given a positive space ${\cal X}$ there is a unique set 
${\cal X}(P)$ of $P$-points of ${\cal X}$.
 It can be defined as 
$$
{\cal X}(P) = \coprod_{\alpha}H_\alpha(P)/(\mbox{identifications $\psi_{\alpha, \beta \ast}$}).
$$ 
For every object $\alpha$ of the coordinate groupoid 
${\cal G}_{\cal X}$ there are functorial (with respect to the maps 
${\cal X}\to {\cal X}'$) isomorphisms
$$
{\cal X}(P) \stackrel{\sim}{=} H_{\alpha}(P) \stackrel{\sim}{=}P^{\rm dim {\cal X}}.
$$
Therefore the fundamental group of the coordinate groupoid 
acts on the set ${\cal X}(P)$.

\vskip 3mm
A  positive space ${ {\cal X}}$ gives rise  
to a prescheme ${\cal X}^{*}$. It is 
obtained by gluing the tori $H_{\alpha}$, where 
$\alpha$ runs through the objects of ${\cal G}_{\cal X}$, according to the 
birational maps $\psi_{f}$ corresponding to morphisms $f: \alpha \to \beta$. 
It, however, may not be separable, and thus may not be a scheme. Each torus $H_{\alpha}$ 
embeds to ${\cal X}^{*}$ as a Zariski open dense subset $\psi_{\alpha}(H_{\alpha})$.

The ring 
of regular functions 
on  ${\cal X}^{*}$ is called the ring of {\it 
universally Laurent polynomials} for 
${\cal X}$ and denoted by ${\Bbb L}({\cal X})$.  In simple terms, 
the ring ${\Bbb L}({\cal X})$ consists of all 
rational functions which are 
regular at every coordinate torus $H_{\alpha}$.

It is often useful to take the affine closure of ${\cal X}^{*}$, 
understood as the spectrum ${\rm Spec}({\Bbb L}({\cal X}))$. 

The positive 
structure on ${\cal X}^{*}$ provides the semifield of all positive rational 
functions on ${\cal X}^{*}$. Intersecting it with the ring 
${\Bbb L}({\cal X})$ we get the semiring $\widetilde {\Bbb L}_+({\cal X})$.
As an example $1-x+x^2 = (1+x^3)/(1+x)$ shows,  a rational 
function can be positive, while the coefficients of the corresponding Laurent polynomial may be not. So we define a smaller semiring ${\Bbb L}_+({\cal X})$ of {\it positive 
universally Laurent polynomials} for 
${\cal X}$ as follows:  an element of  ${\Bbb L}_+({\cal X})$ 
is a rational function on ${\cal X}^{*}$ whose restriction 
to one (and hence any)  of the embedded coordinate tori $\psi_{\alpha}(H_{\alpha})$  is a 
linear combination of characters of this torus with positive integral coefficients.

\subsection{Cluster ensembles: definitions} They are defined by a combinatorial data --- seed --- similar 
\footnote{although different in detail-- we do not include 
the cluster coordinates in the definition of a seed, and give a coordinate free definition} 
to  the one used in the definition of 
 cluster algebras \cite{FZI}.

\paragraph{Seeds and seed tori.}  
Recall that a lattice is a free abelian 
group. 
\begin{definition} \label{esd}
A seed
  is a
 datum
$
(\Lambda, (\ast, \ast), \{e_i\}, \{d_i\}),
$ 
 where

i) $\Lambda$ is a lattice; 

ii) $(\ast, \ast)$ is a skewsymmetric $\Q$-valued bilinear form 
on $\Lambda$; 

iii) $\{e_i\}$ is a basis  of the lattice $\Lambda$,  and $I_0$ 
is a subset of basis vectors, called frozen basis vectors; 

iv) $\{d_i\}$ are positive integers assigned to the basis vectors, such that 
$$
{\varepsilon}_{ij} := (e_i, e_j)d_j \in \Z \quad 
\mbox{unless 
$i,j \in I_0\times I_0$}. 
$$
\end{definition}
The numbers $\{d_i\}$ are called  the {\em multipliers}. 
We assume that their greatest common divisor is $1$. 

\begin{figure}[ht]
\centerline{\epsfbox{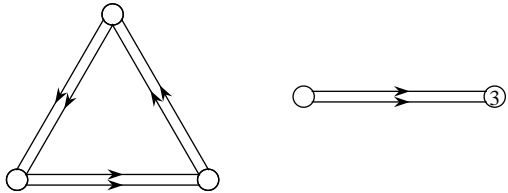}}
\caption{Picturing seeds by quivers - we show $d_i$'s only if they differ from $1$.}
\label{quiver}
\end{figure}

\vskip 2mm
{\it Seeds as quivers.} 
A seed is a version of the notion of a {\it quiver}. 
Precisely, let us assume for simplicity 
that the set of frozen basis vectors is empty. 
A quiver corresponding to a seed is a graph whose set of vertices  $\{i\}$ 
is identified with the set of basis vectors $\{e_i\}$; 
two vertices $i,j$ are connected 
by  $|(e_{i}, e_j)|$ arrows going from $i$ to $j$ 
if  $(e_{i}, e_j)>0$, and from $j$ to $i$ otherwise;  
 the $i$-th vertex is  marked by $d_i$, see Fig. \ref{quiver}. 
The (enhanced by multipliers) quivers we get 
have the following property: all arrows between any two vertices 
are oriented the same way, and there are no arrows from a vertex to itself. 
Clearly any enhanced quiver like this corresponds to a unique seed. 

{\it Lattices and split algebraic tori}. Recall that a lattice $\Lambda$ gives rise to a split algebraic torus 
$$
{\cal X}_{\Lambda}:= {\rm Hom}(\Lambda, {\Bbb G}_m). 
$$
The set of its points with values in a field $F$ is the group ${\rm Hom}(\Lambda, F^*)$. 

An element $v \in \Lambda$ provides a character $X_v$ of ${\cal X}_{\Lambda}$. Its value 
on a homomorphism $x \in {\cal X}_{\Lambda}$ is $x(v)$. The assignment 
$
\Lambda \lra {\cal X}_{\Lambda}
$ 
is a contravariant functor providing an equivalence of categories
$$
\mbox{the dual to the category of finite rank lattices} 
\stackrel{\sim}{\lra} \mbox{the category of split algebraic tori}. 
$$  
The inverse functor assigns to a split algebraic torus ${\cal T}$ its lattice of characters 
${\rm Hom}({\cal T}, {\Bbb G}_m)$.

There is the dual lattice 
 $$
\Lambda^*:= {\rm Hom}(\Lambda, \Z).
$$
An element $a \in \Lambda^*$ gives rise to a cocharacter
$$
\varphi_a: {\Bbb G}_m \lra {\rm Hom}(\Lambda, {\Bbb G}_m).
$$
On the level of $F$-points, $\varphi_a(f)$ is the homomorphism $v \lms f^{a(v)}$.

\vskip 2mm
The {\it seed ${\cal X}$-torus}  
is a split algebraic torus
$
{\cal X}_{\Lambda}:= {\rm Hom}(\Lambda, {\Bbb G}_m) 
$. 
It carries  a  Poisson structure  
provided by the form $(\ast, \ast)$:
$$
\{X_v, X_w\} = (v,w)X_v X_w.
$$
The basis $\{e_i\}$ provides {\it cluster ${\cal X}$-coordinates $\{X_i\}$}. 
They form a basis in the group of  characters of ${\cal X}_{\Lambda}$.

\vskip 2mm

The basis $\{e_i\}$ provides a dual basis $\{e^{\ast}_i\}$ of 
the lattice $\Lambda^*$.
We need a quasidual basis 
 $\{f_i\}$ given by 
\be \la{TT}
f_i = d_i^{-1}e^{\ast}_i. 
\ee
Let $\Lambda^{\circ}\subset \Lambda^*\otimes \Q$ be 
the sublattice spanned by the vectors $f_i$. 
The {\it seed ${\cal A}$-torus}  
is a split algebraic  torus
$$
{\cal A}_{\Lambda}:= {\rm Hom}(\Lambda^\circ, {\Bbb G}_m).
$$
The basis 
 $\{f_i\}$ provides {\it cluster ${\cal A}$-coordinates $\{A_i\}$}.

\vskip 2mm
Let ${\cal O}(Y)^*$ be the group of invertible regular functions on a variety $Y$.
There is a map
$$
d\log \wedge d\log: 
{\cal O}(Y)^*\wedge {\cal O}(Y)^* \lra \Omega^2(Y), \quad f\wedge g \lms d\log(f)\wedge  d\log(g).
$$
The skew-symmetric bilinear form $(\ast, \ast)$, viewed as an element of  
$\bigwedge^2\Lambda^\circ_{\bf i}$, provides an element 
$$
W \in {\cal O}({\cal A}_{\Lambda})^*\wedge {\cal O}({\cal A}_{\Lambda})^*.
$$
Applying the map $d\log \wedge d\log$ to $W$ we get a closed $2$-form $\Omega$ on the torus ${\cal A}_{\Lambda}$:
$$
\Omega := d\log \wedge d\log (W) \in \Omega^2({\cal A}_{\Lambda}). 
$$
\vskip 2mm

There is a non-symmetric bilinear form $[\ast, \ast]$ on the lattice 
$\Lambda$, defined by setting 
\be \la{[]}
[e_i, e_j]:= (e_i, e_j)d_j.
\ee 
There is a natural map of lattices 
$$
p^*: \Lambda \lra \Lambda^\circ, \qquad v\lms \sum_j (v, e_j)e_j^{\ast}  = \sum_j [v, e_j]f_j.
$$ 
It gives rise to a homomorphism of seed tori
\be \la{pi}
p: {\cal A}_{\Lambda} \lra {\cal X}_{\Lambda}. 
\ee
The following lemma is straightforward.

\begin{lemma} \label{1.18.04.1}
The fibers of the map $p$ are the leaves of the null-foliation of
the $2$-form $\Omega$. The subtorus
${\cal U}_{\Lambda}:= {p}({\cal A}_{\Lambda})$ is a
symplectic leaf of the Poisson structure on ${\cal X}_{\Lambda}$. 

The symplectic structure on ${\cal U}_{\Lambda}$ induced by the form $\Omega$ on ${\cal A}_{\Lambda}$ 
coincides with
the symplectic structure given by the restriction of the Poisson structure on ${\cal X}_{\Lambda}$.
\end{lemma}

Summarising, a seed ${\mathbf i}$ gives rise to  seed ${\cal X}$- and ${\cal A}$-tori. 
Although they depend only on the lattice $\Lambda$, 
we denote them by ${\cal X}_{\bf i}$ and ${\cal A}_{\bf i}$ to emphasize 
the cluster coordinates on these tori provided by the seed ${\bf i}$.

\paragraph{Seed mutations.} Set   $[\alpha]_+= \alpha$ if $\alpha\geq 0$ 
and $[\alpha]_+=0$ otherwise. So $[\alpha]_+ = {\rm max}(0, \alpha)$. 

Given a seed ${\bf i}$ and a non-frozen basis vector $e_k$, we define a new seed 
${\bf i}'$, called {\it the seed obtained from ${\bf i}$ by mutation in the 
direction of a non-frozen basis vector $e_k$}. The seed  ${\bf i}'$ 
is obtained by changing the basis $\{e_i\}$ -- the rest of the datum stays 
the same. 
The new basis $\{e'_i\}$ is
\begin{equation} \label{12.12.04.2a}
e'_i := 
\left\{ \begin{array}{lll} e_i + [\varepsilon_{ik}]_+e_k
& \mbox{ if } &  i\not = k\\
-e_k& \mbox{ if } &  i = k.\end{array}\right.
\end{equation}
We denote by $\mu_{e_k}({\bf i})$, or simply by $\mu_{k}({\bf i})$,  
the seed ${\bf i}$ mutated in the direction 
of a basis vector $e_k$. By definition, the frozen/non-frozen basis vectors 
of the mutated seed are the images of the  frozen/non-frozen basis vectors of the original seed. 

The basis $\{f_i\}$ in $\Lambda^\circ$ mutates as follows:
\be \la{qdbas}
f'_i := 
\left\{ \begin{array}{lll} -f_k + \sum_{j\in I}[-\varepsilon_{kj}]_+f_j
& \mbox{ if } &  i = k\\
f_i& \mbox{ if } &  i \not = k.\end{array}\right.
\ee
Therefore,  although the definition 
of the lattice $\Lambda^\circ$ involves a choice of a seed, the lattice 
does not depend on it.  

\vskip 3mm

{\bf Remark}. The basis $\mu_k^2(\{e_i\})$ does not necessarily 
coincide with 
$\{e_i\}$. For example, let $\Lambda$ be a rank two lattice with a basis 
$\{e_1, e_2\}$, and $(e_1, e_2)=1$. Then 
$$
\{e_1, e_2\} \stackrel{\mu_2}{\lra} \{e_1+e_2, -e_2\} \stackrel{\mu_2}{\lra} 
\{e_1+e_2, e_2\}.  
$$
However, although the seeds  $\mu_k^2({\bf i})$ and ${\bf i}$ are different, 
they are canonically isomorphic. 

\paragraph{Coordinate description.} 
\begin{definition} \label{9.29.04.1} A {\em seed} ${\mathbf i}$ 
is a quadruple $(I, I_0, \varepsilon, d)$, 
where

i) $I$ is a finite set, and $I_0$ is a subset of $I$; 

ii) $\varepsilon = \varepsilon_{ij}$ is a $\Q$-valued function on $I \times I$, such that
$\varepsilon_{ij} \in {\mathbb Z}$, unless $(i,j) \in I_0\times I_0$;

iii) $d = \{d_i\}$, where $i \in I$, is a set of positive rational numbers,
such that the function $$
\widehat{\varepsilon}_{ij}=\varepsilon_{ij}d^{-1}_j ~~\mbox{is skew-symmetric: $\widehat{\varepsilon}_{ij} =- \widehat{\varepsilon}_{ji}$}. 
$$ 

\end{definition}

 Definitions \ref{esd} and \ref{9.29.04.1} are equivalent. Indeed, 
given a seed from Definition \ref{9.29.04.1}, we set 
$$
\Lambda := \Z[I], \quad e_i:= \{i\}, ~~i\in I, \quad (e_i, e_j):= 
{\varepsilon}_{ij}~d^{-1}_j. 
$$
The non-symmetric bilinear form is the function ${\varepsilon}$: 
$$
[e_i, e_j]= 
{\varepsilon}_{ij}. 
$$
The function $\varepsilon$
is called the 
{\em exchange function}. The numbers $\{d_i\}$ are the {multipliers}. 
The subset $I_0 \subset I$ is the {\em frozen subset} of $I$, 
and its elements are the {\it frozen elements} of $I$. 
Elements of the set $I$ are often called {\it vertices}. 

\vskip 2mm

The Poisson structure on the torus ${\mathcal X}_{\mathbf i}$ looks in coordinates as follows: 
\begin{equation} \label{4.30.03.2}
\{X_i, X_j\} = \widehat{\varepsilon}_{ij}X_iX_j,\qquad \widehat \varepsilon_{ij }:= \varepsilon_{ij }d_j^{-1}.
\end{equation}

The $2$-form $\Omega$ on the  torus ${\mathcal A}_{\mathbf i}$
 is
\begin{equation} \label{4.30.03.2WP}
\Omega = \sum_{i,j \in I}\widetilde{\varepsilon}_{ij}d\log A_i \wedge d\log A_j, \qquad  \widetilde \varepsilon_{ij}:= d_i
\varepsilon_{ij}.
\end{equation}

The homomorphism $p$ -- see (\ref{pi}) --  is given by  
\begin{equation} \label{9.28.04.2}
p: {\mathcal A}_{\mathbf i} \lra 
{\mathcal X}_{\mathbf i}, \qquad p^*X_i = \prod_{j \in I}A_j^{\varepsilon_{ij}}.
\end{equation}

Given a seed ${\mathbf i}=(I, I_0,\varepsilon, d)$, every non-frozen 
element $k \in I-I_0$ provides  a {mutated in the 
direction $k$} seed $\mu_k({\mathbf i}) = {\mathbf i'} =(I', I'_0,\varepsilon', d')$: one has $I':= I, I_0':= I_0, d':=d$ and
\begin{equation} \label{5.11.03.6}
 \varepsilon'_{ij} := \left\{ \begin{array}{lll} 
- \varepsilon_{ij} & \mbox{ if $k \in \{i,j\}$} \\ 
\varepsilon_{ij} & \mbox{ if $\varepsilon_{ik}
\varepsilon_{kj} \leq 0, \quad k \not \in \{i,j\}$} \\
\varepsilon_{ij} + |\varepsilon_{ik}| \cdot \varepsilon_{ kj}& 
\mbox{ if $\varepsilon_{ik}
\varepsilon_{kj} > 0, \quad k \not \in \{i,j\}.$}\end{array}\right.
\end{equation} 
This procedure is involutive: 
the mutation of  $\varepsilon'_{ij}$ at the vertex $k$ is 
the original function $\varepsilon_{ij}$.

This definition of mutations is equivalent to the coordinate free definition thanks to the following Lemma.
\begin{lemma} \label{12.12.04.1a}
One has 
 $\varepsilon'_{ij} = (e'_i, e'_j)d_j$, where $\varepsilon'_{ij}$ 
is given by formula $(\ref{5.11.03.6})$. 
\end{lemma}

{\bf Proof}. Clearly 
$
(e'_i, e'_k)  = (e_i + [\varepsilon_{ki}]_+e_k, - e_k) 
= - \widehat \varepsilon_{ik} = \widehat \varepsilon'_{ik}.  
$ 
Assume that $k \not \in \{i,j\}$. Then 
$$
(e'_i, e'_j) 
= (e_i + [\varepsilon_{ik}]_+e_k, e_j + [\varepsilon_{jk}]_+e_k) 
= 
\widehat \varepsilon_{ij} + [\varepsilon_{ik}]_+\widehat \varepsilon_{kj}
+ \widehat \varepsilon_{ik}[\varepsilon_{jk}]_+  
$$
$$
=\widehat \varepsilon_{ij} + [\varepsilon_{ik}]_+\widehat \varepsilon_{kj}
+ \varepsilon_{ik}[-\widehat \varepsilon_{kj}]_+ = \widehat
\varepsilon'_{ij}. 
$$
The lemma is proved.

\paragraph{Cluster transformations.} 
This is the heart of the story.  
A seed mutation $\mu_k$ induces positive 
 rational maps between the corresponding 
 seed ${\cal X}$- and ${\cal A}$-tori, denoted by the same
 symbol $\mu_k$. Namely, 
 denote the cluster coordinates related to 
the seed $\mu_k({\bf i})$ by $X_i'$ and $A_i'$. Then we define 
\begin{equation} \label{5.11.03.1x}
\mu_k^*X'_{i} := \left\{\begin{array}{ll} X_k^{-1}& \mbox{ if }  i=k \\
 X_i(1+X_k^{-{\rm sgn} (\varepsilon_{ik})})^{-\varepsilon_{ik}} & \mbox{ if }   i\neq k,
\end{array} \right.
\end{equation}
\begin{equation} \label{5.11.03.1a}
  A_{k}\cdot \mu_k^*A'_{k} := \quad \prod_{j| \varepsilon_{kj} >0} 
A_{j}^{\varepsilon_{kj}} + \prod_{j| \varepsilon_{kj} <0} 
A_{j}^{-\varepsilon_{kj}}; \qquad \mu_k^*A'_{i} =  
A_{i}, \quad i \not = k.
\end{equation}
Here 
if just one of the sets $\{j| \varepsilon_{kj} >0\}$ and 
$\{j| \varepsilon_{kj} < 0\}$ is empty, the corresponding monomial 
is $1$. If $\varepsilon_{kj} =0$ for every $j$, the 
right hand side of the formula is $2$, and $\mu_k^*X'_{k} =X_k^{-2}$, 
$\mu_k^*X'_{i} =X_i$ for $i \not = k$.

Seed isomorphisms $\sigma$ obviously  induce isomorphisms between the corresponding 
 seed tori, which are  denoted by the same
 symbols $\sigma$:
\begin{equation} \label{6.4.07.1}
\sigma^*X'_{\sigma(i)}=X_{i}, \qquad \sigma^*A'_{\sigma(i)} = A_{i}.
\end{equation}

A {\em seed cluster transformation} is a 
composition of seed isomorphisms and mutations. 
It gives rise 
to a {\em cluster transformation} of the corresponding seed
${\cal X}$- or ${\cal A}$-tori. The latter is a rational map 
obtained by the composition of 
 isomorphisms and mutations corresponding to the 
seed isomorphisms and mutations forming  
the seed cluster transformation. 
Given a semifield $P$, cluster transformations induce 
isomorphisms between the sets of 
$P$-points of the corresponding cluster tori. 
 
Two seeds are
called {\em equivalent} if they are related 
by a cluster transformation. The equivalence class of a
seed ${\mathbf i}$ is denoted by $|{\mathbf i}|$.

\vskip 2mm
Mutation formulas (\ref{5.11.03.6})  and (\ref{5.11.03.1a}) were 
invented by Fomin and Zelevinsky \cite{FZI}. 
Clearly  the functions obtained by cluster ${\cal A}$- (resp. ${\cal X}$-) 
transformations from 
the coordinate functions on the initial seed ${\cal A}$-torus 
(resp. seed ${\cal X}$-torus)
are positive rational functions 
on this torus. The rational functions $A_i$ obtained this way 
generate the cluster algebra.

\paragraph{Cluster modular groupoids.} A seed cluster transformation ${\bf i}\to {\bf i}$ 
is called {\it trivial}, if 
the corresponding maps of the seed ${\cal A}$-tori as well as of the seed 
${\cal X}$-tori are the identity maps. \footnote{We conjecture that 
one of them implies the identity of the other.}

We define the {\it cluster modular  groupoid} ${\cal G}_{|{\bf i}|}$ as a groupoid 
whose objects are seeds equivalent to a given seed ${\bf i}$, and 
morphisms are cluster transformations modulo the trivial ones. 
The fundamental 
group $\Gamma_{{\bf i}}$ of this groupoid (based at ${\bf i}$) 
is called 
the {\it cluster modular group}.

\paragraph{The ${\cal A}$- and ${\cal X}$- positive spaces.} 
We have defined three categories. 
The first is the groupoid ${\cal G}_{|{\bf i}|}$. The other two have seed ${\cal A}$-/${\cal X}$-tori as objects and 
cluster transformations of them as morphisms. 
There are canonical functors from the first category to the second and third. 
They provide a pair of positive spaces 
of the same dimension, 
denoted by ${\cal A}_{|{\bf i}|}$ 
and ${\cal X}_{|{\bf i}|}$, which 
share a common coordinate groupoid ${\cal G}_{|{\bf i}|}$. 
We skip the subscript $|{\bf i}|$ whenever possible, writing 
 ${\cal X}$ for ${\cal X}_{|{\bf i}|}$ etc.

\paragraph{Examples of trivial cluster transformations.}  
Given a seed ${\bf i}$, denote by 
$\sigma_{ij}({\bf i})$  a new seed induced by the map of sets 
$I \to I$  
interchanging $i$ and $j$.  

\bp \la{FTFTFT} 
Let $h=2,3,4,6$  when $p = 0, 1, 2, 3$ respectively. 
Then if $\varepsilon_{ij}=-p\varepsilon_{ji}=-p$, 
\begin{equation} \label{K10}
(\sigma_{ij}\circ \mu_i)^{h+2} = \mbox{a trivial cluster transformation}. 
\end{equation}
\ep
Relations (\ref{K10}) are affiliated with the rank two 
Dynkin diagrams, i.e. $A_1 \times A_1, A_2, B_2, G_2$. 
The number $h=2,3,4,6$ is the Coxeter 
number of the diagram. One can present these relations in the form
$$
\mu_i \circ \mu_j \circ \mu_i \circ \mu_j \circ \ldots \stackrel{\sim}{=} \sigma_{ij}^{h+2},
$$
where the number of mutations on the left equals $h+2$. Notice that 
the right hand side is the identity in all but $A_2$ cases. 
We prove Proposition \ref{FTFTFT} in Section \ref{sec2.5}. 
%See Lemma \ref{7.7.04.1r} and the Example after it. 

\paragraph{Special 
cluster modular groupoid and modular groups.} 
{\it Special trivial seed cluster transformations}  
are compositions of 
the one given by (\ref{K10}) and isomorphisms. 
We do not know any other general procedure 
to generate trivial cluster transformations.

\begin{definition} \la{smg} {\rm Special 
cluster modular groupoid $\widehat {\cal G}$} is a connected groupoid 
whose objects are seeds, and 
morphisms are cluster transformations modulo the special trivial ones. 

The fundamental 
group $\widehat \Gamma$ of the groupoid $\widehat {\cal G}$ is called 
the {\rm special  cluster modular group}. 
\end{definition}

So there is a canonical functor $\widehat {\cal G}\to {\cal G}$ 
inducing a surjective map $\widehat \Gamma \to \Gamma$. 

The groupoid $\widehat {\cal G}$ has a natural  geometric interpretation,  
which justifies Definition \ref{smg}. 
Namely, thanks to Theorem \ref{12.9.03.2} the group $\Gamma$ acts, with finite 
 stabilizers, 
on a certain manifold with a polyhedral decomposition. 
So it acts on the dual polyhedral complex $\widehat M$, 
called {\it the modular complex}. The groupoid $\widehat {\cal G}$ 
is identified with the fundamental groupoid of this polyhedral complex, see  
Theorem \ref{11.2.03.1}.

\paragraph{The cluster ensemble.} We show (Proposition \ref{7.7.04.12}) that cluster transformations commute with the map
$p$  -- see (\ref{9.28.04.2}). So the 
map $p$  
 gives rise 
to a monomial morphism of positive spaces 
\begin{equation} \label{9.28.04.3}
p: {\cal A} \lra {\cal X}.
\end{equation}
\begin{definition} \label{9.29.04.3} 
The cluster ensemble related to seed ${\bf i}$  is 
the pair of 
positive spaces ${\cal A}_{|{\bf i}|} $ 
and ${\cal X}_{|{\bf i}|}$,  with common coordinate groupoid 
${\cal G}_{|{\bf i}|}$, 
related by a (monomial) morphism of positive spaces (\ref{9.28.04.3}). 
\end{definition}

 The algebra of regular functions on the ${\cal A}$-space is the same thing as 
the  algebra of universal Laurent polynomials ${\Bbb L}({\cal A})$. The Laurent phenomenon theorem \cite{FZ3} 
implies that the  cluster algebra of \cite{FZI} is a subalgebra of 
${\Bbb L}({\cal A})$. The algebra ${\Bbb L}({\cal A})$ 
 is bigger then the cluster algebra in most cases. 
It coincides with the upper cluster algebra introduced in 
\cite{BFZ}. 

\vskip 2mm
Alternatively, one can describe the above families of birational isomorphisms of seed tori 
by introducing {\it cluster ${\cal X}$- and ${\cal A}$-schemes}. By the very definition, a cluster ${\cal X}$-scheme is the scheme ${\cal X}^{*}$ related to the positive space ${\cal X}$, and similarly the cluster ${\cal A}$-scheme. 
Below we skip the superscript $*$ in the notation for cluster schemes.

 \vskip 2mm

Cluster transformations respect both the Poisson structures and 
the forms $\Omega$. Thus 
${\cal X}$ is a Poisson space, and there is a  $2$-form on the space 
${\cal A}$. (Precisely, ${\cal X}$ is understood as a functor 
from the coordinate groupoid to the appropriate category of Poisson tori.) 
In particular the manifold  ${\cal X}(\R_{>0})$ has a 
$\Gamma$-invariant Poisson structure. 
We show in Sections 3 and \ref{motivic} that the Poisson structure on the space ${\cal X}$ and the 
$2$-form 
on the space ${\cal A}$  are shadows of more sophisticated structures, namely  
a non-commutative $q$-deformation of the ${\cal X}$-space, 
and motivic avatars of 
the form $\Omega$.

\paragraph{The chiral and Langlands duality for seeds.} 
We define the  Langlands dual seed by 
$$
{\bf i}^{\vee}:= (I, I_0, \varepsilon^{\vee}_{ij}, d^{\vee}_i), \qquad 
\varepsilon^{\vee}_{ij}:= - \varepsilon_{ji}, \quad d_i^{\vee}:= d_i^{-1}
D, \qquad D:= {\rm l.c.m.}\{d_i\}\footnote{Here 
$D$ is the least common multiple of the set of positive integers $d_i$.}.
$$
This procedure is evidently involutive. 
Here is an alternative description. 

(i)  We define the {\it transposed} seed $
{\bf i}^t:= (I, I_0, \varepsilon^t_{ij}, d^t_i)$, where 
$\varepsilon^t_{ij}:= \varepsilon_{ji}$ and $d_i^t:= d_i^{-1}D.$

(ii) We define the {\it chiral dual}  seed 
$
{\bf i}^{o}:= (I, I_0, \varepsilon^o_{ij}, d^o_i)$, where 
$\varepsilon^o_{ij}:= - \varepsilon_{ij}$, and $d_i^o:= d_i$. 

Definitions (i)-(ii) are 
consistent with mutations. 
Combining them, we 
get the Langlands duality on  seeds. On the language of lattices. 
the Langlands duality amounts to replacing the bilinear form $[a,b]$ to the one 
$-[b,a]$, and changing the multipliers. 

Here is a natural realization of the Langlands dual seed. 
Let $\Lambda^{\vee}$ be the lattice dual to the lattice $\Lambda^{\circ}$. 
So, given a seed ${\bf i}= \Bigl(\Lambda, [\ast, \ast]_{\Lambda}, \{e_i\}, d_i\Bigr)$,  
there is an isomorphicm of lattices
\be \la{deltamap}
\delta_{\bf i}: \Lambda \lra \Lambda^{\vee}, \qquad e_i \lms e_i^\vee:= d_ie_i. 
\ee
Let us introduce a bilinear form on $\Lambda^\vee$ by setting
$$
[e_i^\vee, e_j^\vee]_{\Lambda^\vee}:= -[e_j, e_i]_{\Lambda}.
$$
\bl \la{8.4.09.1}
The map $\delta_{\bf i}$ provides an isomorphism of the Langlands dual seed ${\bf i}^\vee$ with the 
seed 
$$
\Bigl(\Lambda^\vee, [\ast, \ast]_{\Lambda^\vee}, \{e_i^\vee\}, d^\vee_i\Bigr).
$$ 
This isomorphism is compatible with mutations, i.e. the following diagram is commutative:
$$
\begin{array}{ccc}
\{e_i\}&\stackrel{\delta_{\bf i}}{\lra}&\{e_i^\vee\}\\
\mu_k\downarrow ~~ &&~~~\downarrow \mu^\vee_k \\
\{e'_i\}&\stackrel{\delta_{\bf i'}}{\lra}&\{(e'_i)^\vee\}
\end{array}
$$
\el

{\bf Proof}. The case $i=k$ is obvious. If $i \not = k$, we have 
$\mu_k(d_ie_i) = (d_ie_i) + d_i[\varepsilon_{ik}]_+d_k^{-1}(d_ke_k)$. 
So the Lemma follows from the formula 
\be \la{xc}
d_i\varepsilon_{ik}d_k^{-1} = - \varepsilon_{ki}. 
\ee

\subsection{An example: cluster ${\cal X}$-variety 
structure of the moduli space ${\cal M}_{0,n+3}$} \la{simex}

The moduli space  ${\cal M}_{0,n+3}$ parametrizes configurations 
of $n+3$ distinct points $(x_1, ..., x_{n+3})$ 
on ${\Bbb P}^1$ considered modulo the action of $PGL_2$. 

{\bf Example}. The cross-ratio of four points on ${\Bbb P}^1$, 
normalized by $r^+(\infty, -1, 0, z) = z$, provides an isomorphism 
$$
r^+: {\cal M}_{0,4} \stackrel{\sim}{\lra} {\Bbb P}^1 - \{0, -1, \infty\}, \qquad 
r^+(x_1, x_2, x_3, x_4) = \frac{(x_1-x_2)(x_3-x_4)}{(x_2-x_3)(x_1-x_4)}.
$$ 

The moduli space  ${\cal M}_{0,n+3}$  has a cluster 
${\cal X}$-variety atlas \cite{FG1}, which we  recall now. 
It is determined by a cyclic order of the points $(x_1, ..., x_{n+3})$. 
So although the symmetric group $S_{n+3}$ acts by automorphisms 
of ${\cal M}_{0,n+3}$, only its 
cyclic subgroup 
$\Z/(n+3)\Z$ will act by automorphisms of the cluster structure.

Let $P_{n+3}$ be a convex polygon with vertices $p_1, ..., p_{n+3}$. 
We assign the points $x_i$ to the vertices $p_i$, 
so that the order of points $x_i$ is compatible with  the 
clockwise  cyclic order of the vertices.
The cluster coordinate systems 
are parametrized by the set ${\cal T}_{n+3}$ of complete triangulations 
of the polygon $P_{n+3}$. 
Given such a triangulation $T$, the coordinates are assigned to the diagonals of $T$.
The coordinate $X^T_E$ corresponding to a triangulation $T$ and its diagonal
 $E$ is defined as follows. 
There is a unique rectangle formed by the sides and diagonals of the triangulation, 
with the diagonal given by  $E$. 
Its vertices provide a cyclic configuration of four points on ${\Bbb P}^1$. 
We order them starting from a vertex of the $E$, 
getting a configuration of four points 
$(x_1, x_2, x_3, x_4)$ on ${\Bbb P}^1$.  
Then we set 
$$
X^T_E:= r^+(x_1, x_2, x_3, x_4).
$$
\begin{figure}[ht]
\begin{picture}(216,72)
\centerline{\epsfbox{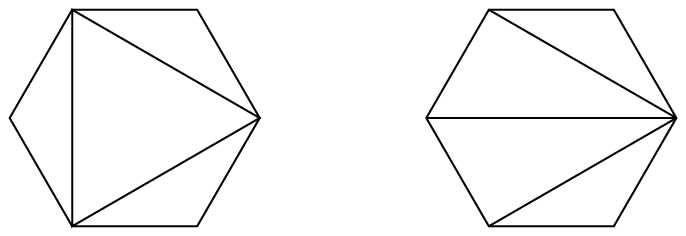}}
\put(-238,34){$x_1$}
\put(-358,34){$x_1$}
\put(-145,34){$x_4$}
\put(-265,34){$x_4$}
\put(-220,-2){$x_6$}
\put(-340,-2){$x_6$}
\put(-162,-2){$x_5$}
\put(-282,-2){$x_5$}
\put(-220,70){$x_2$}
\put(-340,70){$x_2$}
\put(-162,70){$x_3$}
\put(-282,70){$x_3$}
\put(-316,34){$E$}
\end{picture}
\caption{The two triangulations of the hexagon are related by the flip at the edge $E$.}
\label{clusp11}
\end{figure}
%\begin{figure}[ht]
%\centerline{\epsfbox{fliphex.eps}}
%\caption{The two triangulations of the hexagon are related by a flip.}
%\label{clusp11}
%\end{figure}
There are exactly two ways to order the points as above, which differ by a 
cyclic shift by two. Since the cyclic shift by one changes the cross-ratio to its inverse, 
the rational function $X^T_E$ is well defined. For example the 
diagonal $E$ on Fig \ref{clusp11} provides the function $r^+(x_2, x_4, x_6, x_1)$.

We define the cluster seed assigned to a triangulation $T$ as follows. The lattice 
$\Lambda$ is the free abelian group generated by the diagonals of the triangulation, with  
a basis is given by the diagonals. 
The bilinear form is given by the adjacency matrix. 
Namely, two diagonals $E$ and $F$ of the triangulation are called adjacent 
if they share a vertex, and 
there are no diagonals of the triangulation between them. 
We set $\varepsilon_{EF} =0$ if $E$ and $F$ are not adjacent. 
If they are, $\varepsilon_{EF} =1$ if $E$ is before $F$ 
according to the clockwise orientation of the diagonals at the vertex $v$ shared by 
$E$ and $F$, and $\varepsilon_{EF} = - 1$ otherwise. 

\begin{figure}[ht]
\centerline{\epsfbox{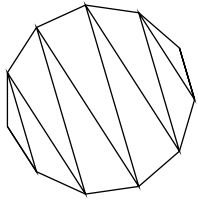}}
\caption{A zig-zag triangulations of the $11$-gon.}
\label{clusp13}
\end{figure}

\vskip 2mm
{\bf Example}. Let us consider a  triangulation of $P_{n+3}$ which has the following 
property: every triangle of the triangulation contains at least one side of the 
polygon. Then it provides a seed of type $A_n$. 
For example, a zig-zag triangulation, see Fig \ref{clusp13}, has this property.  

One shows that a mutation at a diagonal $E$ corresponds to 
the flip of the diagonal, see Fig \ref{clusp11}. 
This means that formula (\ref{5.11.03.6}) describes 
the adjacency matrix of the flipped triangulation. 
This way we get a cluster ${\cal X}$-variety 
atlas. One easily sees that the zig-zag triangulation provides 
a quiver of type $A_n$.

The cyclic order of the points 
$(x_1, ..., x_{n+3})$ 
provides a connected component ${\cal M}^0_{0,n+3}(\R)$ 
in ${\cal M}_{0,n+3}(\R)$, parametrising configurations of points on 
${\Bbb P}^1(\R)$ whose cyclic order is compatible 
with an orientation of ${\Bbb P}^1(\R)$. 
The space of positive points of the cluster ${\cal X}$-variety 
defined above coincides with ${\cal M}^0_{0,n+3}(\R)$.

\vskip 2mm
{\bf Remark.} We show in \cite{FG4} that the Knudsen-Deligne-Mumford moduli space 
$\overline {\cal M}_{0,n+3}$ can be recovered in a natural way as a cluster compactification 
of the cluster ${\cal X}$-variety of type $A_n$.

\subsection{Cluster ensemble structures}
Below we outline the structures related to a cluster ensemble. 
\vskip 2mm

A cluster ensemble  gives rise to the following data:
\vskip 2mm
i) A  pair of real manifolds 
$
{\cal A}(\R_{>0})$ and ${\cal X}(\R_{>0})$, provided 
by the positive structures on ${\cal A} $ and  ${\cal X}$, and a map $p:{\cal A}(\R_{>0}) \to 
{\cal X}(\R_{>0})$.  
For a given seed ${\mathbf i}$,  the functions 
$\log|A_i|$ (resp. $\log|X_i|$)  
provide  diffeomorphisms 
$$
\alpha_{\mathbf i}:  {\cal A}(\R_{>0}) 
 \stackrel{\sim}{\lra} \R^{I}; \qquad \beta_{\mathbf i}: 
{\cal X}(\R_{>0})  \stackrel{\sim}{\lra} \R^{I}.
$$
Similarly there are sets  of points of  ${\cal A}$ and  ${\cal X}$ 
with values in the {tropical semifields} $\Z^t$, $\Q^t$,  $\R^t$.

\vskip 2mm
ii) A modular group  ${\Gamma}$  acts by automorphisms of the whole structure. 

\vskip 2mm
iii) A $\Gamma$--invariant Poisson structure $\{*, *\}$ 
on ${\cal X}$. In any 
${\cal X}$--coordinate system $\{X_{i}\}$   it  is the  
quadratic Poisson structure given by (\ref{4.30.03.2}).  
\vskip 2mm
iv) A $\Gamma$--invariant $2$--form $\Omega$ on ${\cal A}$, which in any ${\cal A}$--coordinate system $\{A_{i}\}$  is given by  
(\ref{4.30.03.2WP}). 
It can be viewed as a presymplectic structure on ${\cal A}$. 

\vskip 2mm
v) A pair of split algebraic tori of the same dimension, 
$H_{\cal X}$ and $H_{\cal A}$. 
The torus $H_{\cal A}$ acts freely on ${\cal A}$.  
The orbits are the fibers of the map $p$, 
and 
 the leaves of the null foliation for the 
$2$-form $\Omega$. 
Thus 
$
{\cal U}:= p({\cal A})
$ 
 is a positive symplectic space. 
Dually, there is a canonical projection 
$
\theta: {\cal X} \to H_{\cal X}. 
$ Its fibers are 
the symplectic leaves of the  Poisson structure. Moreover 
${\cal U} = \theta^{-1}(e)$, where $e$ is 
the unit of $H_{\cal X}$. 
So the natural embedding
$i: {\cal U} \hra {\cal X}$ is a Poisson map. 
\vskip 2mm
vi) A 
{\it quantum space}   ${\cal X}_q$. It is a non-commutative 
 $q$-deformation 
of the positive space ${\cal X}$, equipped with an involutive 
antiautomorphism $\ast$,  understood 
as a functor
$$
\psi_q: \widehat {\cal G} \lra {\rm QPos}^*
$$
where ${\rm QPos}^*$ is the category of quantum tori with involutive
antiautomorphism $*$.   Precisely, the category ${\rm QPos}^*$ is 
the opposite category to  
the category whose objects are quantum tori algebras, and 
morphisms are positive rational 
$*$-maps.\footnote{Quantum torus algebra satisfies Ore's condition, 
so its non-commutative fraction field and hence rational functions as its elements are defined.  
The reader may skip positivity from the above definition.}  
For a   seed
${\mathbf i}$,   the corresponding quantum torus algebra 
$\psi_q({\mathbf i})$ is  the algebra ${\bf T}_{ {\mathbf i}}$ 
generated by the elements $X_i$, $i \in I$, subject to the relations
\begin{equation} \label{4.28.03.11x}
q^{-\widehat \varepsilon_{ij }} X_i X_j = 
q^{-\widehat \varepsilon_{ji }} X_j X_i, \qquad \ast X_i = X_i, \quad \ast q = q^{-1}. 
\end{equation}
We denote by ${\Bbb T}_{ {\mathbf i}}$ its non-commutative field of fractions. 
Given a mutation ${\mathbf i} \to {\mathbf i}'$, 
the birational map $\psi^{{\mathbf i}, {\mathbf i}'}_q: {\Bbb T}_{ {\mathbf i}'} \to 
{\Bbb T}_{\mathbf i} $ is a $q$-deformation 
of the mutation map (\ref{5.11.03.1x}) from the definition  
of the positive space ${\cal X}$. 
It is given by 
the conjugation by the quantum dilogarithm.

There is a canonical projection  
$
\theta_q: {\cal X}_q \lra H_{\cal X}.
$ 
 The inverse images of the characters of $H_{\cal X}$ are ``quantum Casimirs'': 
they are in the center 
of  ${\cal X}_q$  and generate it for generic $q$.   

\vskip 2mm
{\it The quantum space ${\cal X}_q$ at roots of unity.} Now suppose that $q^{DN}=1$, where 
$D$ is the least common multiple of $d_i$'s.  Then, under certain assumption on the exchange function, 
see Theorem \ref{QFROB},  there is a 
{\it quantum Frobenious map} of quantum spaces 
$$
{\Bbb F}_N: {\cal X}_q \lra {\cal X},
$$
which in any cluster coordinate system acts on the cluster coordinates 
$Y_i$ on 
${\cal X}$ as 
 ${\Bbb F}_N^*Y_i = X_i^N$. Here  $X_i$ the coordinates on ${\cal X}_q$. 
Notice that ${\Bbb F}_N^*$ is a ring homomorphism since $q^N=1$. 
So there is a   diagram:
$$
\begin{array}{ccc}
{\cal X}_q& \stackrel{\theta_q}{\lra}&H_{\cal X}\\
&&\\
{\Bbb F}_N \downarrow ~~&&\\
&&\\
{\cal X}&&
\end{array}
$$
The center 
of the algebra of regular functions ${\Bbb L}({\cal X}_q)$ is 
generated by the inverse images of the functions on ${\cal X}$ and 
$H_{\cal X}$. 
 
\vskip 2mm
vii) {\it Motivic data.}\footnote{The reader is advised to look at Section 
\ref{motivic}.2 for the 
background related to the motivic data.} There are two levels of understanding:

 a) {\it $K_2$--avatar of $\Omega$}. It is given by 
a $\Gamma$--invariant 
class 
${W} \in K_2({\cal A})^{\Gamma}$.

b) {\it Motivic dilogarithm class.} For any  
seed ${\mathbf i}$ 
the class $W$ can be lifted to an element 
$$
{W}_{\bf i} = \sum_{i, j}\widetilde \varepsilon_{ij }
A_{i}\wedge A_{j} \in {\bigwedge}^2\Q({\cal A})^*. 
$$
It has  the zero tame symbol at every divisor. It is $H_{\cal A}$-invariant, 
and thus is a lift of an element of $\Lambda^2\Q({\cal U})^*$ by the map $p$. However elements 
${W_{\bf i}}$ are not equivariant under the action of the cluster modular groupoid, and in particular are not 
$\Gamma$-invariant. Their behavior under the action of the cluster modular groupoid 
is described by a class, called {\it motivic dilogarithm class of cluster ensemble}: 
$$
{\Bbb W} \in H^2_{\Gamma}({\cal U}, \Q(2)_{\cal M})
$$
in the weight two $\Gamma$--equivariant motivic cohomology 
of the scheme ${\cal U}$. We define 
the weight two motivic cohomology via the dilogarithm complex,  
also known as the Bloch-Suslin complex.

\vskip 2mm
The simplest way to see the dilogarithm in our story is the following. 
Mutations act by Poisson automorphisms of ${\cal X}(\R_{>0})$. 
The generating function   
describing a mutation ${\mathbf i} \to {\mathbf i}'$ at a vertex $k$ 
is Roger's version of the 
dilogarithm, applied to the  coordinate function $X_k = e^{x_k}$.

\vskip 2mm
In \cite{FGII}, which is the second part of this paper, we 
introduce one more ingredient of the data 
which, unlike everything else above, is of analytic nature: 

 viii) {\it $\ast$-representations}. a) A series of $\ast$-representations by unbounded operators 
in Hilbert spaces of the 
{\it modular double} 
\be \la{md}
{\cal X}_{q, |{\bf i}|} \times {\cal X}_{q^{\vee}, |{\bf i^{\vee}}|},  ~~~~~~
q = e^{i\pi h}, ~q^{\vee}= e^{i\pi /h}, \quad h \in \R_{\geq 0}, 
\ee
 where  
${\bf i^{\vee}}$ is the { Langlands dual} seed. 
  It is constructed explicitly using 
the quantum dilogarithm. 

b) The  modular double of the quantum space 
${\cal X}_{q, |{\bf i}|} \times {\cal X}_{q, {|{\bf i^o}|}}$, where 
${{\bf i^o}}$ defines  the {chiral dual}  seed, 
has {\it canonical} unitary projective representation in the Hilbert space 
$L_2({\cal A}(\R_{>0}))$. 

%There is a 
%$\Gamma_{|{\bf i}|}$--equivariant 
%gerb 
%on the space ${\cal A}_{|{\bf i}|}$. 
%The objects of this gerb are certain holomorphic line bundles with connection whose curvature 
%is the form $\Omega_{|{\bf i}|}$. 
%The equivalence class of this gerb is described by the class ${\Bbb W}_{|{\bf i}|}$.
\vskip 2mm

ix) Finally, in \cite{FG2}, which continues this paper, we define the third 
ingredient of the cluster ensemble, the {\it symplectic double}. It is the quiasiclassical counterpart of the 
canonical representation  from the part b) of viii). We define its non-commutative $q$-deformation. 
We construct a series of infinite dimensional unitary projective 
representations of the special cluster modular groupoid, 
parametrised by unitary characters of the torus $H_{\cal X}(\R_{>0})$. We show that they intertwine  
the $\ast$-representations of (\ref{md}) in certain Schwartz type spaces, 
and using this prove all claims from  viii).

\vskip 2mm
The described cluster 
quantization is a rather general quantization scheme, which we hope has many 
applications. 
\vskip 2mm

Remarkably the
quantization is governed by the motivic avatar of the Weil-Petersson form on
the ${\cal A}$-space.  
This and the part b) of vii) show that although the ${\cal X}$-space seems to
be the primary part of a cluster ensemble - in our basic example it 
gives the Teichm\"uller space, while the
${\cal A}$-space gives only its
decorated 
unipotent part -  one has to
study the  ${\cal X}$ and ${\cal A}$ spaces in a package.

\subsection{Our basic example} Let $G$ be a split semi-simple simply-connected 
algebraic group over $\Q$. Denote by $G'$ the quotient of $G$ modulo the center. 
Let $S$ be a hyperbolic surface with non-empty boundary and $m$ distinguished points on the boundary. 
We defined in \cite{FG1} the a pair of moduli spaces (${\cal X}_{G',S}, {\cal A}_{G,S}$), 
and proved that for $G=SL_m$ it 
gives rise to a cluster ensemble, leaving the case of  general $G$ 
to the sequel of that paper. We proved that,  regardless 
of the cluster ensemble structure, this 
pair of 
moduli spaces  for general $G$ 
have all  described above classical structures. 
Here is a more detailed account. 
The references are made to chapters of \cite{FG1}. 
The example discussed in Section \ref{simex} is the special case when 
$S$ is a disc with $n+3$ marked points on the boundary, and $G = PGL_2$. 

\vskip 3mm
The pair of positive spaces $({\cal X}, {\cal A})$  
is provided by the pair of positive 
stacks  ${\cal X}_{G',S}$ and ${\cal A}_{G,S}$.  

i) The corresponding pair of 
positive real spaces  is the higher 
Teichm{\"u}ller space
 ${\cal X}_{G', S}(\R_{>0})$ and its decorated version
 ${\cal A}_{G, S}(\R_{>0})$. 
The 
 ${\Bbb A}^t$-points of ${\cal A}_{SL_2, S}$ give Thurston's laminations  
(Chapter 12). 
The space of positive real points of 
${\cal X}_{PGL_2, S}$ is a version of the classical Teichm{\"u}ller space 
 on $S$, 
and the one of ${\cal A}_{SL_2, S}$ is Penner's \cite{P1} decorated Teichm{\"u}ller space  (Chapter 11).  

\vskip 3mm
ii) There is a modular group $\Gamma_{G,S}$ provided by the cluster ensemble 
structure of the pair $({\cal X}_{G',S}, {\cal A}_{G,S})$. 
The positive spaces ${\cal X}_{G',S}$ and ${\cal A}_{G,S}$ are 
$\Gamma_{G,S}$-equivariant positive
spaces. The group $\Gamma_{G,S}$ 
contains as a subgroup the modular group 
 $\Gamma_S$ of the surface $S$. 
If $G = SL_2$, these two groups coincide. Thus the cluster modular group  
is a generalization  of the classical modular group. 
Otherwise  $\Gamma_{G,S}$ is bigger then 
$\Gamma_S$. For example, the cluster modular group $\Gamma_{G,S}$ 
where $G$ is of type $G_2$ and $S$ is a disc with three points on the boundary 
was calculated in \cite{FG3}: it is (an infinite quotient of)
 the braid group of type $G_2$, while its classical counterpart is $\Z/3\Z$.

The quotient ${\cal M}:= {\cal X}(\R_{>0})/
\Gamma$ 
is an analog of the moduli space of complex structures on a surface. 
We conjecture that the space 
${\cal M}_{G',S}$ is related to the $W$-algebra for the group $G$ 
just the same way the classical moduli space (when $G=SL_2$)
is related to the Virasoro algebra. 
We believe that ${\cal M}_{G',S}$ is the moduli space 
of certain objects, $W$--structures, but can not define them.

\vskip 3mm
iii) There is a canonical projection from the moduli space ${\cal X}_{G',S}$ 
to the moduli space of $G'$-local systems on $S$. 
The Poisson structure on ${\cal X}_{G',S}$ 
is the inverse image of 
the standard Poisson structure on the latter by this map. 
%(\cite{FR}). 
%See \cite{FG} for a description of 
%its restriction to the Teichm{\"u}ller space 
%for $G' = PGL_2$. 

\vskip 3mm
iv) The form $\Omega$ on the space ${\cal A}_{G,S}$ was defined in  Chapter 15. 
For $G = SL_2$ its restriction to the 
decorated Teichm{\"u}ller space ${\cal A}_{SL_2, S}(\R_{>0})$
is the Weil-Petersson form studied by Penner \cite{P1}. 
\vskip 3mm
v) {\it The tori $H_{\cal A}$ and $H_{{\cal X}}$}. 
Let $H$ be the Cartan group of  $G$, and  $H'$ the 
Cartan group of $G'$. The canonical projection $\theta: 
{\cal X} \lra H_{{\cal X}}$
and the action of the torus $H_{{\cal A}}$ on the ${\cal A}$-space 
generalize similar structures defined in Chapter 2  for a
hyperbolic surface $S$:  the canonical projection  
$$
{\cal X}_{G', S} \lra H'^{\{\mbox{punctures of $S$}\}}
$$ 
and the action of 
$H^{\{\mbox{punctures of $S$}\}}$ on the moduli space ${\cal A}_{G, S}$. 
\vskip 3mm
vi) The results
of this paper plus the cluster ensemble structure  of the pair 
$({\cal X}_{PGL_m,S}, {\cal A}_{SL_m,S})$ (Chapter 10) provide  a  
quantum space  ${\cal X}^q_{PGL_m,S}$. 
For $m=2$ it is equivalent to the one in \cite{FCh}.  

\vskip 3mm
vii) The motivic data for the pair 
$({\cal X}_{G',S}, {\cal A}_{G,S})$ was defined in Chapter 15.  It was 
previously missing even 
for the classical Teichm{\"u}ller space. 
In the case $G = SL_m$ an explicit cocycle representing the class ${\Bbb W}$ 
is obtained from the explicit cocycle representing 
the second motivic Chern class 
of the simplicial classifying space $BSL_m$ defined in \cite{G2}. 
The investigation of this cocycle for ${\Bbb W}$ led us to discovery 
of the whole picture. 
\vskip 3mm
viii) Replacing the Dynkin diagram of the group $G$ by its Langlands dual we get the
  Langlands dual cluster ensemble. Changing the orientation on $S$ we
 get the chiral dual cluster ensemble. 

\vskip 3mm
The classical Teichm{\"u}ller space was quantized, independently, 
 in \cite{K} and in \cite{FCh}: 
the Poisson manifold 
${\cal X}_{PGL_2, S}(\R_{>0})$ was quantized in \cite{FCh}, and  
its symplectic leaf ${\cal U}_{PGL_2, S}(\R_{>0})$ 
in \cite{K}. 
\vskip 3mm
The  principal embedding $SL_2 \hra G$,   
defined up to a conjugation,  
leads to  natural embeddings
$$
{\cal X}_{PGL_2, S} \hra {\cal X}_{G, S}, \quad {\cal A}_{SL_2, S} \hra 
{\cal A}_{G, S}
$$
and their counterparts for the Teichm{\"u}ller, lamination and moduli spaces. 
However  since the cluster modular 
group $\Gamma_{G,S}$ is bigger then the modular group 
$\Gamma_S$, it is hard to expect 
natural $\Gamma_{G,S}$-equivariant 
projections like  ${\cal X}_{G, S} \lra {\cal X}_{PGL_2, S}$. 
We do not know whether  $W$-structures on $S$ can be defined 
as a complex structure plus some extra data on $S$.
 
\vskip 3mm
The elements $p^*(X_i)$ of the cluster algebra were considered 
by Gekhtman, Shapiro and Vainshtein \cite{GSV1} who studied various  
Poisson 
structures on a cluster algebra. 
The form $\Omega$ and the connection between Penner's   
decorated Teichm\"uller spaces to cluster algebras were 
independently discovered in 
\cite{GSV2}. The relation of the form $\Omega$ to the 
Poisson structures is discussed there. 

\vskip 3mm
After the first version of this paper appeared 
in the ArXiv (math/0311245), Berenstein and Zelevinsky released a 
paper \cite{BZq}, where they 
defined and studied 
$q$-deformations of cluster algebras. 
In general 
there is a family of such 
$q$-deformations, matching the Poisson structures 
on cluster algebras defined in \cite{GSV1}. 
The cluster modular group does not preserve 
individual $q$-deformations.  
However if ${\rm det }\varepsilon_{ij} \not = 0$, 
the $q$-deformation of cluster algebra
is unique, and thus $\Gamma$-invariant.

\subsection{The structure of the paper} Cluster ensembles are 
studied in Section 2. 
We discuss the cluster nature of the Teichm\"uller theory on 
a punctured torus, as well as  the cluster structure of the 
pair of universal Teichm\"uller spaces. In the latter case 
the set $I$ is     
the set of edges of the Farey triangulation of the hyperbolic plane, and 
the modular group is the Thompson group.

In Section 3 we define the 
 non-commutative  ${\cal X}$-space and establish its properties. 

In Section 4 we present our duality conjectures. 

In Section 5 we furnish some evidence for the duality conjectures in the finite type case:  
We define, in the finite type case, a canonical pairing 
between the tropical points of dual cluster varieties 
and one of the two canonical maps. 
Our definitions do not depend on the Classification Theorem, and do not use root systems etc. 
As a byproduct, we obtain a canonical 
decomposition of the 
space of real tropical points of a finite type cluster ${\cal X}$-variety. 
It is   dual to 
 the  generalised associahedra defined in \cite{FZ}.

%In Section 5 we define the tropical boundary of a positive space. 
%We discuss it for the 
%Teichm\"uller space on the punctured torus.
%We define special completions of cluster ${\cal X}$-varieties. 
%In Section 6 we introduce operations of amalgamation and folding. 
%They are used extensively in \cite{FG2} and in the investigation of cluster structure of moduli spaces 
%${\cal X}_{G, S}$ and ${\cal A}_{G, S}$. 
%The appendix contains some comparison results. 
%In Section 6 we investigate cluster varieties of the finite type $A_n$. 
%We prove the duality conjectures for them. 
%We prove Theorems \ref{compact2} and  \ref{compact1}. 

Section \ref{motivic} we introduce motivic structures related to a cluster ensemble. 
They are defined using 
the  dilogarithmic  motivic complex, and 
play a key role in our understanding of cluster ensembles.

\vskip 3mm
In \cite{FGII}  we started a program of 
quantization of cluster ensembles
using the quantum dilogarithm 
intertwiners. It is a quantum version of the motivic data 
from Section \ref{motivic}.   It was completed in \cite{FG2}.

\paragraph{Acknowledgments.} 
This work was developed  when 
the first author visited Brown University. He would like to thank Brown University 
for hospitality and support. He was supported by the grants CRDF 2622; 2660. 
The second author was 
supported by the  NSF grants DMS-0099390,   DMS-0400449 and  DMS-0653721. 
The first draft of this paper was completed when the second author enjoyed the 
hospitality of IHES (Bures sur Yvette)  
and MPI (Bonn)
during the Summer and Fall of 2003. 
He is grateful to NSF, IHES and MPI for the support. 

We are grateful to Sergey Fomin for pointing out 
some errors in the first version of
the paper, and to Benjamin Enriquez who read the text carefully, made a lot of 
suggestions regarding the exposition,  
and corrected several errors. We thank the referees for the comments.

\section{Cluster ensembles and their properties}

\subsection{Cluster ensembles revisited}

\paragraph{Cluster transformations of cluster seed tori.}  
A seed mutation  $\mu_k : {\mathbf i} \to {\mathbf i}'$ provides 
birational isomorphisms   
\be \la{ax}
\mu^x_{k}: {\cal X}_{{\mathbf i}} \lra {\cal X}_{{\mathbf i}'}\quad \mbox{ and } \quad 
\mu^a_{k}: {\cal A}_{{\mathbf i}} \lra {\cal A}_{{\mathbf i}'} 
\ee
acting on cluster coordinates by formulas (\ref{5.11.03.1x})
and (\ref{5.11.03.1a}), respectively. Seed isomorphisms  provide  
isomorphism of tori, see (\ref{6.4.07.1}).  
So a seed cluster transformation  ${\bf c}: {\mathbf i} \to 
{\mathbf i}'$  
 gives rise to  birational isomorphisms 
\begin{equation} \label{5.11.03.1rtsd}
{\bf c}^a: {\cal A}_{{\mathbf i}} \lra {\cal A}_{{\mathbf i}'}, \qquad 
{\bf c}^x: {\cal X}_{{\mathbf i}} \lra {\cal X}_{\mathbf i'}
\end{equation}

Recall the coordinate groupoid ${\cal G}$. 
The following lemma results from the very definition. 
\begin{lemma} \label{7.2.03.6a} 
There are well defined functors 
$$
\psi_{\cal A}: {\cal G} \lra {\rm Pos}, \quad 
\psi_{\cal A}({\mathbf i}) :=
{\cal A}_{\mathbf i}, \quad \psi_{\cal A}(\mu_k):= \mu^a_{k},
$$
$$
\psi_{\cal X}: {\cal G} \lra {\rm Pos}, 
\quad \psi_{\cal X}({\mathbf i}) :=
{\cal X}_{\mathbf i}, \quad \psi_{\cal X}(\mu_k):= \mu^x_{k}.
$$
\end{lemma}

\vskip 2mm
Let  ${\cal A}$ and 
${\cal X}$ be the positive spaces defined by the functors 
from Lemma \ref{7.2.03.6a}. 
These spaces are related as follows. Given a seed ${\mathbf i}$, the map $p$ 
looks in coordinates as follows: 
\begin{equation}\label{7.7.04.13df}
p_{\mathbf i}: {\cal A}_{\mathbf i} \lra {\cal X}_{\mathbf i}, \qquad 
p_{\mathbf i}^*X_i := \prod_{j \in I}A_j^{\varepsilon_{ij}}.
\end{equation}
\begin{proposition} \label{7.7.04.12} 
The maps of the seed tori (\ref{7.7.04.13df}) give rise to 
a  map of positive spaces $p: {\cal A} \lra {\cal X}$.
\end{proposition}

We will give a proof after a discussion of  decomposition of mutations.

\paragraph{Decomposition of mutations.} The seed tori 
${\cal A}_{\bf i}$ and ${\cal A}_{\bf i'}$  (respectively ${\cal X}_{\bf i}$ and ${\cal X}_{\bf i'}$) 
are canonically identified with the torus ${\cal A}_{\Lambda}$
(respectively ${\cal X}_{\Lambda}$). Therefore there are tautological isomorphisms
$$
\mu'_k:{\cal A}_{\bf i} \stackrel{\sim}{\lra} {\cal A}_{\bf i'}, \qquad 
\mu'_k:{\cal X}_{\bf i} \stackrel{\sim}{\lra} {\cal X}_{\bf i'}. 
$$
These isomorphisms, however, do not respect the cluster coordinates on these tori. 
Therefore there are two ways to write the mutation transformations:

\vskip 1mm
(i) Using the cluster coordinates assigned to the seeds ${\bf i}$ and ${\bf i'}$, or

(ii) Using the cluster coordinates assigned to the seeds ${\bf i}$ only. 
\vskip 1mm
\noindent
Equivalently, in the approach (ii) we present 
mutation birational isomorphisms (\ref{ax}) as  compositions
$$
\mu^a_k = \mu'_k\circ \mu^{\sharp}_k, \qquad  \mu^{\sharp}_k: {\cal A}_{\bf i}\lra {\cal A}_{\bf i}, 
\quad \mu'_k: {\cal A}_{\bf i}\stackrel{\sim}{\lra} {\cal A}'_{\bf i},
$$
$$
\mu^x_k = \mu'_k\circ \mu^{\sharp}_k, \qquad  \mu^{\sharp}_k: {\cal X}_{\bf i}\lra {\cal X}_{\bf i}, 
\quad \mu'_k: {\cal X}_{\bf i}\stackrel{\sim}{\lra} {\cal X}'_{\bf i},
$$
and then look at the birational isomorphisms $\mu_k^{\sharp}$ only.

We usually use the approach (i). In particular formulas (\ref{5.11.03.1x})
and (\ref{5.11.03.1a}) were written this way. However the approach (ii) 
leads to simpler formulas, which are 
  easier to deal with, especially for the  ${\cal X}$-space. 
What is more important, the conceptual meaning of the map 
$\mu_k^{\sharp}$ in the ${\cal X}$-case 
becomes crystal clear when we go to the $q$-deformed spaces, see Section 4. 

Below we work out these formulas, i.e. compute mutation birational automorphisms 
$\mu_k^{\sharp}$ in the cluster coordinates 
assigned to the seed ${\bf i}$. 

 It is handy to employ the following notation:
$$
{\Bbb A}_k^+:= \prod_{i|\varepsilon_{ki}>0}A_i^{\varepsilon_{ki}}, \qquad
{\Bbb A}_k^-:= \prod_{i|\varepsilon_{ki}<0}A_i^{-\varepsilon_{ki}}.
$$
Then 
\be \la{qa}
\frac{{\Bbb A}_k^+}{{\Bbb A}_k^-} = \prod_jA_j^{\varepsilon_{kj}} = p^*X_k.
\ee
\bp 
Given a  seed ${\bf i}$, the birational 
 automorphism $\mu_k^{\sharp}$ of the seed 
${\cal A}$-torus  (respectively ${\cal X}$-torus) 
acts on the cluster coordinates $\{A_i\}$ (respectively $\{X_i\}$) related to the seed ${\bf i}$ as follows: 
\begin{equation} \label{f3**}
A_i \lms A^{\sharp}_{i} :=   A_i(1+ p^*X_k)^{-\delta_{ik}}  =
\left\{\begin{array}{lll} A_i& \mbox{ if } 
& i\not =k, \\ A_k(1+ {\Bbb A}_k^+/{\Bbb A}_k^-)^{-1} & \mbox{ if } &  i= k. \\
\end{array} \right.
\end{equation}
\begin{equation} \label{f3**as}
X_i \lms X^{\sharp}_{i} := X_i(1+ {X}_k)^{-\varepsilon_{ik}} 
\end{equation}
\ep

{\bf Proof}. 
Let $\{A_i'\}$ be the cluster coordinates in the function field of ${\cal A}_{\Lambda}$ 
assigned to the mutated seed ${\bf i'}$. They are related to the cluster coordinates $\{A_i'\}$ 
assigned to the seed ${\bf i}$ as follows:  
\begin{equation} \label{11.18.06.10sdf}
A'_{i} \lms \left\{\begin{array}{lll} 
A_i& \mbox{ if } & i\not =k, \\
    {\Bbb A}_k^-/A_k
 & \mbox{ if } &  i= k. \\
\end{array} \right.
\end{equation}
These formulas describe the action of the tautological mutation isomorphism $\mu_k'$ 
on the cluster coordinates. They reflect the action of the seed mutation 
on the quasidual basis $\{f_i\}$, see (\ref{qdbas}). 

Then the transformation  $\mu'_k\circ \mu^{\sharp}_k$ acts on the coordinates $A_i'$ as follows: 
$$
(\mu^{\sharp}_k)^*\circ (\mu'_k)^*: A'_{k} \lms {\Bbb A}_k^-/A_k= 
{\Bbb A}_k^-(1+ \frac{{\Bbb A}_k^+}{{\Bbb A}_k^-})A_k^{-1}
 = \frac{{\Bbb A}_k^+ + {\Bbb A}_k^-}{A_k}.
$$
This coincides with the action of the mutation $\mu_k$ on the coordinate $A_k'$. This proves the ${\cal A}$-part of the proposition.

Similarly,  the tautological mutation isomorphism $\mu_k'$ acts on the coordinates by 
\begin{equation} \label{11.18.06.10hr}
X'_{i} \lms \left\{\begin{array}{lll} X_k^{-1}& \mbox{ if } & i=k, \\
    X_i(X_k)^{[\varepsilon_{ik}]_+} & \mbox{ if } &  i\neq k. \\
\end{array} \right.
\end{equation}
This reflects the action of the seed mutation 
on the basis $\{e_i\}$. 

Therefore the transformation  $\mu'_k\circ \mu^{\sharp}_k$ acts on the coordinates $X_i'$ as follows: 
$$
X_k' \lms X_k^{-1} \lms X_k^{-1},
$$ 
and, if $i \not = k$, 
$$
X_i' \lms X_i(X_k)^{[\varepsilon_{ik}]_+}  \lms X_i(1+X_k)^{-\varepsilon_{ik}} 
(X_k)^{[\varepsilon_{ik}]_+} = X_i(1+X_k^{-{\rm sgn}(\varepsilon_{ik})})^{-\varepsilon_{ik}}.
$$
This coincides with the action of the mutation $\mu_k$ on the coordinate $X'_i$. 
The proposition is proved.

\paragraph{Proof of Proposition \ref{7.7.04.12}.} It is 
equivalent to the following statement\footnote{As pointed out by a referee, relation between 
$p^*X_i$ and $p^*X_i'$ is
 equivalent to  Lemma 1.2 in \cite{GSV1}.}: 
For each mutation $\mu_k$ of the seed ${\mathbf i}$ there is a commutative diagram
\begin{equation}\label{7.7.04.13}
\begin{array}{ccc}
{\cal A}_{\Lambda}& \stackrel{\mu^\sharp_{k}}{\lra} &{\cal A}_{\Lambda}\\
%&&\\
p\downarrow && \downarrow p\\
%&&\\
{\cal X}_{\Lambda} &\stackrel{\mu^\sharp_{k}}{\lra}&{\cal X}_{\Lambda}.
\end{array}
\end{equation}
Let us prove this statement. 
Going up and to the left we get
$$
X_i\lms \prod_{j\in I}A_j^{\varepsilon_{ij}} \lms 
\prod_{j\in I}A_j^{\varepsilon_{ij}} \cdot 
(A_k^{\sharp}/A_k)^{\varepsilon_{ik}}. 
$$ 
 Going to the left and up we get the same:
$$
X_i \lms X_i(1+X_k)^{-\varepsilon_{ik}} \lms 
\prod_{j\in I}A_j^{\varepsilon_{ij}} 
(1+p^*X_k)^{-\varepsilon_{ik}}.
$$
\vskip2mm

\begin{definition} \label{8.15.04.3}
The space ${\cal U}$ is the image of the space 
${\cal A}$ under the map $p$.
\end{definition}

We leave to the reader to check that the space 
${\cal U}$ is indeed a positive space. 

\begin{corollary} \label{8.24.05.1x} Assume 
that ${\rm det}~\varepsilon_{ij} \not = 0$. Let ${\bf c}: {\bf i} \to {\bf i}$ 
be a seed cluster transformation. It gives rise to cluster 
transformations ${\bf c}^a$ and ${\bf c}^x$ of the ${\cal A}$ and ${\cal X}$ spaces. 
Then 
${\bf c}^a = {\rm Id}$ implies ${\bf c}^x = {\rm Id}$.
\end{corollary}

\vskip 2mm
{\bf Proof}. Assume that ${\rm det}~\varepsilon_{ij}\not = 0$. 
Then the map of algebras 
$%\begin{equation} \label{8.26.05.2x}
p^*: \Z[{\cal X}_{\mathbf i}] \lra \Z[{\cal A}_{\mathbf i}] 
$ %\end{equation}
is an injection, and commutes with the cluster transformations thanks to Proposition \ref{7.7.04.12}.  This implies the claim.

\paragraph{A Poisson structure on the ${\cal X}$-space}  
\begin{lemma} 
Cluster transformations preserve the  Poisson structure on the seed ${\cal X}$-tori. 
Therefore the space ${\cal X}$ has a Poisson structure.
\end{lemma}

{\bf Proof}.   
This can easily be checked directly, and also 
follows from a similar but stronger statement 
about the $q$-deformed cluster ${\cal X}$-varieties 
proved, independently of the Lemma,  in Lemma 
\ref{4.28.03.1}.

\vskip2mm

{\it A Poisson structure on the real tropical ${\cal
    X}$-space}. 
Given a seed, we define a Poisson bracket $\{*, *\}$ 
on the space ${\cal X}(\R^t)$ by
$
\{x_i, x_j\}:= \widehat \varepsilon_{ij}. 
$  
Since mutations are 
given by piecewise linear transformations, it makes sense 
to ask whether it is invariant under mutations -- the invariance of the Poisson structure should be understood on the domain of differentiability.   
It is easy to check that 
this Poisson bracket 
does not depend on the choice of the seed.

\subsection{The ${\cal X}$-space is 
fibered over the torus $H_{{\cal X}}$}
Consider the left kernel of the form (\ref{[]}): 
\begin{equation} \label{11.14.03.1}
{\rm Ker}_L[\ast, \ast] :=  \left\{l \in \Lambda  ~~
  |~~  [l, v] =0 \quad \mbox{for every $v \in\Lambda
  $} \right\}. 
\end{equation}
Given a seed ${\mathbf i}$, there is an isomorphism 
$$
{\rm Ker}_L[\ast, \ast] = \left\{\{\alpha_i\} \in \Z^I ~~
  |~~  \sum_{i\in I}\alpha_i \varepsilon_{ij} =0 \quad \mbox{for every $j \in
  I$} \right\}.
$$
Denote by $H_{{\cal X}}$ the split torus with the 
group of characters ${\rm Ker}_L[\ast, \ast]$.  
The tautological inclusion 
${\rm Ker}_L[\ast, \ast] \hra \Lambda$ provides a surjective homomorphism 
$$
\theta: {\cal X}_{\Lambda} \to  H_{{\cal X}}.
$$
Denote by $\chi_{\alpha}$ the character  of the torus 
$H_{{\cal X}}$ 
corresponding to $\alpha \in {\rm Ker}_L[\ast, \ast]$. 
In the cluster coordinates assigned to a seed ${\bf i}$ we have  
$\theta^*\chi_{\alpha} = \prod_{i \in I} X_i^{\alpha_i}$.

\begin{lemma} \label{11.14.03.3} The following diagram is commutative: 
$$
\begin{array}{ccc}
{\cal X}_{\Lambda}& \stackrel{\mu^x_{k}}{\lra} &{\cal X}_{\Lambda}\\
\theta\downarrow && \downarrow \theta\\
H_{{\cal X}}& \stackrel{\sim}{\lra} & H_{{\cal X}}
\end{array}
$$
\end{lemma}

{\bf Proof}. Follows from the  quantum version, 
see Lemma \ref{1.09.04.100a}.

\vskip2mm

Let us interpret the torus $H_{{\cal X}}$ as a 
tautological positive space, i.e. as a  
functor 
\be \la{thet}
\theta: {\cal
  G} \lra \mbox{ the category of split algebraic tori}, 
\ee
sending objects to the torus $H_{{\cal X}}$, and morphisms to the identity map. 
Then  Lemma \ref{11.14.03.3} implies

\begin{corollary} \label{8.1.04.11}
There is  a unique  map of positive spaces
$
\theta: {\cal X} \lra H_{{\cal X}}
$ 
such that for any seed
$$
\theta^*(\chi_{\alpha}):= \prod_{i \in I} X_i^{\alpha_i}.
$$
\end{corollary}

 Let $e$ be the unit of $H_{{\cal X}}$. 
Thanks to Lemma \ref{11.14.03.3} the fibers $\theta^{-1}(e)$ of the maps $
\theta: {\cal X}_{\Lambda} \lra H_{{\cal X}}$ 
are  glued 
into a positive space. It is the space ${\cal U}$ from
Definition \ref{8.15.04.3}. 

Similarly, for a general $h$,  the fibers $\theta^{-1}(h)$ 
can be glued into an object generalizing
positive space. We are not going to develop   this point of view, observing
only that the  $\R_{>0}$-points of the fibers make sense as   manifolds. 

\begin{proposition} \label{1.8.04.12} 
a) The fibers of the map $\theta$ are the symplectic leaves
of the Poisson space ${\cal X}$. 

b) In particular the fibers of the map 
$\theta: {\cal X}(\R_{>0}) \lra {H}_{{\cal X}}(\R_{>0})$  
are the symplectic leaves of the Poisson manifold 
${\cal X}(\R_{>0})$. 
\end{proposition}

{\bf Remark}. Here in a) by the fibers we mean the collection of varieties 
$\theta^{-1}(h)\subset {\cal X}_{{\mathbf i}}$ and birational isomorphisms between them provided 
by Lemma \ref{11.14.03.3}. The claim is that they 
are symplectic leaves in the tori ${\cal X}_{{\mathbf i}}$, and the gluing maps 
respect the symplectic structure.

{\bf Proof}. Follows immediately from Lemma \ref{1.18.04.1}.

\vskip3mm
\subsection{The torus ${H}_{{\cal A}}$ acts on the ${\cal A}$-space} 
Consider  the right kernel of the form $[\ast, \ast]$: 
\be \la{Atorus}
{\rm Ker}_R[\ast, \ast]:= \{l \in \Lambda ~~ | ~~
 [v, l] = 0 \quad \mbox{for any $v \in \Lambda$}\}. 
\ee
Given a seed ${\bf i}$, there is an isomorphism
$$
{\rm Ker}_R[\ast, \ast] = \left\{\{\beta_j\}\in \Z^I 
\stackrel{\sim}{=} \Lambda ~~ | ~~
 \sum_{j \in I}\varepsilon_{ij}\beta_j =0 \quad \mbox{for any $i \in I$}\right\}.
$$

Recall the map $\delta_{\bf i}$, see (\ref{deltamap}). 
Thanks to Lemma \ref{8.4.09.1},  the lattice 
\be \la{KA}
{\rm K}_{\cal A}:= \delta_{\bf i}({\rm Ker}_R[\ast, \ast]) \subset \Lambda^\vee
\ee 
does not depend on the choice of ${\bf i}$. It is, of course, isomorphic 
to the lattice ${\rm Ker}_R[\ast, \ast]$. 

Let ${H} _{{\cal A}}$ 
be the torus with the group of cocharacters (\ref{KA}). 
Observe that the lattice $\Lambda^\vee$ is the group of cocharacters of the torus ${\cal A}_\Lambda$. 
Thus  
the homomorphism 
$$
{\rm K}_{\cal A} \times \Lambda^{\vee}  \lra \Lambda^{\vee}
$$ 
given by action  of the lattice  ${\rm K}_{\cal A}$ 
on the lattice $\Lambda^{\vee}$
gives rise to a homomorphism of tori
\be \la{po}
{H}_{{\cal A}}\times {\cal A}_{\Lambda} \to {\cal A}_{\Lambda}, \quad 
\chi_{\beta}(t) \times (a_1, ..., a_n) \lms (t^{\beta_1}a_1, ..., t^{\beta_n}a_n), 
\ee
where  $\chi_{\beta}: {\Bbb G}_m \to
{H}_{{\cal A}}$  is the cocharacter
assigned to $\beta$.

\begin{lemma} \label{11.14.03.3a}
a) The maps (\ref{po}) glue into an action of the torus ${H}_{{\cal A}}$ 
on the  ${\cal A}$-space. 

b) The projection $p: {\cal A}\lra {\cal X}$ is the factorization by the action of the  torus 
${H}_{{\cal A}}$. Moreover, there is an ``exact sequence'':
$$
{\cal A}\stackrel{p}{\lra} {\cal X}\stackrel{\theta}{\lra} {H}_{{\cal X}}\lra 1, \qquad 
{\rm Im}\, p = \theta^{-1}(1).
$$ 
%which amounts for every seed ${\bf i}$ to an exact sequence of tori 
%$$
%1 \lra  {H}_{{\cal A}}\lra{\cal A}_{\bf i}\stackrel{p}{\lra} {\cal X}_{\bf i}\lra {H}_{{\cal X}}\lra 1. 
%$$
\end{lemma}

{\bf Proof}. a) Taking into account the decomposition of the mutations, 
the claim amounts to the commutativity of the following diagram:
$$
\begin{array}{ccc}
{H}_{{\cal A}}\times {\cal A}_{\Lambda} & \lra & {\cal A}_{\Lambda} \\
{\rm Id} \downarrow \mu_k^\sharp && \downarrow \mu_k^\sharp\\
{H}_{{\cal A}}\times {\cal A}_{\Lambda} &\lra &{\cal A}_{\Lambda}
\end{array}
 $$
The coordinate $A_k$ under the composition of the right and top arrows
 transforms as follows:
$$
A_i \lms A_i(1+p^*X_k)^{-\delta_{ik}} \lms t^{\beta_i}A_i(1+p^*X_k)^{-\delta_{ik}}. 
$$
The other composition gives the same. Indeed, since $\sum_j \varepsilon_{kj}\beta_j=0$, 
 the transformation 
$A_j \lms t^{\beta_j}A_j$ does not change expression (\ref{qa}). 
So the two way to compute the transformation of the coordinate $A_k$ in the 
diagram lead to the same result. 

b) Clear. The lemma is proved.

\begin{lemma} \label{11.14.03.3c}
There is a canonical isomorphism
$
{H}_{{\cal A}}\otimes  \Q = {H}_{{\cal X}}\otimes  \Q 
$. 
\end{lemma}

{\bf Proof}. This just means that 
$
({\rm Ker}_R[\ast, \ast]\otimes  \Q)^{\vee}
 = {\rm Ker}_L[\ast, \ast]\otimes  \Q
$.

\vskip2mm

{\bf Remark}. Observe that ${\rm Ker}_R[\ast, \ast]= 
{\rm Ker}_L[\ast, \ast]^t$. 
There is a canonical isomorphism
\begin{equation}\label{1.9.04.10}
X_*(H_{\cal A}) = X^*(H_{{\cal X}^{\vee}}). 
\end{equation}
Indeed, both abelian groups are identified with ${\rm Ker}_R[\ast, \ast]$. It
 plays an essential role in Section 4.

\vskip3mm

\begin{lemma} \label{1.8.04.12c} 
There are canonical group homomorphisms
$$
\Gamma \lra {\rm Aut}\left({\rm Ker}_L[\ast, \ast] \right), 
\qquad 
\Gamma \lra {\rm Aut}\left({\rm Ker}_R[\ast, \ast] \right).
$$ 
\end{lemma}

{\bf Proof}. Clear from the very definition.

\subsection{Cluster modular groups revisited} 

{\bf The simplicial complex ${\Bbb S}$ \cite{FZI}.} 
Let ${\mathbf i}$ be a seed, $n:=|I|$, $m:=|I-I_0|$. 
Let $S$ be an $(n-1)$-dimensional
simplex, equipped with a bijection (decoration) 
\begin{equation} \label{17:00a}
\{\mbox{codimension one faces of $S$}\} \stackrel{\sim}{\lra} I.
\end{equation} 
We call it a ${I}$-decorated, or simply decorated simplex. 

Take a decorated simplex, and glue to it $m$ other decorated simplices as follows. 
To each codimension one face of the initial simplex
 decorated by an element $k\in I-I_0$ we glue a new decorated 
simplex along its  codimension one face decorated by the same $k$. Then to each of the remaining codimension one faces decorated by the 
elements of $I-I_0$ we glue new  decorated 
simplices, and so on, repeating this construction infinitely many times. 
We get a simplicial complex ${\Bbb S}$. 
Let ${\mathbf S}$ be the set of all its simplices.

We connect two elements of ${\mathbf S}$ by an edge if the 
corresponding simplices share a common codimension one face. 
We get an $m$-valent tree with 
the set of vertices ${\mathbf S}$. 
Its edges are decorated  by 
the elements of the set $I-I_0$. We denote it by ${\rm Tr}$.

\begin{lemma}
There are canonical bijections:
$$
\{\mbox{Seeds equivalent to a seed ${\mathbf i}$}\} \leftrightarrow  
\{\mbox{The set ${\mathbf S}$ of 
simplices of the simplicial complex ${\Bbb S}$}\},
$$
$$
\{\mbox{compositions of seed  mutations}\} \leftrightarrow 
\{\mbox{Paths on the tree ${\rm Tr}$}\}.
$$
\end{lemma}

{\bf Proof}. The seed ${\mathbf i}$ is  assigned to the original simplex $S$. 
Given any other simplex $S'$ of the simplicial complex ${\Bbb S}$, there is a unique path 
on the tree ${\rm Tr}$ 
connecting $S$ with $S'$. It gives rise to a sequence of mutations parametrised by the edges of the path, so that 
the edge decorated by $k$ gives rise to the mutation in the direction $k$. 
Mutating the seed ${\mathbf i}$ by this sequence of mutations, we get the seed assigned to the simplex $S'$. The lemma follows.

\vskip 3mm
{\bf Remark}. 
Cluster ${\cal A}$-coordinates are
assigned to vertices of the simplicial complex ${\Bbb S}$. 
Cluster ${\cal X}$-coordinates are assigned
to  cooriented faces of ${\Bbb S}$.  
 Changing coorientation  
amounts to inversion of the corresponding cluster ${\cal X}$-coordinate.
Mutations are parametrized by codimension one 
faces of ${\Bbb S}$.

\paragraph{Another look at the cluster modular groups.} 
Let $ {\bf F}({\Bbb S})$ be the set of all pairs $(S,F)$ where $S$ 
is a simplex of ${\Bbb S}$, and $F$ is a codimension one face of $S$. 
Pairs of faces belonging to the same simplex are parametrized 
by the fibered product 
\begin{equation} \label{5.20.03.11}
{\bf F}({\Bbb S})\times_{{\bf S}} {\bf F}({\Bbb S})
\end{equation}
The collection of exchange functions $\varepsilon_{ij}$ can be viewed as a single 
function ${\cal E}$ on the set (\ref{5.20.03.11}).

\vskip 2mm
Let  ${\rm Aut}({\Bbb S})$ be the automorphism group  of the 
simplicial complex ${\Bbb S}$. It contains the 
subgroup ${\rm Aut}_0({\Bbb S})$ respecting the decorations (\ref{17:00a}). 
The group ${\rm Aut}({\Bbb S})$ is a semidirect product:
$$
0 \lra {\rm Aut}_0({\Bbb S}) \lra  {\rm Aut}({\Bbb S})
\stackrel{}{\lra}  {\rm Per}
\lra 0.
$$ 
Here ${\rm Per}$ is the group of automorphisms of the pair $(I, I_0)$. 
Given ${\mathbf i} \in {\bf S}$, the 
group ${\rm Per}$ is realized as a subgroup ${\rm Aut}({\Bbb S})$ 
 permuting the faces of $S_{{\mathbf i}}$.

The group ${\rm Aut}({\Bbb S})$ 
acts on the set 
(\ref{5.20.03.11}), and hence on the set of  exchange functions ${\cal E}$. 

\begin{definition} \label{5.19.03.1}  
Let ${\cal E}$ be 
the exchange function corresponding to a seed  ${\bf i}$.  

The group $D$ is the subgroup of  ${\rm Aut}({\Bbb S})$
 preserving ${\cal E}$: 
\begin{equation} \label{5.19.03.10}
D:= \{\gamma \in {\rm Aut}({\Bbb S}) \quad | 
\quad \gamma^* ({\cal E}) = 
{\cal E}\}; 
\end{equation}

The group $\Delta$  is the subgroup of $D$
  preserving the cluster ${\cal A}$- and ${\cal X}$-coordinates:
 \begin{equation} \label{5.19.03.11}
\Delta:= \{\gamma \in 
D \quad | \quad \gamma^*A_{j} = A_j, \quad \gamma^*X_{j} = X_j\};
\end{equation}

The {cluster complex} $C$ is the quotient of 
${\Bbb S}$ by the action of the group $\Delta$:
$$
C: 
= {\Bbb S}/\Delta. 
$$
\end{definition} 
Clearly $\Delta$ is a normal subgroup of $D$. 
\begin{lemma} \label{7.2.03.2}
The quotient group
$
\Gamma:= D/\Delta  
$ 
is the   cluster modular 
 group.  
\end{lemma}

{\bf Proof}.  
For any two simplices $S_{\mathbf i}$ and $S_{\mathbf i'}$ of ${\Bbb S}_{I, I_0}$ there exists a unique  
element of the group $ {\rm Aut}_0({\Bbb S}_{I, I_0})$ transforming 
$S_{\mathbf i}$ to $S_{\mathbf i'}$. So given an element $d\in D$, 
there is a cluster transformation ${\bf c}_d:{\bf i} \to {\bf i'}$. 
Then by (\ref{5.19.03.10}) there is a seed isomorphism  
$\sigma_d: {\bf i'} \to {\bf i}$. 
Thanks to (\ref{5.19.03.11}) the cluster transformation
$\sigma_d\circ {\bf c}_d:{\bf i} \to {\bf i}$ is 
trivial.

\vskip 3mm

{\bf Variants}. In Definition \ref{5.19.03.1} and Lemma \ref{7.2.03.2} 
we looked how  the ${\cal A}$- and ${\cal X}$-coordinates 
behave under cluster transformations. There are similar definitions 
using either 
${\cal A}$-coordinates, or  ${\cal X}$-coordinates. 
This way we get the groups $\Delta_{?}$, $\Gamma_{?}$, 
and the simplicial complexes $C_{?}$, where $?$ stands, respectively, 
 for 
${\cal A}$ and  
${\cal X}$. The cluster  complex ${C}_{\cal A}$ was defined in [FZII]. 
Corollary 
\ref{8.24.05.1x} immediately implies 

\begin{lemma} \label{9.5.05.2} Assume that ${\rm det}\varepsilon_{ij} \not = 0$. Then $\Delta_{{\cal A}} = \Delta \subset \Delta_{{\cal X}}$, 
$C_{{\cal A}} = C$, and 
$\Gamma_{{\cal A}} = \Gamma$. 
\end{lemma}

\subsection{Example: Cluster transformations for  
${\cal X}$-varieties of types $A_1\times A_1, A_2, B_2, G_2$.} \la{sec2.5}

The results of Section 2.5 play a crucial role in Section 5.2. 
We start by an 
elaboration of the cluster transformations for cluster 
${\cal X}$-variety of type $ B_2$. Its main goal is to tell the reader 
how we picture mutations, quivers etc. 
The obtained formulas, however,  
do not seem very illuminating. 

We show  that the situation changes dramatically when we go to the tropicalised 
cluster transformations. 
Notice that they  
contain just the same  information as the original cluster transformations. 
The advantage of the tropicalised formulas 
stems from the fact that they are piecewise 
linear transformations, and thus can be perceived geometrically. 
We demonstrate this idea by working out 
tropicalisations of cluster transformations 
for  every finite type cluster 
${\cal X}$-variety of rank two, uncovering an interesting geometry 
standing behind.

\paragraph{Quivers and cluster transformations in the $B_2$ case.} 
We picture a seed 
by a quiver  with two vertices. The cluster ${\cal X}$-coordinates 
assigned to a seed are the functions 
written near the corresponding vertices. Every two neighboring seeds 
are related by a horizontal arrow, associated with one of the vertices 
of the left quiver. It shows a mutation in the direction of that vertex. 
The exchange function $\varepsilon$ is determined as follows. 
Denote by $b$ and $t$ the bottom and top vertices. 
Then for the very left quiver $\varepsilon_{bt} = -2$, $\varepsilon_{tb} = 1$. 
For the next one, $\varepsilon_{bt} = 2$, $\varepsilon_{tb} = -1$, and so on. 
This sequence of mutations has period $6$. 
Similar calculations 
can be done for the seeds of types $A_1\times A_1, A_2, G_2$. 

$$
\begin{array}{ccccccccc}
~ &
x & \to &
x^{-1} & ~ &
x^{-1}(1+y+xy)^2 & \to &
x(1+y+xy)^{-2} ~ \\
~ & \includegraphics{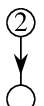} &
~ & \includegraphics{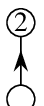} &
~ & \includegraphics{b2.eps} &
~ & \includegraphics{b2up.eps} \\
\to &
y & ~ &
y(1+x) & \to &
y^{-1}(1+x)^{-1} & ~ &
x^{-1}y^{-1}(1+2y+y^2+xy^2) & \to \\
\end{array}
$$
$$
\begin{array}{ccccccccc}
~ &
x^{-1}(1+y^{-1})^2 & \to &
x(1+y^{-1})^{-2} & ~ &
x& \to \\
~ & \includegraphics{b2.eps} &
~ & \includegraphics{b2up.eps} &
~ & \includegraphics{b2.eps}
~ \\
\to &
xy(1+2y+y^2+xy^2)^{-1} & ~ &
y^{-1} & \to &
y& ~ \\
\end{array}
$$

\paragraph{Geometry of the tropicalised cluster transformations in the finite type 
rank two case.} Take  a finite type rank two seed ${\bf i}$, with  $I= \{1,2\}$. 

\begin{figure}
\begin{center}
\raisebox{50pt}{$A_2$:}\includegraphics{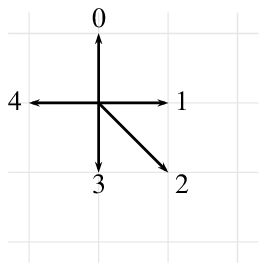}
\raisebox{50pt}{$B_2$:}\includegraphics{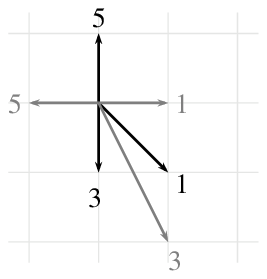} 
\includegraphics{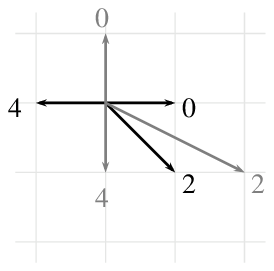}
\end{center}

\begin{center}
\raisebox{70pt}{$G_2$:}\includegraphics{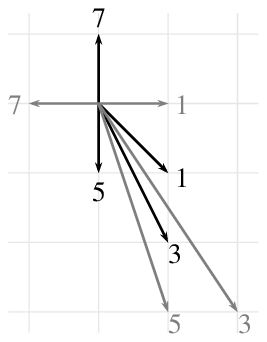} 
\includegraphics{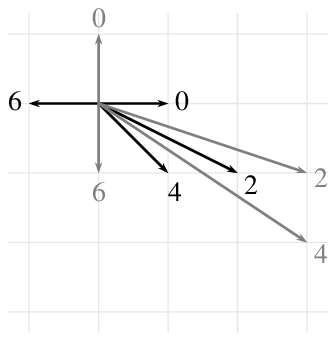}
\end{center}
\caption{Orbits of mutations.}
\label{Orbits}
\end{figure}

\textit{Case $A_2$.}  In this case all seeds are isomorphic. 
Consider a cluster transformation 
$$
\mu:= \sigma_{1,2} \circ \mu_1: {\bf i} \lra \sigma_{1,2}({\bf i}) \sim {\bf i}.
$$
Let $P$ be the tropicalisation of the cluster 
${\cal X}$-torus corresponding to the seed ${\bf i}$. 
It is a plane  with  coordinates $(x,y)$. 
The cluster transformation $\mu$ induces a map  
$$
\mu^t: P \lra P, \quad x \lms y+{\rm max}(0,x), \quad y\lms -x;
$$
The plane $P$ is decomposed into $5$ sectors as shown on Fig \ref{Orbits}. 
Three of them are coordinate quadrants. The other two are obtained by 
subdividing the remaining  quadrant into two sectors. 
We order the 
sectors clockwise cyclically, starting from the 
positive quadrant $\{(x,y)| x,y\geq 0\}$. 

\bl \la{OHK1} The map $\mu^t$ moves the $i$-th sector to the $(i+1)$-st. 
Its restriction to any sector is  linear. 
The sectors are the largest domains in $P$ on which 
any power of the map $\mu^t$ is linear.  
\el

{\bf Proof}. The vectors $(0,1), (1,0), (1,-1), (0,-1), (-1,0)$ 
are the primitive integral vectors generating the boundary arrows 
of the domains. The map $\mu^t$ moves 
them cyclically clockwise. The Lemma follows easily from this. 

\bc
One has $\mu^5 = {\rm Id}$. The element $\mu$ 
generates the modular group of the cluster ${\cal X}$-variety of type $A_2$, 
and  identifies it with $\Z/5\Z$. 
\ec

\textit{Cases $A_1\times A_1$.} This is the simplest case. We have 
$$
\mu^t: P \lra P, \quad x \lms y, \quad y\lms -x;
$$
There are four sectors in this case, given by the coordinate quadrangles. 

\bl \la{OHK0} The map $\mu^t$ moves the $i$-th sector to the $(i+1)$-st. 
Its restriction to any sector is  linear. 
The sectors are the largest domains in $P$ on which 
any power of the map $\mu^t$ is linear.  

One has $\mu^4 = {\rm Id}$. The element $\mu$ 
generates the modular group of the cluster ${\cal X}$-variety of type $A_1\times A_1$, 
and  identifies it with $\Z/4\Z$.
\el

\textit{Cases $B_2$ and  $G_2$.} In these cases there are 
two non-isomorphic seeds, denoted ${\bf i}_-$ and ${\bf i}_+$. 
Consider  cluster transformations 
$$
\mu_-:= \sigma_{1,2} \circ \mu_1: {\bf i}_- \lra \sigma_{1,2}({\bf i}_+), \qquad 
\mu_+:= \sigma_{1,2} \circ \mu_1: {\bf i}_+ \lra \sigma_{1,2}({\bf i}_-).
$$
There are two tropical planes  $P_-$ and $P_+$ with coordinates $(x,y)$,
which are tropicalisations of the cluster 
${\cal X}$-tori corresponding to the seeds ${\bf i}_-$ and ${\bf i}_+$. 
The cluster transformations $\mu_{\pm}$ induce the maps 
$$
\mu^t_-: P_- \lra P_+, \quad x \lms y+{\rm max}(0,x), \quad y\lms -x;
$$
$$
\mu^t_+: P_+ \lra P_-, \quad 
x \lms y+c ~{\rm max}(0,x), \quad y\lms -x.
$$
Here $c=2,3$ for the Dynkin diagrams $B_2, G_2$ respectively. 
The maps $\mu_{\pm}^t$ have the following geometric description. 
Each of the planes $P_{\pm}$ is decomposed into a union of $h+2$ sectors. 
These sectors  include all 
coordinate quadrants but the bottom right one. The remaining sectors 
subdivide that quadrant as shown on Fig \ref{Orbits}. 
The boundaries of these sectors 
are arrows whose directing vectors are:
$$
B_2: \qquad P_-: (1,-1), (1,-2); \qquad P_+: (2, -1), (1,-1). 
$$
$$
G_2: \qquad P_-: (1,-1), (2, -3) (1,-2), (1, -3);  \qquad P_+: 
(3, -1), (2, -1), (3, -2), (1,-1). 
$$
The remaining four  directing vectors are $(0,1), (1,0), (0,-1), (-1,0)$ in both cases. 
Let us order the 
sectors clockwise cyclically, starting from the 
positive coordinate quadrant.

\bl \la{OHK2} The map $\mu_-^t$ moves the $i$-th sector on $P_-$ to the $(i+1)$-st sector on 
$P_+$. The map $\mu_+^t$ moves the $i$-th sector on $P_+$ to the $(i+1)$-st sector on 
$P_-$. They are linear maps on the sectors. 
The sectors are the largest domains in $P_\pm$ on which 
any composition of the map $\mu^t_\pm$ is linear.  
\el

{\bf Proof}. In the  $B_2$ case we get two sequences of vectors, 
shown by black and grey on Fig \ref{Orbits}:
\be \la{per}
 (0, 1) \stackrel{\mu_-}{\lra} (1,0)  \stackrel{\mu_+}{\lra} (1, -1)
\stackrel{\mu_-}{\lra}  (1, -1)  \stackrel{\mu_+}{\lra}  (0, -1)
\stackrel{\mu_-}{\lra}   (-1,0) \stackrel{\mu_+}{\lra}(0, 1). 
\ee
\be \la{vt}
 (1, 0) \stackrel{\mu_-}{\lra} (2, -1)  \stackrel{\mu_+}{\lra} (1, -2)
\stackrel{\mu_-}{\lra}  (0, -1)  \stackrel{\mu_+}{\lra}  (-1, 0)
\stackrel{\mu_-}{\lra}   (0, 1) \stackrel{\mu_+}{\lra}(1, 0). 
\ee
In the $G_2$ case we also get two  sequences of vectors, 
shown by black and grey on Fig \ref{Orbits}:
$$
 (0, 1) \stackrel{\mu_-}{\lra} (1,0)  \stackrel{\mu_+}{\lra} (1,-1)
\stackrel{\mu_-}{\lra}  (2, -1)  \stackrel{\mu_+}{\lra}  (1,-2)
\stackrel{\mu_-}{\lra}   (1,-1) \stackrel{\mu_+}{\lra}(0, -1)
\stackrel{\mu_-}{\lra}   (-1, 0) \stackrel{\mu_+}{\lra}(0, 1). 
$$
$$
 (1,0) \stackrel{\mu_-}{\lra} (3,-1)  \stackrel{\mu_+}{\lra} (2,-3)
\stackrel{\mu_-}{\lra}  (3, -2)  \stackrel{\mu_+}{\lra}  (1,-3)
\stackrel{\mu_-}{\lra}   (0,-1) \stackrel{\mu_+}{\lra}(-1, 0)
\stackrel{\mu_-}{\lra}   (0, 1) \stackrel{\mu_+}{\lra}(1, 0). 
$$
One checks that the cluster transformations 
$\mu_{\pm}$ are linear on each of the sectors

\bc One has $(\mu_+\mu_-)^3 = {\rm Id}$ in the case $B_2$,  
and $(\mu_+\mu_-)^4 = {\rm Id}$ in the case $G_2$. 

The modular group of the cluster 
${\cal X}$-variety of type $B_2$ is $\Z/3\Z$. Its generator is $\mu_+\mu_-$.

The modular group of the cluster 
${\cal X}$-variety of type $G_2$ is $\Z/4\Z$. Its generator is $\mu_+\mu_-$.
\ec

\paragraph{Proof of Proposition \ref{FTFTFT}.}
It is easy to check that the ${\cal A}$-coordinates 
in the rank two  cases have the same period.  
This settles the Proposition for the rank two case. 
The claim in general for the ${\cal A}$-coordinates as well as 
the exchange functions $\varepsilon_{ij}$ was proved in \cite{FZI}. 
The claim for the ${\cal X}$-coordinates can be reduced to it  via the following trick. 
One can find a set $I'$ containing $I$ and a skew-symmetrisable 
exchange function $\varepsilon'_{ij}$ on $I' \times I'$ extending 
$\varepsilon_{ij}$ on $I\times I$ 
such that ${\rm det}\varepsilon'_{ij}\not = 0$. 
Then the claim follows from Corollary \ref{8.24.05.1x}, since 
the composition of the ${\cal A}$-mutations assigned 
to the standard $(h+2)$-gon is trivial for any seed.

\subsection{Modular complexes, 
modular orbifolds and 
special modular groups}  A simplex of the simplicial complex $C_{?}$ 
is of {\it finite type} if the set of all simplices of $C_{?}$ 
containing it is finite.  Here $?$ stands for 
${\cal A}$, ${\cal X}$, or  
no label at all. 

\begin{definition} \label{12.9.03.1} The reduced cluster complex 
$C^*_{?}$ is the union of finite type simplices of $C_{?}$. 
\end{definition}

The reduced cluster complex 
is not a simplicial complex:
 certain 
faces of its
simplices may not belong to it. But it  has a  
 topological realization. 

\begin{theorem} \label{12.9.03.2} 
Topological realizations 
of the reduced cluster complexes $C^*$ and  $C^*_{\cal A}$ 
are homeomorphic to manifolds.
\end{theorem} 

{\bf Proof}. We give a proof for the cluster complex $C^*$ --  
the case of $C^*_{\cal A}$ similar, and a bit simpler. 

A simplicial complex 
is of {\it finite type} if it has has a finite number of simplices. 
According to the 
Classification Theorem \cite{FZII} cluster algebras of finite type, i.e. the ones with 
the cluster complexes of finite type, 
are classified  by the Dynkin diagrams of type $A, B, ..., G_2$. 
The cluster complex of type $A_n$ is a  Stasheff polytope. The cluster complexes 
corresponding to other finite type cluster algebras are the 
{\it generalized associahedra}, or generalized Stasheff polytopes  
\cite{FZ}.

We need the following crucial lemma.

\begin{lemma} \label{12.9.03.3}
Let $S'_{\mathbf i}$ be a simplex of finite type in the cluster complex 
$C$. Then 
the set of all simplices containing $S'_{\mathbf i}$ is naturally identified with the set of 
all faces of a generalized Stasheff polytope, so that the codimension 
$i$ simplices correspond to the $i$-dimensional faces. 
\end{lemma}

{\bf Proof}. Let $I_{\mathbf i}$ be the 
set of vertices of a top dimensional simplex $S_{\mathbf i}$ of 
the simplicial complex $C$. The exchange function $\varepsilon$ is a 
function on 
$I_{\mathbf i} \times I_{\mathbf i}$.  
Let $S'_{\mathbf i}$ be a finite type simplex 
contained in the simplex $S_{\mathbf i}$. 
Let $I'_{\mathbf i} \subset I_{\mathbf i}$ be 
the subset of the   vertices of $S'_{\mathbf i}$. 
The set of top dimensional simplices of $C$ 
containing $S'_{\mathbf i}$ is obtained from the simplex $S_{\mathbf i}$ 
as follows. Consider the seed 
defined by the function $\varepsilon'$ 
with the frozen variables parametrized by 
the subset $I'_{{\mathbf i}}$. Recall that this means 
 that we do mutations only at the vertexs 
of $I_{\mathbf i} - I'_{\mathbf i}$. Since the 
simplex $S'_{\mathbf i}$ is finite type, it gives rise to 
a finite type cluster algebra. Indeed, since by the very definition 
$\Delta$ is a subgroup of $\Delta_{\cal A}$, the simplicial complex 
$C$ has more simplices than  $C_{\cal A}$. So if $C$ is of finite type, 
$C_{\cal A}$ is also of finite type. Therefore 
the matrix $\varepsilon_{ij}$ is non-degenerate by the 
Classification Theorem. 
So by Corollary \ref{8.24.05.1x} in our case $C_{\cal A} = C$. 
The corresponding cluster complex 
is the generalized Stasheff polytope corresponding 
 to the Cartan matrix assigned to the 
 exchange function $\varepsilon'$ on the set 
$(I_{\mathbf i} - I'_{\mathbf i})^2$. It is a convex polytope  \cite{CFZ}. 
%See Theorem \ref{PiR} for a simple proof. 
This proves the lemma. 

\vskip2mm

Let us deduce the theorem from this lemma. Consider a convex polyhedron $P$. 
Take the dual decomposition of its boundary, and connect each of 
the obtained polyhedrons with a point inside of $P$ by straight lines. 
We get a conical decomposition of $P$. 
Let us apply this construction to the generalized Stasheff polytope. Then 
the product of the interior part of the simplex $S'_{\mathbf i}$ and the defined above 
conical decomposition of the generalized Stasheff polytope 
corresponding to the 
exchange function $\varepsilon'$ on 
$(I_{\mathbf i}-I'_{\mathbf i})^2$ gives the link of the interior part of the simplex $S'_{\mathbf i}$. 
In particular 
a neighborhood of any interior point of $S'_{\mathbf i}$ is topologically a ball. 
The theorem is proved.

\begin{conjecture} \label{7.7.04.1}
The simplicial complex $C_{{\cal X}}$ is of finite type 
if and only if $C$ is of finite type. 
\end{conjecture} 

Proposition \ref{FTFTFT} implies that 
Conjecture \ref{7.7.04.1} is valid if $|I - I_0|=2$.
\vskip 2mm

%\vskip 2mm
%\noindent
%{\bf Type} $A_1 \times A_1$. Then $b=c=0, h=2$ and 
%$$
%x_3 = x_1^{-1}, \quad x_4 = x_2^{-1}.
%$$
%{\bf Type} $A_2$. Then $b=c=1, h=3$ and 
%$$
%x_3 = \frac{1+x_2}{x_1}, \quad x_4 = \frac{1+x_1+x_2}{x_1x_2}, \quad 
%x_5 = \frac{1+x_1}{x_2}.
%$$
%{\bf Type} $B_2$.  Then $b=1, c=2, h=4$ and 
%$$
%x_3 = \frac{1+x_2}{x_1}, \quad x_4 = \frac{(1+x_1+x_2)^2}{x^2_1x_2}, \quad 
%x_5 = \frac{(1+x_1)^2+x_2}{x_1x_2},\quad x_6 = \frac{(1+x_1)^2}{x_2}.
%$$
%{\bf Type} $G_2$.  Then $b=1, c=3, h=6$ and 
%$$
%x_3 = \frac{1+x_2}{x_1}, \quad x_4 = \frac{(1+x_1+x_2)^3}{x^3_1x_2}, \quad 
%x_5 = \frac{(1+x_1)^3+ x_2(3x_1+x_2+2)}{x^2_1x_2},
%$$
%$$
%x_6 = 
%\frac{((1+x_1)^2+x_2)^3}{x_1^3x^2_2},\quad x_7 = 
%\frac{(1+x_1)^3 +x_2}{x_1x_2},\quad x_8 = 
%\frac{(1+x_1)^3 }{x_2}.
%$$

{\it The cluster modular complex}. 
Suppose that we have a decomposition of a manifold on simplices, although 
some faces of certain simplices may not belong to 
the manifold. Then  
the dual polyhedral decomposition of the  manifold is a polyhedral complex.  
Its topological realization is homotopy equivalent to the manifold. 

Thanks to Theorem \ref{12.9.03.2} the reduced cluster complex $C^*$ is homeomorphic to a manifold. 
Therefore the dual polyhedral complex for $C^*$ is a polyhedral complex whose topological realization 
is homotopy equivalent 
to a manifold. This motivates the following two definitions. 

\begin{definition} \label{12.9.03.4}
The {\em cluster modular complex $\widehat M$} is the dual polyhedral complex 
for the reduced cluster complex 
$C^*$. 
\end{definition}

{\it The modular orbifold}. 
The   cluster modular group $\Gamma$ acts on $C^*$, and hence on 
 $\widehat M$. The stabilizers of points are finite groups.

\begin{definition} \label{12.8.03.1} The {\em cluster modular orbifold} $M$  is the orbifold $ 
\widehat M/\Gamma$. 
\end{definition}

The  fundamental groupoid of a polyhedral complex ${P}$ is a
groupoid whose objects are vertices of $P$, and 
morphisms are homotopy classes 
of paths between the vertices. 

\begin{theorem} \label{7.2.03.4}
The special modular groupoid $\widehat {\cal G}$ 
is the  fundamental groupoid of 
the cluster modular orbifold $M$. 
The special modular  group $\widehat {\Gamma}$ 
is the fundamental group of the orbifold $M$, centered at a vertex of $M$
\end{theorem} 

{\bf Proof}. The morphisms in the 
fundamental groupoid of a polyhedral complex ${P}$ can be described by 
generators and
relations as follows. The generators are given by the edges of 
$P$. The relations correspond to
the two dimensional cells of $P$. 
So to prove the theorem, we describe 
the $2$-skeleton of the  polyhedral complex
$\widehat M$. 
%, which are  the standard $(h+2)$-gins by Corollary \ref{11.2.03.1}
We start from a reformulation of Lemma \ref{12.9.03.3}:

\begin{corollary}  \label{11.2.03.1d}
Any cell of the polyhedral complex 
$\widehat M$ is  
isomorphic to the generalized Stasheff polytope corresponding to a Dynkin diagram 
from the Cartan-Killing classification.
\end{corollary}

This implies the following description of 
the $2$-skeleton of the  polyhedral complex
$\widehat 
M$. 
The $1$-skeleton of $\widehat M$ is 
the quotient ${\rm Tr}/\Delta$. Let us describe 
the $2$-cells of $\widehat M$. 
As was discussed in Section \ref{sec2.5}, 
 if $\varepsilon_{ij}$ is one of the following matrices 
\begin{equation} \label{12.31.03.1}
\left (\matrix{0& 0\cr 0&0\cr}\right ), \quad 
\pm \left (\matrix{0& 1\cr -1&0\cr}\right ), \quad 
\pm \left (\matrix{0& 1\cr -2&0\cr}\right ), \quad \pm \left (\matrix{0& 1
\cr -3&0\cr}\right ) 
\end{equation} 
then performing 
mutations at the vertexs $i, j, i, j, i,
...$ we get an $(h+2)$-gon, where $h$ is the Coxeter number 
of the Dynkin diagram of type $A_1 \times A_1, A_2, B_2, G_2$ respectively,
i.e. $h=2, 3, 4, 6$.   
These  $(h+2)$-gons  are called the {\it standard 
$(h+2)$-gons} in  $\widehat M$. 
We get

\begin{corollary}  \label{11.2.03.1}
Any $2$-cell of the polyhedral complex $\widehat M$ is a standard
$(h+2)$-gon, $h=2, 3, 4, 6$. 
\end{corollary}

Therefore the edges and $2$-cells of the polyhedral complex $\widehat M$ 
match the generators and the relations of the special modular groupoid from Definition \ref{smg}.
The theorem is proved. 

\vskip 2mm
{\bf Remark}. The vertices 
of the modular orbifold are parametrized by the functions 
on the set $I \times I$ up to permutation equivalence 
obtained by mutations of an  initial exchange function 
$\varepsilon_{ij}$. 
So  
 the number of cells is infinite if and only if the 
absolute value of the cluster 
function  ${\cal E}$ is unbounded. 
\vskip 2mm

\begin{hypothesis} \label{7.2.03.2aq}
The  reduced 
cluster complex $C^*$ (or, equivalently, the modular complex 
$\widehat M$)   is simply connected. 
\end{hypothesis} 
Hypothesis \ref{7.2.03.2aq} is equivalent to the one that 
the  canonical epimorphism $\widehat \Gamma 
\to \Gamma$ is an isomorphism. So in the cases when it is satisfied   
we have a transparent description of the modular group 
$\Gamma$: all relations come from the standard 
$(h+2)$-gons. 

A cluster ensemble is of {\it finite type} if the cluster complex $C$ 
is of finite type.  Lemma  \ref{9.5.05.2} implies that 
classification of finite 
type cluster ensembles is the same as the one 
for cluster algebras. 
Hypothesis \ref{7.2.03.2aq} is valid for cluster ensembles of finite type 
if and only if $|I|>2$.  
Indeed, in the finite type case the cluster complex $C$ 
is the boundary of a generalized Stasheff polytope, which is a convex 
polytope \cite{CFZ}. 
So the topological realization of $C$  is homeomorphic to a 
sphere.

\vskip 3mm
{\bf Remark}. Just the same way as  we prove 
Theorem \ref{12.9.03.2},  
Conjecture  \ref{7.7.04.1} implies that 
the topological realization 
of the reduced cluster complex $C^*_{{\cal X}}$ 
is homeomorphic to a manifold.

%Therefore we can make  the following definitions.

%\begin{definition} \label{12.8.03.1x} Assuming Conjecture  \ref{7.7.04.1}: 

%i) The polyhedral complex $\widehat M_{\cal X, |{\bf i}|}$ is the dual 
%for the reduced cluster complex $C^*_{\cal X, |{\bf i}|}$. 

%ii) $M_{{\cal X}, {|{\bf i}|}}:= 
%\widehat M_{{\cal X}, {|{\bf i}|}}/\Gamma_{{\cal X}, {{\bf i}}}$ is 
%the {\it ${\cal X}$-modular orbifold} of the cluster ensemble.
%\end{definition}

%Thanks to Lemma \ref{7.7.04.1r} the topological realization 
%of $C^*_{{\cal X}, {|{\bf i}|}}$ is homeomorphic to a manifold in codimension $\leq 2$. 
% Therefore one can define the $2$-skeleton ${\rm Sk}_{\leq 2}\widehat M_{\cal X, |{\bf i}|}$ 
%and prove that its cells are the standard $(h+2)$-gons, where $h = 2,3,4,6$. 
%Thus, regardless of Conjecture \ref{7.7.04.1},  we have the following definitions: 

%\begin{definition} \label{12.8.03.1xx} 
%i) ${\cal G}_{{\cal X}, {|{\bf i}|}}$  
%is the  fundamental groupoid of $({\rm Sk}_{\leq 2}
%\widehat C^*_{\cal X, |{\bf i}|})/\Gamma'_{{\cal X}, {{\bf i}}}$. 

%ii) The cluster ${\cal X}$-modular group ${\Gamma}_{{\cal X}, {{\bf i}}}$ 
%is the fundamental group of ${\cal G}_{{\cal X}, {|{\bf i}|}}$.
%\end{definition}

\subsection{Cluster nature of the classical Teichm\"uller space}

\begin{figure}[ht]
\begin{picture}(144,96)(0,-24)
\centerline{\epsfbox{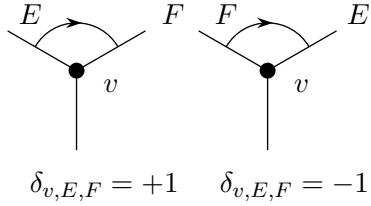}}
\put(-304,54){$E$}
\put(-180,54){$E$}
\put(-230,54){$F$}
\put(-250,54){$F$}
\put(-272,28){$v$}
\put(-200,28){$v$}
\put(-300, -10){$\delta_{v,E,F}=+1$}
\put(-228, -10){$\delta_{v,E,F}=-1$}
\end{picture}
\caption{The function $\delta_{E,F}$.}
\label{cluste00}
\end{figure}

\paragraph{Cluster data for 
the Teichm\"uller space of an oriented hyperbolic surface $S$ with punctures.} 
It was defined in Chapter 11 of \cite{FG1}. 
Consider a trivalent tree $T$ embedded into $S$, homotopy equivalent to 
$S$.  
Let $\Lambda$ be the lattice generated by the edges of $T$. It has a basis given by the edges. 
The skew-symmetric matrix $\varepsilon_{EF}$, where $E$ and $F$ run
through the set of edges of the tree $T$, is 
defined as follows.
Each edge $E$ of   $T$ determines two 
flags, defined as pairs $(v, E)$ where $v$ is 
a vertex of an edge $E$. Given two flags $(v, E)$ and  
$(v, F)$ sharing the same vertex, we define 
$\delta_{v, E, F} \in \{-1, 1\}$ as follows:
$\delta_{v, E, F} = +1$ (respectively $\delta_{v, E, F} = -1$) 
if the edge $F$ goes right after (respectively right before) 
the edge $E$ according to the orientation of the surface, see Fig. \ref{cluste00}. 
For each pair $(E, F)$ of the edges of $T$, consider the set 
$v(E,F)$ of their common vertices. It has at most two elements. 
We set

%Figure 4
\begin{equation} \label{7.5.04.11}
\varepsilon_{EF}:= \sum_{v \in v(E,F)}\delta_{v, E, F} \in \{\pm 2, \pm 1, 0\}. 
\end{equation} 

We defined  
in Chapter 3 of \cite{FG1} a polyhedral complex ${\Bbb G}_{S}$, 
called the {\it modular complex}. Its 
 vertices are 
parametrized by the isotopy classes of trivalent graphs on $S$ which are  homotopy equivalent to $S$. 
Its dimension $k$ faces correspond to the isotopy classes of graphs $G$ on $S$, homotopy equivalent to $S$, 
such that valency ${\rm val }(v)$ of each vertex $v$ of $G$ is $\geq 3$, and $k = \sum_{v}({\rm val }(v)-3)$. 

The modular complex ${\Bbb G}_{S}$ can be identified with
 the cluster modular complex 
for the exchange function (\ref{7.5.04.11}). Since ${\Bbb G}_{S}$ is known to be contractible, 
Hypothesis \ref{7.2.03.2aq} 
is valid in this case. This suggests that Hypothesis \ref{7.2.03.2aq} 
may be valid in a large class of examples. 
The  faces of the modular complex ${\Bbb G}_{S}$ are  the Stasheff polytopes or their products: 
this illustrates Theorem \ref{12.9.03.2}. 
Since ${\Bbb G}_{S}$ is contractible, the 
two versions $T$ and $\widehat T$ 
of the cluster modular group are isomorphic, and 
identified 
with the modular group of $S$.

\begin{figure}[ht]
\centerline{\epsfbox{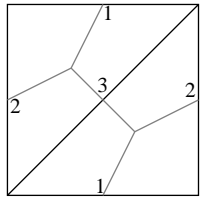}}
\caption{The trivalent ribbon graph corresponding to a punctured torus.}
\label{clus2}
\end{figure}

\paragraph{Cluster data for the Teichm\"uller space of the punctured torus $S$.} 
There is a unique up to isomorphism trivalent
ribbon graph 
corresponding to a punctured torus. It is shown on Fig \ref{clus2} embedded in the
punctured torus: the puncture is at the identified vertices of the square. 
Let us number its edges by $\{1, 2, 3\}$. 
The general recipe in the case of the punctured torus leads to 
a exchange function $\varepsilon_{ij}$ given by 
the skew-symmetric matrix
\begin{equation} \label{1.3.04.1}
\varepsilon_{ij} = \left (\matrix{0& 2& -2\cr -2&0&2\cr2&-2&0
\cr}\right),
\end{equation}
Its quiver is shown on the left of Fig. \ref{quiver}. 
Mutations change the sign of the function $\varepsilon_{ij}$:
$\varepsilon'_{ij} = -\varepsilon_{ij}$.

\begin{figure}[ht]
\centerline{\epsfbox{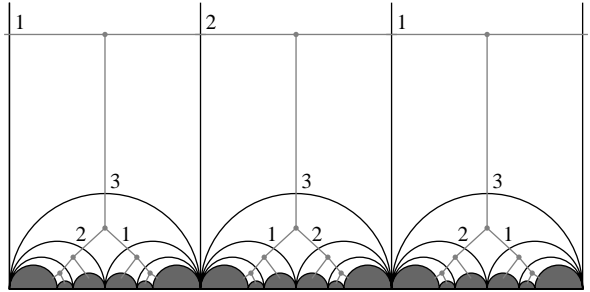}}
\caption{The modular triangulation of the upper half plane, and the dual tree.}
\label{clus1}
\end{figure}

It is well known that in this case the modular group 
is $PSL_2(\Z)$. Consider the classical modular triangulation of the upper
half plane, obtained by reflections of the geodesic triangle with 
vertices at $0, 1, \infty$. 
If we mark the vertices of one of the
triangles by elements of the set $\{1, 2, 3\}$ then there is a unique way to mark 
the vertices of the modular triangulation by the elements of the same
set so that the vertices of each triangle get distinct marks, see Fig \ref{clus1}. 
The set of  vertices of the modular
triangulation is the set of the cusps, identified with $P^1(\Q)$. 
The set of the edges of the modular triangulation  
inherits a decoration by the elements of the set $\{1, 2, 3\}$ such that a vertex 
of each of the modular triangles and the side opposite  to this vertex are
labeled by the same element.  

The modular triangulation is the simplicial complex ${\Bbb S}$. The dual graph of this
triangulation without the vertices is a trivalent tree. 
Its edges inherit labels by the elements of the set $\{1, 2,
3\}$. So it is an $\{1, 2,
3\}$-decorated tree. It is the polyhedral complex ${\Bbb G}_{S}$ 
for the punctured torus $S$: 
an embedded graph as on Fig \ref{clus2} corresponds to a vertex of this 
tree. 
 A flip at an edge of this graph corresponds to the flip at the
 corresponding edge of the 
tree.

\bl
For the cluster ensemble related to the Teichm\"uller space on the punctured torus there are canonical isomorphisms
$$
{\rm Aut}({\Bbb S})= PGL_2(\Z), \quad D = PSL_2(\Z), \quad \Delta = \{e\}, \quad \Gamma = PSL_2(\Z), \quad 
C - C^* = P^1(\Q) .
$$
\el

{\bf Proof}. Recall that $PGL_2(\R)$ acts on $\C -
\R$. This action commutes with the complex conjugation $c: z \to \overline z$
acting on $\C -
\R$. Thus $PGL_2(\R)$ acts on the quotient $(\C -
\R)/c$, which is identified with the upper
half plane.  
The subgroup $PGL_2(\Z)$ 
preserves the modular picture. 
The canonical homomorphism $p: {\rm Aut}({\Bbb S})
 \to {\rm Perm}(I)$ is the projection 
$PGL_2(\Z) \to S_3$ provided by the action on the set $\{1,2,3\}$. 
 The subgroup $PSL_2(\Z)$ of $PGL_2(\Z)$ preserves  the 
function $\varepsilon_{ij}$. So it is the cluster subgroup $D$. 
The cluster subgroup $\Delta$ is trivial in this case. 
This agrees with the fact that the matrix $\varepsilon_{ij}$ has no principal 
$2\times 2$ submatrix from the list (\ref{12.31.03.1}), and thus the 
dual polyhedral complex is reduced to a tree. 
 So the
modular group is $PSL_2(\Z)$. 
The modular triangulation of the upper half plane
coincides with the topological realization of the 
reduced cluster complex $C^*$. One has   $C - C^* = P^1(\Q) $. 
The lemma is proved

%We continue discussion of the cluster nature of the Teichm\"uller space 
%for the punctured torus in Section \ref{Sec6.3}, where we give 
%a cluster description of Thurston's boundary of this Teichm\"uller space.

\vskip 2mm
Below 
we 
explain the cluster nature of the universal Teichm\"uller space, 
from which the case of a hyperbolic surface $S$ can be obtained by 
taking the $\pi_1(S)$-invariants.

\subsection{Cluster nature of  
universal Teichm{\"u}ller spaces and 
the Thompson group}

\begin{definition} \label{8.17.04.3} The {\it universal Teichm{\"u}ller space} ${\cal X}^+$ 
is the space of $PGL_2(\R)$-orbits on the set of maps
\begin{equation} \label{1.29.04.rt} 
\beta: {\Bbb P}^1(\Q) \lra {\Bbb P}^1(\R)
\end{equation} 
respecting the natural cyclic order of  both sets.
\end{definition}
A generalization to an
arbitrary split simple Lie group with trivial center $G$ is in \cite{FG1}. 
The name and relationship with Teichm{\"u}ller spaces 
${\cal T}^+_S$ for surfaces $S$ with punctures is 
explained below. 

Let ${\Bbb S}^1(\R)$ be the set of all rays in $\R^2-\{0,0\}$. 
There is a $2:1$ cover ${\Bbb S}^1(\R) \to {\Bbb P}^1(\R)$. 
Let $s$ be the antipodal involution. It is the unique non-trivial automorphism of this covering.  
Let ${\Bbb S}^1(\Q)$ be the set of its rational points, given by the rays with rational slopes.

\begin{definition} \label{8.17.04.2} Consider the set of all maps
\begin{equation} \label{1.29.04.rtsx} 
\alpha: {\Bbb S}^1(\Q) \lra \R^2 - \{0,0\}\quad \mbox{satisfying the condition $\alpha (s(p)) = - \alpha (p)$}
\end{equation} 
such that composing $\alpha$ with the projection 
$\R^2 - \{0,0\} \lra {\Bbb S}^1(\R)$ we get a map 
$
\overline \alpha: {\Bbb S}^1(\Q) \lra {\Bbb S}^1(\R)
$ 
respecting the natural cyclic order of  both sets.  
The {\it universal decorated Teichm{\"u}ller space} ${\cal A}^+$ 
is the quotient of this set by the natural action of the group $SL_2(\R)$ on it. 
\end{definition}

{\it The Thompson group ${\Bbb T}$}. It is 
the group of all piecewise $PSL_2(\Z)$ automorphisms of
${\Bbb P}^1(\Q)$: for  every 
$g \in {\Bbb T}$ there exists a decomposition of ${\Bbb P}^1(\Q)$
into a union of finite number of segments, which may overlap only at the
ends,  such that the restriction of $g$ to each segment is
given by an element of $PSL_2(\Z)$.

\begin{figure}[ht]
\centerline{\epsfbox{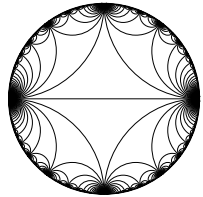}}
\caption{The Farey ideal triangulation of the hyperbolic disc.}
\label{Far}
\end{figure}

Consider the Farey triangulation $T$ of the hyperbolic
plane ${\cal H}$ shown on Fig \ref{clus1}. \ref{Far}
Let $T$ be the dual trivalent tree. We have  canonical identifications 
\begin{equation} \label{1.29.04.6} 
{\Bbb P}^1(\Q) = \Q \cup \infty = 
\{\mbox{vertices of the Farey
  triangulation}\}. 
\end{equation} Let
$$
I_{\cal F} := \{\mbox{edges of the Farey triangulation}\}.
$$
It is identified with the edges of the dual tree $T$. 
Therefore applying  (\ref{7.5.04.11}) to the latter, 
we get  a skew-symmetric function $\varepsilon_{ij}: 
I_{\cal F} \times I_{\cal F} \to
\{0, \pm 1\}$.    
Observe that although $I_{\cal F}$ is an infinite set, 
for any $i \in I_{\cal F}$ the function $I_{\cal F} \to \Z$, $j 
\lms \varepsilon_{ij}$  has
a finite support. It is easy to see that in such situation all the
constructions above work. 
So we get the corresponding seed and the  cluster
ensemble, called {\it Farey cluster ensemble}. 

\begin{theorem} \label{2.22.04.1} Let $({\cal X}_{\cal F}, {\cal
      A}_{\cal F})$ be the 
  Farey cluster ensemble. Then 

i) ${\cal X}_{{\cal F}}(\R_{>0})$ is identified with the universal Teichm{\"u}ller
space ${\cal X}^+$. 

ii) ${\cal A}_{{\cal F}}(\R_{>0})$ is identified with the universal decorated Teichm{\"u}ller
space ${\cal A}^+$. 

iii) The modular group $\Gamma$ 
of the Farey cluster ensemble is the Thompson group. It is isomorphic to the group $\widehat \Gamma$. 
\end{theorem} 

{\bf Proof}. The proof of the parts i) and ii) is very similar to the proofs in the finite genus 
case given in Chapter 11 of \cite{FG1}. Let us outline the proof of i). 
Let us identify once and for all
 the sets ${\Bbb P}^1(\Q)$ in (\ref{1.29.04.6}) and  
(\ref{1.29.04.rt}). 

\begin{lemma} There is a canonical isomorphism 
$% \begin{equation} \label{1.29.04.7} 
\varphi: {\cal X}^+ \stackrel{\sim}{\lra} \R^{I}_{>0}.
$% \end{equation} 
\end{lemma}

{\bf Proof}. It assigns to a map $\beta$ a function $\varphi(\beta)$
on $I$ defined as follows. Let $E$ be an edge of the Farey
triangulation. Denote by $v_1, v_2, v_3, v_4$ the vertices of the
$4$-gon obtained by taking the union of the two triangles sharing
$E$. We assume that the vertices 
follow an orientation of the circle, and the vertex $v_1$ does not belong 
to the edge $E$. Recall the
cross-ratio $r^+$ normalized by $r^+(\infty, -1, 0, x)=x$. Then $$
\varphi(\beta)(E):= r^+(\beta(v_1), \beta(v_2), \beta(v_3),
\beta(v_4))
$$ 
To prove  that every positive valued function  on the set $I$ is
realized we glue one by one triangles of the triangulation 
to the initial one in the hyperbolic disc so that the cross-ratio
corresponding to 
each edge $E$ by the above formula is the value of the given function
at $E$. The lemma is proved. 

It remains to check that flips at the edges of the Farey triangulation are given by the same 
  formulas as the corresponding mutations in the cluster ensemble. This is a 
straightforward check, 
left to the reader. 

ii) Let us define a map ${\cal A}^+ \to \R_{>0}^{I}$. 
Take an edge $E$ of the Farey triangulation. Let $p_1(E)$, $p_2(E)$ be the endpoints 
of the edge $E$, considered as the points of ${\Bbb P}^1(\Q)$. 
The map (\ref{1.29.04.rtsx}) assigns to them vectors $v_1(E)$, $v_2(E)$, 
each well defined up to a sign. The coordinate 
corresponding to the edge $E$ is the absolute value 
of the area of the parallelogram in $\R^2$ generated by these vectors. 
One checks that the exchange relation follows from the Pl\"ucker relation.

iii) Here is another way to look at the Thompson group. 
The Farey triangulation has  a
  distinguished oriented edge, 
connecting $0$ and $\infty$.  
The Thompson group contains the 
following elements, called flips at the edges: Given an edge $E$
of the Farey triangulation $T$, we do a flip at an
edge $E$ 
obtaining a new triangulation $T'$ with a distinguished
oriented edge. This edge is the old one if $E$ is not the distinguished
oriented edge, and it is the flip of the distinguished
oriented edge otherwise. Observe that the ends of the trivalent trees dual to
the triangulations $T$ and $T'$ are
identified, each of them with $P^1(\Q)$. On the other hand, there exists unique isomorphism of the
plane trees $T$ and $T'$ which identifies  their distinguished oriented
edges. It provides a map of the ends of these trees, and hence an
automorphism of $P^1(\Q)$, which is easily seen to be piece-wise
linear. The Thompson group is generated by 
flips at the edges (\cite{I}). It remains to check that the relations in the Thompson group corresponds 
to the standard pentagons in the cluster complex. 
Indeed, these pentagons are exactly the pentagons of the Farey triangulation, 
it is well known \cite{CFP} that they give rise to  relations in the Thompson group, and 
all relations are obtained this way. 

Let us prove that $\widehat \Gamma = \Gamma$. 
First,  $\widehat \Gamma = \Gamma$ is true for the cluster ensemble 
related to a triangulation of an $n$-gon. Indeed, in this case the 
cluster modular complex is nothing else then the Stasheff polytope, 
which is simply connected if $n>4$. This implies that 
the same is true for the modular triangulation. 
 The theorem is proved.

\vskip2mm

The universal decorated Teichm\"uller space ${\cal A}^+$ is 
isomorphic to the one defined by Penner \cite{P1}. 
\vskip3mm
{\it Relation with the Teichm\"uller spaces of  surfaces}.
 Given a torsion free subgroup $\pi \subset 
PSL_2(\Z)$, set $S_{\pi}:= {\cal H}/\pi$. The Teichm{\"u}ller space ${\cal X}^+_{S_{\pi}}$ 
is embedded into ${\cal X}^+$ as the subspace of
$\pi$-invariants:
$$
{\cal X}^+_{S_{\pi}} = ({\cal X}^+)^{\pi}.
$$
Let us define this isomorphism. The Teichm{\"u}ller space 
${\cal X}^+_{S}$ of a surface with punctures $S$ has canonical coordinates 
corresponding to an ideal triangulation of $S$. 
We have a natural triangulation on $S_{\pi}$, the image of the Farey
triangulation under the projection $\pi_{\pi}: {\cal H} \to
S_{\pi}$. So ${\cal X}^+_{S_{\pi}}$ is identified with the
 $\R_{>0}$-valued  functions on $I_{\cal F}/\pi$, i.e. with $\pi$-invariant 
$\R_{>0}$-valued functions on $I_{\cal F}$.

\section{A non-commutative $q$-deformation of the ${\cal X}$-space}

\vskip 3mm
\subsection{Heisenberg groups and quantum tori} Let $\Lambda$ be a lattice 
equipped with  a skew-symmetric bilinear form 
$(\ast, \ast): \Lambda \times \Lambda \to \Z$. 
We associate to this datum a Heisenberg group ${\cal H}_{\Lambda}$. 
It is a central extension 
$$
0 \lra \Z \lra {\cal H}_{\Lambda}\lra \Lambda \lra 0.
$$
The composition law is given by the rule
$$
\{v_1, n_1\} \circ \{v_2, n_2\} = \{v_1+v_2, n_1+n_2 + (v_1, v_2)\}, \qquad v_i \in \Lambda, n_i \in \Z.
$$

\begin{definition} \label{9/16/04/1} Let $\Lambda$ be a 
lattice equipped with a skew-symmetric bilinear form 
$(\ast, \ast): \Lambda \times \Lambda \to \Z$. 
The corresponding quantum torus algebra 
${{\bf T}}_{\Lambda}$ is the group ring of the 
Heisenberg group ${\cal H}_{\Lambda}$. 
\end{definition} 

Let $q$ be the element of the group algebra 
corresponding to the central element $(0, 1) \in {\cal H}_{\Lambda}$.
 Denote by $X_{v}$ the element of the group 
algebra corresponding to the element $(v, 0)\in {\cal H}_{\Lambda}$. Then 
$$
q^{-(v_1, v_2)}X_{v_1} X_{v_2} = X_{v_1+v_2}.
$$
In particular the left hand side is symmetric in $v_1, v_2$. 
There is an involutive antiautomorphism 
$$
\ast: {{\bf T}}_{\Lambda}\lra {{\bf T}}_{\Lambda}, \qquad  
\ast(X_{v}) = X_{v}, ~
\ast(q) = q^{-1}.
$$
\vskip 2mm

Choose a basis $\{e_i\}$ of the lattice $\Lambda$. Set $X_i:= X_{e_i}$. Then
 the algebra  ${{\bf T}}_{\Lambda}$ is identified with the algebra 
of non-commutative polynomials in $\{X_i\}$ over 
the ring $\Z[q, q^{-1}]$ subject to the relations 
\begin{equation} \label{4.28.03.11xe}
q^{-\widehat \varepsilon_{ij}} X_i X_j = 
q^{-\widehat \varepsilon_{ji}} X_jX_i, \qquad \widehat \varepsilon_{ij}:= (e_i, e_j).
\end{equation}
Let us choose an order $e_1, ..., e_n$ of the basis of 
$\Lambda$. Then given a vector $v = \sum_{i=1}^n a_ie_i$ of $\Lambda$, one has
\begin{equation} \label{4.28.03.11xewr}
X_{v} = q^{-\sum_{i<j}a_{i}a_{j}\widehat \varepsilon_{ij}}\prod_{i=1}^nX_i^{a_i}.
\end{equation}
In particular the right hand side does not depend on the choice of the order.

The above construction gives rise to a functor from the category of 
lattices with skew-symmetric forms to the category of non-commutative algebras with an  involutive antiautomorphism $\ast$. 

\vskip 2mm
There is a version of this construction where $q$ is a complex number with 
absolute value $1$ and $\ast$ is a semilinear antiautomorphism 
preserving the generators $X_i$. There is a specialization homomorphism 
of $\ast$--algebras sending the formal variable $q$ to its value.

\paragraph{Center of the quantum torus algebra at roots of unity.} 
Let us start with an example. Let $\Lambda$ be a lattice spanned by 
elements $e_1, e_2$ with $(e_1, e_2) =1$. Then the quantum torus algebra 
${\bf T}_{\Lambda}$ is isomorphic to the algebra of non-commutative Laurent polynomials 
satisfying the relation $q X_1X_2 = q^{-1}X_2X_1$. Suppose now that 
$q$ is an $N$th root of unity. Due to the $q$-binomial formula \cite{Mi} we have
$$
(X_1+X_2)^N = X_1^N+X_2^N. 
$$
Obviously $X_i^N$ are central elements. One easily proves that the center 
of the algebra ${\bf T}_{\Lambda}$ is the algebra of Laurent polynomials in $X_1^N, X_2^N$. 

One can reformulate this using the geometric language as follows. Let $Y_1, Y_2$ 
be the coordinate functions on the two dimensional torus ${\Bbb G}_m \times {\Bbb G}_m$. 
Then there is a natural homomorphism of algebras
$$
{\Bbb F}^*_N: {\cal O}({\Bbb G}_m \times {\Bbb G}_m) \hra {\bf T}_{\Lambda}, \qquad 
Y_i \lms X_i^N. 
$$
Its image is the center of ${\bf T}_{\Lambda}$. we may think about this map 
as of the {\it quantum Frobenious map}
 $$
{\Bbb F}_N: {\rm Spec}({\bf T}_{\Lambda}) \lra {\Bbb G}_m \times {\Bbb G}_m. 
$$
\vskip 2mm
Here is a generalization of this example.
Let $\Lambda_0$ be the kernel of the 
form $(\ast, \ast)$. Then ${\bf T}_{\Lambda_0}$ is in the center  of ${{\bf T}}_{\Lambda}$. 
For generic $q$ the center coincides with ${\bf T}_{\Lambda_0}$. Geometrically, 
the center is generated by the 
preimages of the characters under the Casimir  map
$$
\theta_q:  {\rm Spec}({\bf T}_{\Lambda}) \lra {\rm Spec}({\bf T}_{\Lambda_0}) = {\rm Hom}(\Lambda_0, {\Bbb G}_m)
$$
 In the $q^N=1$ case we have in addition to this the Frobenious map
$$
{\Bbb F}_N: {\rm Spec}({\bf T}_{\Lambda}) \lra {\rm Hom}(\Lambda, {\Bbb G}_m). 
$$
The center in this case is generated by the preimages of the characters under these maps. 

%Let us now specialize $q$ to a non-zero element of a field $F$. 
%For an integer $N$, the sublattice 
%$
%\Lambda_0[N]:= N\Lambda + \Lambda_0 \subset \Lambda  
%$ 
%inherits from $\Lambda$ the skew-linear 
%form, also denoted by $\widehat \varepsilon$. Then if $q^N=1$, then 
%${\bf T}_{\Lambda_0[N]}$ evidently lies in the center of ${{\bf T}}_{\Lambda}$.  

%\begin{lemma} \label{9.6.04.1}
%1. Suppose that $q$ is not a root of unity. Then ${\bf T}_{\Lambda_0}$ is the center  of ${{\bf T}}_{\Lambda}$. 

%2. If $q$ is a primitive $N$-th root of unity, then ${\bf T}_{\Lambda_0[N]}$ is the center of ${\bf T}_{\Lambda}$. 
%\end{lemma} 

%{\bf Proof}. An easy exercise left to the reader. 

\vskip 2mm
 Below we show how to glue quantum tori, 
so that at roots of unity the quantum Frobenious map will be preserved. 

\subsection{The quantum dilogarithm} 
Consider the following 
formal power series, a version of the inverse of the {\it Pochhammer symbol}: 
 \be \la{psi}
{\bf \Psi}_q(x):= \frac{1}{(-qx;q^2)_{\infty}}= \prod_{a=1}^{\infty}(1+q^{2a-1}x)^{-1} = 
\frac{1}{(1+qx)(1+q^3x)(1+q^5x)(1+q^7x)\ldots}.
\ee
It is a {\it $q$-analog of the gamma function}. 
It is characterized, up to a constant, by a  
difference relation
\begin{equation} \label{11.19.06.20}
{\bf \Psi}_q(q^2x) = (1+qx){\bf \Psi}_q(x), \quad 
\mbox{or, equivalently,}\quad  
{\bf \Psi}_q(q^{-2}x) = (1+q^{-1}x)^{-1}{\bf \Psi}_q(x).
\end{equation} 
It is also called the 
{\it $q$-exponential}. The  name is 
justified by the power series expansion
$$
{\bf \Psi}_q(x)= \sum_{n=0}^{\infty}\frac{q^{-\frac{n(n-1)}{2}}x^n}
{(q-q^{-1}) (q^2-q^{-2}) \ldots (q^{n}-q^{-n})}.
$$
It is easily checked by using 
the difference relation. There is a power series expansion of the inverse 
of ${\bf \Psi}_q(x)$:
$$
{\bf \Psi}_q(x)^{-1} = \sum_{n=0}^\infty\frac{q^{n^2}x^n}{(1-q^2)(1-q^4) \ldots (1-q^{2n})}.
$$

The difference relation 
 immediately implies the following property of the $q$-exponential power series:
\begin{equation} \label{1ssw}
{\bf \Psi}_{q^{-1}}(x) = {\bf \Psi}_q(x)^{-1}.
\end{equation} 
Indeed, both parts of the equation satisfy the equivalent difference relations (\ref{11.19.06.20}). 

Formal power series (\ref{psi})
are also known by the name the {\it  quantum dilogarithm power series}. 
To justify the name, recall a version of the 
classical dilogarithm function:
$$
{\rm L}_2(x):= \int_0^x\log(1+t)\frac{dt}{t} = -{\rm Li}_2(-x).
$$
It has a $q$-deformation, called the {\it $q$-dilogarithm power series}, given by 
$$
{\rm L}_2(x;q):= \sum_{n=1}^{\infty}\frac{x^n}{n(q^{n}-q^{-n})}.
$$
One has the identity
$$
\log {\bf \Psi}_q(x)  = {\rm L}_2(x;q).
$$
It is proved  easily by using  the difference relations (\ref{11.19.06.20}) characterizing 
${\bf \Psi}_q(x)$.

\vskip 2mm
The precise relation with the classical dilogarithm is the following. 
If $|q|<1$ the power series ${\bf \Psi}_q(x)^{-1}$ converge, providing 
an analytic function in $x \in \C$. If in addition to this $|x|<1$, the 
$q$-dilogarithm power series also converge. There are asymptotic expansions when $q\to 1^-$:
\begin{equation} \label{sine}
{\rm L}_2(x;q) \sim \frac{{\rm L}_2(x)}{\log q^2}, \qquad 
{\bf \Psi}_q(x) \sim {\rm exp}\Bigl(\frac{{\rm L}_2(x)}{\log q^2}\Bigr). 
\end{equation}

\subsection{The quantum space ${\cal X}_q$} 
According to Definition \ref{esd}, 
a seed ${\mathbf i}$ includes a lattice $\Lambda$ with a skew-symmetric bilinear form $(\ast, \ast)$, and thus 
determines   a quantum torus $\ast$-algebra ${\bf T}_{\Lambda}$, denoted  
${\bf T}_{\mathbf i}$. 
Using the basis $\{e_i\}$ it is described by generators and relations, see (\ref{4.28.03.11xe}). 
Denote by $ {\Bbb T}_{\mathbf i}$ the non-commutative fraction field 
of ${\bf T}_{\mathbf i}$. 

The quantum mutation map $\mu_k^q$ is an isomorphism of skew fields
$$
\mu^q_k: {\Bbb T}_{\bf i'} \lra {\Bbb T}_{\bf i}. 
$$

The simplest way to define it  employs the following fact: 
The algebras  
${\bf T}_{\mathbf i}$ for the seeds ${\bf i}$ related by seed 
cluster transformations are canonically isomorphic, 
since each of them is identified with the algebra ${\bf T}_\Lambda$,

\begin{definition} \label{11.18.06.111} 
The mutation homomorphism $\mu^q_k: {\Bbb T}_{\bf i'} \lra {\Bbb T}_{\bf i}$ 
is  the conjugation by 
${\bf \Psi}_{q_k}(X_k)$, where $X_k = X_{e_k}$ is a basis element for the seed ${\bf i}$: 
$$
\mu_k^{q}:= 
{\rm Ad}_{{\bf \Psi}_{q_k}(X_k)}, \qquad q_k := q^{1/d_k}.
$$ 
\end{definition} 
In other words, the map $\mu^q_k$ is defined as the unique map making the following diagram commutative.
Here the vertical maps are the canonical isomorphisms. 
$$
\begin{array}{ccc}
{\Bbb T}_{\bf i'}& \stackrel{\mu^q_k}{\lra}&{\Bbb T}_{\bf i}\\
\sim \downarrow && \downarrow \sim\\
{\Bbb T}_{\Lambda} & \stackrel{{\rm Ad}_{{\bf \Psi}_{q_k}(X_k)}}{\lra} &{\Bbb T}_{\Lambda} 
\end{array}
$$

Although  ${\bf \Psi}_{q_k}(X_k)$ is not a rational function, 
we show in Lemma \ref{11.18.06.3} that $\mu_k^{q}$
is a rational map. 

\begin{lemma} \label{4.28.03.1}
The map $\mu^q_{k}$ is a homomorphism of $\ast$-algebras.  
\end{lemma}

{\bf Proof}. The map 
$\mu^{q}_k$ is given by the conjugation. So is a homomorphism of algebras. It commutes with the involution $\ast$ thanks to (\ref{1ssw}). Indeed, since $\ast X_i = X_i$, we have 
$$
\ast \Bigl({\bf \Psi}_{q}(X_k)X_i{\bf \Psi}_{q}(X_k)^{-1}\Bigr) = 
{\bf \Psi}_{q^{-1}}(X_k)^{-1}X_i{\bf \Psi}_{q^{-1}}(X_k) \stackrel{(\ref{1ssw})}{= } 
{\bf \Psi}_{q}(X_k)X_i{\bf \Psi}_{q}(X_k)^{-1}. 
$$

\paragraph{Decomposition of quantum mutations.} Although the algebras ${\bf T}_{\mathbf i'}$ and 
${\bf T}_{\mathbf i}$ are canonically isomorphic, 
they are equipped with  different sets of 
the generators -- the cluster coordinates -- $\{X_{e'_i}\}$ and $\{X_{e_i}\}$. 
Let us write the mutation map in the 
cluster coordinates. 
Then we have 
$$
\mu^q_k = \mu^\sharp_k\circ \mu'_k, \qquad \mu^\sharp_k:= 
{\rm Ad}_{{\bf \Psi}_{q_k}(X_k)}: {\bf T}_{\mathbf i}\lra {\bf T}_{\mathbf i}, 
$$
where $\mu'_k$ is a map which tells how the coordinates related to the basis $\{e'_i\}$ 
are related to the ones related to the basis $\{e_i\}$. It is given in the cluster coordinates as follows: 
\be \la{muk}
\mu'_k: {\bf T}_{\bf i'} \lra {\bf T}_{\bf i}, \quad 
 X_{e_i'} \lms X_{e_i'} = X_{e_i + [\varepsilon_{ik}]_+e_k} = 
q^{-\widehat\varepsilon_{ik}[\varepsilon_{ik}]_+} 
X_{e_i}X^{[\varepsilon_{ik}]_+}_{e_k}. 
\ee
 So although it is the identity map after the canonical 
identification of algebras ${\bf T}_{\mathbf i'}$ and 
${\bf T}_{\mathbf i}$, it looks as a non-trivial map when written 
in the cluster coordinates.

\paragraph{An explicit computation of the automorphism $\mu_k^{\sharp}$.} 

\begin{lemma} \label{11.18.06.3} The automorphism $\mu_k^{\sharp}$ is 
given on the generators by the formulas
\begin{equation} \label{f3*}
X_i \lms X^{\sharp}_{i} := \left\{\begin{array}{lll} %X_k& \mbox{ if } & i=k, \\
X_i (1+q_kX_k)(1+q^3_kX_k)
  \ldots (1+q_k^{2|\varepsilon_{ik}|-1}X_k)& 
\mbox{ if } &  \varepsilon_{ik}\leq 0, 
\\
    X_i \left((1+q^{-1}_kX_k)(1+q_k^{-3}X_k)
  \ldots (1+q_k^{1-2|\varepsilon_{ik}|}X_k)\right)^{-1}& \mbox{ if } &  \varepsilon_{ik}\geq 0. \\
\end{array} \right.
\end{equation} 
\end{lemma}

{\bf Proof}. 
For any formal power series $\varphi(x)$ the relation 
$
q^{- \widehat \varepsilon_{ki}}X_kX_i = q^{- \widehat \varepsilon_{ik}}X_iX_k
$ 
implies 
\begin{equation} \label{hghghgw3} 
\varphi(X_k)X_i = 
X_i\varphi(q^{-2\widehat \varepsilon_{ik}}X_k).
\end{equation}
The  difference equation (\ref{11.19.06.20})  implies that  the 
formula  (\ref{f3*}) can be rewritten as
\begin{equation} \label{hghghgw} 
X^\sharp_i =
X_i \cdot {\bf \Psi}_{q_k}(q_k^{-2\varepsilon_{ik}}X_k)
{\bf \Psi}_{q_k}(X_k)^{-1} .
\end{equation}
Using   (\ref{hghghgw3}) and 
$q_k^{-2\varepsilon_{ik}} =q^{-2\widehat \varepsilon_{ik}}$, 
we get 
$$
{\bf \Psi}_{q_k}(X_k)X_i{\bf \Psi}_{q_k}(X_k)^{-1} =
X_i {\bf \Psi}_{q_k}(q^{-2\widehat \varepsilon_{ik}}X_k) 
{\bf \Psi}_{q_k}(X_k)^{-1}  \stackrel{(\ref{hghghgw})}{=} X^\sharp_i.
$$ 
 The lemma is proved.

\vskip 3mm
Let $a\geq 0$ be an integer
and  
$$
G_a(q; X) := {\Psi}_{q}(q^{2a}X){\Psi}_{q}(X)^{-1} = 
\left\{ \begin{array}{ll}\prod_{i=1}^a(1+q^{2i-1}X) & a>0\\
1 & a=0.
\end{array}\right.
$$ 
Let $\mu_k := {\mathbf i} \to {\mathbf i}'$ be a mutation.
The following Lemma follows easily from Lemma \ref{11.18.06.3}. 
\bl \la{QM}
The {\it quantum mutation} homomorphism 
$$
\mu^q_{k}: 
{\Bbb T}_{{\mathbf i}'} \lra {\Bbb T}_{\mathbf i}
$$
is given in the cluster coordinates by the formula 
$$
\mu^q_{k}: X_i' \lms  \left\{ \begin{array}{ll} X_i F_{ik}(q; X_k)& \mbox{ if $k\not =i$}\\
X^{-1}_i & \mbox{ if $k=i$,}\end{array}\right.
\qquad
F_{ik}(q; X) = \left\{ \begin{array}{ll}
G_{|\varepsilon_{ik}|}(q_k; X) 
& \mbox{ if $\varepsilon_{ik} \leq 0$}\\
G_{|\varepsilon_{ik}|}(q_k; X^{-1})^{-1}& 
\mbox{ if $\varepsilon_{ik}\geq 0$}.
 \end{array}\right. 
$$
\el
In particular, it implies 

\bc Setting $q=1$ we recover 
the ${\cal X}$-mutation formulae. 
\ec

{\it The Poisson structure on ${\cal X}$}. 
The quasiclassical limit of the non commutative space  
 ${\cal X}_q$ is described by a Poisson structure on the ${\cal X}$-space. 
This Poisson structure in any cluster coordinate system $\{X_i\}$ 
is 
given by the formula 
$%\begin{equation} \label{4.30.03.2i}
\{X_{i}, X_{j}\} = 2\widehat \varepsilon_{ij} X_{i} X_{j}.
$ %\end{equation}
Lemma \ref{4.28.03.1} implies that it 
is independent of the choice of coordinate system.

\begin{lemma} \label{5.2.03.101} 
We have $(\mu^q_k)^2 = {\rm Id}$ for quantum mutations.
\end{lemma}

{\bf Proof}.  Suppose that $\varepsilon_{ik}=a>0$. Then, 
using (\ref{hghghgw3}) and difference equation (\ref{11.19.06.20}),  
we have 
$$
{\rm Ad}_{\Psi_{q_k}(X_k^{-1})}{\rm Ad}_{\Psi_{q_k}(X_k)} X_i
 = X_iG_a(q_k; X^{-1}_k) G_a(q^{-1}_k, X_k) = q_k^{a^2}X_iX^{-a}_k.
$$
Being composed with (\ref{muk}), this gives the identity map. 
The other case is reduced to this one. The lemma is proved.

\begin{proposition} \label{11.9.03.21}
The collection of the quantum tori 
${\bf T}_{\mathbf i}$
and the quantum mutation maps $\mu^q_{k}$ 
provide a functor 
$
\widehat {\cal G} \lra {\rm QPos}^*.
$ 
\end{proposition}
The quantum space ${\cal X}_q$ is understood as this functor. 

\vskip 3mm
{\bf Remark}. We expect to have 
a functor ${\cal G} \lra {\rm QPos}^*$. The problem is that we do not know 
the relations in the groupoid ${\cal G}$ explicitly. 
\vskip 3mm 

{\bf Proof}. We have to check that the composition of maps corresponding to the
boundary of any standard $(h+2)$-gon equals to the identity. It  can be checked 
by a calculation. We present its crucial step as 
Lemma \ref{5.22.03.2} below. 
The proposition is proved. 
\vskip 3mm

{\bf Remark}. 
Another proof follows from Lemma 2.22 in \cite{FG2} plus the trick used 
in the proof of Proposition \ref{FTFTFT}: one embeds the seed 
${\bf i}$ in a bigger seed ${\bf i'}$ with 
${\rm det}~\varepsilon'_{ij} \not = 0$, and observes that  cluster 
transformations (\ref{K10}) remain trivial on the classical level after 
extension of the seed, so by Lemma 2.22 in \cite{FG2} 
they are trivial on the quantum level.

\paragraph{Examples of quantum relations.} Consider a sequence of mutations 
at the vetrices  $i, j, i, j, i, ...$. We picture it by a 
polygon, whose vertices  match the mutations. 
The seeds are the sides of the polygon, and the ${\cal X}$-coordinates 
for a given seed ${\bf i}$ are assigned to the flags 
(a vertex of the side, the side). The ${\cal X}$-coordinates 
assigned to the flags sharing a vertex are opposite to each other. 
Below we calculate the sequence of the ${\cal X}$-coordinates 
assigned to the flags oriented the same way,  clockwise. 
They determine the set of all ${\cal X}$-coordinates 

The  ${\cal X}$-coordinates on the set of all clockwise oriented flags 
are obtained from the initial
${\cal X}$-coordinates $x_1, x_2$ by the following
inductive procedure. 

\vskip 2mm
{\bf Classical case}. Let $F$ be a field, $x_1, x_2 \in F^*$ and $b,c $ are non-negative integers. 
Consider the recursion
\begin{equation} \label{5.22.03.1a}
x_{m-1}x_{m+1} =  \left\{ \begin{array}{ll}(1+x_{m})^b & \mbox{ $m$: even }\\
(1+x_{m})^c
& \mbox{ $m$: odd }.\end{array}\right.
\end{equation}
According to Chapters 2 and 6 of [FZ1], this sequence is periodic if and only if 
 the Cartan matrix $\left (\matrix {2& -b\cr -c&2\cr}\right )$ 
 or its transpose is of finite type, 
 i.e. 
 $b=c=0$ or $1 \leq |bc| \leq 3$. 
 Therefore, up to a shift $x_i \lms x_{i+1}$,
 there are only 
 four periodic sequences, corresponding to the root systems $A_1\times A_1, A_2, B_2, G_2$. 
The period is $h+2$, where $h$ 
 is the Coxeter 
 number of the root system. 

\vskip 2mm
{\bf Quantum case}. Let $(\varepsilon_{ij})=
\left (\matrix{0& -1\cr c&0\cr}\right )$. 
%Then  $( d_1,  d_2) =  (c, 1)$. 
So the commutation relations are 
$q^{-2c}X_iX_{i+1} = q^{2c}X_{i+1}X_i$. We have:
\begin{equation} \label{5.22.03.1}
X_{m-1}X_{m+1} =  \left\{ \begin{array}{ll}(1+q^cX_{m}) & \mbox{ $m$: even }\\
(1+qX_{m})\cdot(1+q^3X_{m})\cdot ... \cdot (1+q^{2c-1}X_{m}) 
& \mbox{ $m$: odd }\end{array}\right.
\end{equation}

Let $h$ be the Coxeter number for the 
Cartan matrix $\left (\matrix{2& -1\cr -c&2\cr}\right )$ of finite type, i.e. 
 $c=0, 1,2,3$. 
So $h=2$  for $c=0$; $h=3$  for $c=1$; $h=4$ for $c=2$; and $h=6$ for $c=3$. 
\begin{lemma}\label{5.22.03.2} 
For any integer $m$
one has $X_{m+h+2} = X_m$. 
\end{lemma}

{\bf Proof}. Compute the elements $X_m$:

{\bf Type} $A_1 \times A_1$. Then $b=c=0, h=2$ and 
$
X_3 = X_1^{-1}, \quad X_4 = X_2^{-1}
$.

{\bf Type} $A_2$. Then $b=c=1, h=3$ and 
$$
X_3 = X_1^{-1}(1+qX_2), \quad X_4 = (X_1X_2)^{-1}\Bigl(X_1+q(1+qX_2)\Bigr)
, \quad X_5 =X_2^{-1}(1+q^{-1}X_1),  \quad X_6 = X_1.
$$

{\bf Type} $B_2$.  Then $b=1, c=2, h=4$ and 
$$
X_1X_3 = 1+q^2X_2, \qquad 
X_1^2X_2X_4 = \Bigl(X_1 + q(1+q^{6}X_2)\Bigr)\Bigl(X_1 + q^3(1+q^{2}X_2)\Bigr),
$$
$$
X_1X_2X_5 = q^2\Bigl((1+q^{-1}X_1)(1+q^{-3}X_1) + q^2X_2\Bigr), \quad X_2X_6 
= (1+q^{-1}X_1)(1+q^{-3}X_1), \quad X_7 = X_1.
$$

{\bf Type} $G_2$.  Then $b=1, c=3, h=6$ and 
$$
X_1X_3 = 1+q^3X_2, \qquad X^3_1X_2X_4 = \Bigl(X_1 + q(1+q^{15}X_2)\Bigr)
\Bigl(X_1 + q^3(1+q^{9}X_2)\Bigr)\Bigl(X_1 + q^5(1+q^{3}X_2)\Bigr),
$$
$$
X_1^2X_2X_5 = (1+q^{-1}X_1)^3 + (q^3X_2)^2 + (1+q^6)q^3X_2 +  
3q^{-1}X_1q^3X_2, \quad \ldots \quad , 
$$
$$
X_2X_8 = (1+q^{-1}X_1)(1+q^{-3}X_1)(1+q^{-5}X_1), \quad X_9 = X_1.
$$
%The proofs of these 
%identities follow the same pattern, but are  rather lengthy. 
% The method is clear from the two examples given below.

%{\it Proof of the formula for $X_5$ in the $A_2$ case}: 
%$$
%X_2X_5= X_2X_3^{-1}(1+qX_4) = X_2(X_1X_3)^{-1}X_1(1+qX_4) = 
%q^{-2}(1+qX_2)^{-1} X_1X_2(1+qX_4) = 
%$$
%$$
%q^{-2}(1+qX_2)^{-1} \Bigl(X_1X_2 +  qX_1(1+qX_3)\Bigr)
% = q^{-2}(1+qX_2)^{-1}\Bigl(q^2X_2X_1 +qX_1 +q^2(1+qX_2)\Bigr) = q^{-1}X_1+1
%$$

%{\it Proof of the formula for $X_5$ in the $B_2$ case}: 
%$$
%X_1X_2X_5 = X_1X_2X_3^{-1}(1+q^2X_4) = X_1X_2(X_1X_3)^{-1}X_1(1+q^2X_4) = 
%X_1X_2(1+q^2X_2)^{-1}X_1(1+q^2X_4) = 
%$$
%$$
%q^{-4}(1+q^6X_2)^{-1}X_1^2X_2(1+q^2X_4) = 
%q^2\Bigl((1+q^{-1}X_1)(1+q^{-3}X_1) + q^2X_2\Bigr)
%$$

%{\it Proof of the formula for $X_4$ in the $G_2$ case}. 
%$$
%X_1^3X_2X_4 = X_1^3(1+qX_3)(1+q^3X_3)(1+q^5X_3) = 
%\Bigl(X_1 + q(1+q^{15}X_2)\Bigr)
%\Bigl(X_1 + q^3(1+q^{9}X_2)\Bigr)\Bigl(X_1 + q^5(1+q^{3}X_2)\Bigr)
%$$

%The Theorem \ref{11.9.03.21} is proved. 

{\bf Remark}. This way we get 
just the half of all ${\cal X}$-coordinates. They are related to the 
ones in the Example as follows: we get only the ${\cal X}$-coordinates 
assigned to the mutating vertices before the mutations;
the initial coordinates are $x_1: = y^{-1}, x_2: = x$.

\subsection{The quantum Frobenius map}
\paragraph{Center of the quantum space when $q$ is not a root of unity.} 
Let $\alpha \in {\rm Ker}_L[\ast, \ast]$.
Then  the element $X_{\alpha}$ %, see (\ref{9.11.04.2}), 
is in  the center of the quantum torus ${\rm T}^q_{\mathbf i}$. 
The torus $H_{\cal X}$
can be treated as a commutative quantum positive space, see (\ref{thet}). 
Recall the character $\chi_{\alpha}$ of the 
torus $H_{\cal X}$ corresponding to $\alpha$. 

\begin{lemma} \label{1.09.04.100a}
There exists a unique map of quantum positive spaces
$\theta_q: {\cal X}_q \to H_{\cal X}$ such that for any seed we have
$\theta_q^*\chi_{\alpha}= X_{\alpha}$, where  $\alpha \in {\rm Ker}_L[\ast, \ast]$. 
\end{lemma}

{\bf Proof}. Since the element $X_\alpha$ of the quantum torus 
${\bf T}_{\Lambda}$ corresponding to a  vector 
$\alpha \in {\rm Ker}_L[\ast, \ast]$ lies in the center, conjugation by 
$\Psi_q(X_k)$ acts on them as the identity. The Lemma follows. 
\vskip 3mm

\paragraph{Center of the quantum space when $q$ is a root of unity.} 
When $q$ is a root of unity, the quantum space has a much larger center, 
which we are going to describe now. 
Let 
  $\widehat \varepsilon_{ij}\in \Z$. 
For a seed ${\bf i}$,  
denote by $\{Y_i\}$ (respectively $\{X_i\}$) the corresponding cluster coordinates 
on  ${\cal X}$ (respectively on ${\cal X}_{q}$). 
If $q^N=1$ then $X_i^N$ are in the center of the quantum torus algebra 
${\bf T}_{\mathbf i}$. 

\begin{theorem} \label{QFROB} Let $q=\zeta_{N}$ be a
primitive $N$-th root of unity. Let us assume that  $q^{d_k}$ is a primitive 
$N$-th root of
  unity, and for every seed ${\mathbf i}$ the corresponding function 
$\varepsilon_{ij}$
  satisfies 
$(2 \varepsilon_{ij}, N) =1$. Then there exists a map of positive spaces, 
called the {\it quantum Frobenius map}, 
$$
{\Bbb F}_N: {\cal X}_{q} \lra {\cal X}\quad \mbox{such that} \quad 
 {\Bbb F}^*_N(Y_i):= X_i^N
$$
in any cluster coordinate system.
\end{theorem}

{\bf Proof}. To check that the quantum Frobenius map commutes with 
a mutation $\mu_k: {\mathbf i} \to {\mathbf i}'$, 
%This is obvious for the monomial part $\mu_k'$ of the mutation. 
%This will require only the conditions on $d_k$ and
%$\varepsilon_{ij}$ at the seed ${\mathbf i}$. 
%Consider the following diagram
%$$
%\begin{array}{ccc}
%Y_i'& \stackrel{\mu_{e_k}}{\lra} & \left\{ \begin{array}{ll}
%Y_i(1+Y_k)^{-\varepsilon_{ik}}
%& \mbox{ if $ i
%\not = k, \quad \varepsilon_{ik} \leq 0$}\\
%Y_i(1+Y^{-1}_k)^{-\varepsilon_{ik}}& \mbox{ if $ i
%\not = k, \quad \varepsilon_{ik} \geq 0$}
% \end{array}\right.\\
%&&\\
%{\Bbb F}^*_N\downarrow && \downarrow {\Bbb F}^*_N\\
%&&\\
%(X_i')^N & \stackrel{\mu_{e_k}}{\lra}& \left\{ \begin{array}{ll}
%X_i^N(1+X_k^N)^{-\varepsilon_{ik}}
%& \mbox{ if $ i
%\not = k, \quad \varepsilon_{ik} \leq 0$}\\
%X_i^N(1+X_k^{-N})^{-\varepsilon_{ik}}& \mbox{ if $ i
%\not = k, \quad \varepsilon_{ik} \geq 0$}.
% \end{array}\right.\\
% \end{array}
%$$
%We need to prove that this diagram is commutative.
%The case $\varepsilon_{ik}=0$ is trivial.
it is sufficient to check that it commutes with the conjugation by $\Psi_q(X_k)$. 
Here we consider the generic $q$, and only after the conjugation 
specialize 
$q$ to a root of unity. 

Let us assume  that  $\varepsilon_{ik} = -a\leq  0$. Then we have to show that, 
specialising $q_k=1$ in 
$$
\Psi_{q_k}(q_k^{2aN}X_k)\Psi_q(X_k)^{-1},
$$
we get $(1+X_k^N)^a$.
%\begin{equation} \label{12.22.03.1}
%(X_iG_a(q_k; X_k))^N = X_i^N (1+X_k^N)^a.
%\end{equation}
%Observe that
%$
%G_a(q_k; X_k)X_i = X_i G_a(q_k; q^{-2\widehat \varepsilon_{ik}}X_k)
%$.
%Therefore the left hand side of (\ref{12.22.03.1}) equals to
%$$
%X_i^N \prod_{b=0}^{N-1}G_a(q_k; q_k^{-2b \varepsilon_{ik}}X_k). 
%$$
The statement is equivalent to the identity
\begin{equation} \label{12.22.03.2}
\prod_{b=0}^{N-1}G_a(q_k; q_k^{2b a}X_k) = (1+X_k^N)^a.
\end{equation}
Notice that $q_k$ is a primitive $N$-th root of unity and,
since $(2 a, N) =1$, the set $\{-2ab\}$, when $b\in \{1, ..., N-1\}$,  
consists of all residues modulo $N$ except zero.
Thus each factor of the product
$$
G_a(q_k; X_k) = \prod_{i=1}^a(1+q_k^{2i-1}X_k)
$$
contributes $(1+X_k^N)$ thanks to the formula
$
\prod_{c=0}^{N-1}(1+q_k^{2i-1 +c}Z) = 1+Z^N.
$

The argument in the case $\varepsilon_{ik}=a>0$ is similar. 
In fact it can be reduced to the previous case using $(\mu^q_k)^2 ={\rm Id}$ 
and $\varepsilon'_{ik}= - \varepsilon_{ik}$. 
%Now let us assume that $\varepsilon_{ik}=a>0$. We have to compute
%$
%(X_iG_a(q_k^{-1}; X_k)^{-1})^N.
%$  
%Observe that
%$$
%G_a(q_k; X_k^{-1})^{-1}X_i = X_i G_a(q_k; q^{2\widehat \varepsilon_{ik}}
%X_k^{-1})^{-1}.
%$$
%Thus we get
%$
%X_i^N \prod_{b=0}^{N-1} G_a(q_k; q^{2b\widehat \varepsilon_{ik}}X_k^{-1})^{-1}.
%$ 
%Then the argument is the same as above. 
The theorem is proved.

\vskip 2mm 
{\bf Examples}. a) Let $\widehat S$ be a marked hyperbolic surface. Then the
pair of moduli spaces $({\cal X}_{PGL_2, \widehat S}, {\cal A}_{SL_2, \widehat
  S})$ 
has a cluster ensemble structure 
with $\varepsilon_{ij} \in \{0, \pm 1, \pm 2\}$
(\cite{FG1}, Chapter 10; \cite{FG}). Thus it satisfies the assumptions of the theorem for
any odd $N$. 

b) A cluster ensemble of finite type satisfies the assumptions of the
theorem for any odd $N$ in all cases except $G_2$, 
where the condition is $(N, 6) =1$. 

\vskip 2mm 
{\bf Remark}. Sometimes it makes sense to restrict the functor defining the 
space  ${\cal X}_q$ to a subgroupoid $\widehat {\cal G}'$ 
of $\widehat {\cal G}$, restricting therefore the set
of values of the exchange function. For example 
for the pair of moduli spaces $({\cal X}_{PGL_m, S}, {\cal A}_{SL_m, 
  S})$ one may consider only those mutations which were introduced 
in Chapter 10 of \cite{FG1} to decompose flips. Then the restricted cluster
function takes values in $\{0, \pm 1, \pm 2\}$, $d_k=1$, 
and the fundamental group
of the restricted groupoid $\widehat {\cal G}'$ contains 
the classical modular group of $S$.  So the quantum Frobenius map in
this case commutes with the action of the classical modular group. 
If the 
modular group $\Gamma$ is finitely generated, 
we can always restrict to a subgroupoid $\widehat {\cal G}'$ 
of $\widehat {\cal G}$ which has the same fundamental group 
and a bounded set of values $|\varepsilon_{ij}|$.

\section{Duality and canonical pairings: conjectures} 
Below we denote by 
${\cal A}$ and ${\cal A}^{\vee}$ the positive spaces 
${\cal A}_{|{\bf i}|}$ and ${\cal A}_{|{\bf i^{\vee}}|}$, and similarly for the 
${\cal X}$-spaces. 

In this Section we show how to extend  
to  cluster ensembles the philosophy of duality 
between the ${\cal X}$ and ${\cal A}$ positive spaces developed in 
\cite{FG1} in the context  of the two  moduli spaces 
related to a split semisimple group $G$ and a surface $S$. 
We suggest that there exist several types of closely related 
canonical pairings/maps  
between the positive spaces ${\cal X}$ and 
${\cal A}^{\vee}$. An example provided 
by the cluster 
ensemble related to the classical Teichm{\"u}ller theory was  elaborated 
in Chapter 12 of \cite{FG1}. It was extended to the pair of Teichm\"uller spaces related to a surface $S$ with $m>0$ distinguished points on the boundary in \cite{FG}. In particular, when $S$ is a disc with $m$ marked points on the boundary, we cover
 the case of the cluster ensemble of finite type $A_m$. 
The canonical map ${\Bbb I}_{\cal X}$ for cluster ensembles of an arbitrary finite type 
is defined in Section 4.6. Other examples  can be obtained using 
the work \cite{SZ} 
on the rank two finite and affine cluster algebras.  

Our main conjectures are Conjecture \ref{10.10.03.10*} and its quantum version, Conjecture \ref{10.10.03.10fgf}. 

Conjecture \ref{10.10.03.10} is a variation on the theme of Conjecture \ref{10.10.03.10*}. 
We show that, under some assumptions, the latter can be deduced from the former.

\paragraph{Background.}  
  Let $L$ be a set. Denote by $\Z_+\{L\}$ the abelian 
semigroup generated by $L$. Its elements are expressions 
$\sum_i n_i \{l_i\}$ where $n_i \geq 0$, the sum is finite, and 
$\{l_i\}$ is the generator corresponding to  $l_i \in L$. 
Similarly  $\Z\{L\}$ is the abelian group generated by $L$. 

Let ${\cal X}$ be a positive space. Observe that the ring of regular functions on a split torus $H$ is the ring of Laurent polynomials in characters of $H$. 
Recall  (Chapter 1.1 and 
 \cite{FG1}, Section 4.3) 
that a {\it universally Laurent polynomial} on ${\cal X}$ 
is a regular function  on one of
the coordinate tori  $H_{\alpha}$ defining ${\cal X}$
 whose restriction to any other  
coordinate torus  $H_{\beta}$  is a 
regular function there.  
${\Bbb L}({\cal X})$ denotes 
the ring of all universally Laurent polynomials, and ${\Bbb L}_+({\cal
  X})$ 
the semiring of {\it  universally positive Laurent polynomials} 
obtained by imposing the positivity condition on  
coefficients of universally Laurent polynomials.    Let  
${\bf E}({\cal X})$ be the set of extremal elements, that is 
 universally positive Laurent polynomials which can not be decomposed into a sum of 
two non zero 
 universally positive Laurent polynomials with positive coefficients. 

We use the notation  
${\cal X}^+ := {\cal X}(\R_{>0})$. 
Recall that for a given seed ${\mathbf i}$, there are canonical coordinates 
$\{X_i\}$ and $\{A_i\}$ on the ${\cal X}$ and ${\cal A}$ 
spaces. By the very definition, their restrictions to   
${\cal X}^+$ and ${\cal A}^+$ are positive, 
so we have the corresponding logarithmic coordinates $x_i:= \log X_i$ 
and $a_i:= \log A_i$. The coordinates on  the tropicalizations 
${\cal X}({\Bbb A}^t)$ and ${\cal A}({\Bbb A}^t)$ are   
also denoted by $x_i$ and $a_i$.

A convex function on a lattice $L$ is a function $F(l)$ such that 
$F(l_1 +l_2) \leq F(l_1) + F(l_2)$. Let ${\cal X}$ be a positive space. 
A convex function on ${\cal X}({\Bbb A}^t)$ or ${\cal X}(\Q)$ is a function 
which is convex in each of the coordinate systems from the defining atlas 
of ${\cal X}$. 
 \vskip 3mm
\subsection{Canonical maps in the classical setting} 
We consider a partial order on monomials $X^a = \prod X_i^{a_i}$ such that 
$X^a \geq X^b$ if $a_i \geq b_i$ for all $i \in I$. We say that 
$X^a$ is the highest term of a Laurent polynomial $F$ if 
$X^a$ is bigger (i.e. $\geq $)   then any other monomial in $F$. 
Recall the map $p: {\cal A} \to {\cal X}$. 

\begin{conjecture} \label{10.10.03.10*}
There exist   $\Gamma$-equivariant isomorphisms of sets 
\begin{equation} \label{prod1a} 
{\cal A}(\Z^t) = {\bf E}({\cal X}^{\vee}), 
\qquad {\cal X}(\Z^t) = {\bf E}({\cal A}^{\vee}). 
\end{equation}
These isomorphisms give rise to $\Gamma$-equivariant isomorphisms
\begin{equation} \label{10/9/03/1a}
{\Bbb I}_{\cal A}: \Z_+\{{\cal A}(\Z^t)\} \stackrel{\sim}{\lra} 
{\Bbb L}_+({\cal X}^{\vee}), \qquad 
{\Bbb I}_{\cal X}: \Z_+\{{\cal X}(\Z^t)\} \stackrel{\sim}{\lra}  
{\Bbb L}_+({\cal A}^{\vee}), 
\end{equation}

These maps have the following properties:
\begin{enumerate}
\item \label{high} Let $(a_1, ..., a_n)$ be the coordinates of $l \in 
{\cal A}(\Z^t)$. Then the highest  term of 
${\Bbb I}_{\cal A}(l)$ is $\prod_i X_i^{a_i}$.

\item \label{posx} If in a certain cluster coordinate system 
the coordinates $(x_1, ..., x_n)$ of an element
  $l \in {\cal X}(\Z^t)$ are non negative numbers,  then 
$$
{\Bbb I}_{\cal X}(l) = \prod_{i \in I}A_i^{x_i}.
$$

\item\label{symmetry-addalg} 
Let 
$l \in {\cal A}(\Z^t)$ and $m$ is a point of ${\cal A}^{\vee}$. Then 
$$
{\Bbb I}_{\cal A}(l)(p(m)) = {\Bbb I}_{\cal X}(p(l))(m),
$$ 
where ${\Bbb I}_{\cal A}(l)(p(m))$ means the value of the function 
${\Bbb I}_{\cal A}(l)$ on $p(m)$

\item Extending the coefficients from $\Z_+$ to $\Z$, we arrive at isomorphisms
\begin{equation} \label{10/9/03/2a}
{\Bbb I}_{\cal A}: \Z\{{\cal A}(\Z^t)\} \stackrel{\sim}{\lra}  
{\Bbb L}({\cal X}^{\vee}), \qquad 
{\Bbb I}_{\cal X}: \Z\{{\cal X}(\Z^t)\} 
\stackrel{\sim}{\lra}  {\Bbb L}({\cal A}^{\vee}), 
\end{equation}
\end {enumerate}
\end{conjecture}

%{\bf Remark}. Conjecture \ref{10.10.03.10*} should contain as a special case 
%canonical bases in the universal enveloping algebras constructed by Lusztig and Kashiwara. 
%It is known, however, that there are some bases very closely related to the canonical 
% bases. 
%A similar situation should happen in the general cluster set-up. 
%Precisely, the set ${\bf E}$ of 
%extremal points may not be a simplicial set in general, i.e. it may have faces 
%which are not simplices. 
%Nevertheless we expect in the most general case that there  exists 
%a $\Z$-linear basis parametrised by the tropical points of the dual space. 
%However it may not be a $\Z_{>0}$-linear basis.  

%\vskip 2mm

The isomorphisms (\ref{prod1a}) imply that one should have
\begin{equation} \label{prod1} {\Bbb I}_*(l_1) {\Bbb I}_*(l_2)=\sum_l c_*(l_1, l_2;l)
 {\Bbb I}_*(l), 
\end{equation} 
where the  coefficients $c_*(l_1, l_2;l)$
are positive integers and the sum is finite. 
Here $*$ stands for either ${\cal A}$ or ${\cal X}$. 
It follows from part \ref{high} of Conjecture \ref{10.10.03.10*} 
and (\ref{prod1}) that $c_*(l_1, l_2; l_1+ l_2)=1$. 

In addition, the canonical maps are expected to 
satisfy the following additional properties. 

\vskip 3mm
i) {\it Convexity conjecture}. The structural constants $c_*(l_1, l_2;l)$ can be viewed as  maps
$$
c_{\cal A}: {\cal A}(\Z^t)^3 \lra \Z, \quad c_{\cal X}: {\cal X}(\Z^t)^3 \lra \Z.
$$
\begin{conjecture} \label{convexitya}
In any cluster coordinate system, 
the supports of the functions $c_{\cal A}$ and $c_{\cal X}$ are convex
polytopes. 
\end{conjecture}

ii) {\it Frobenius Conjecture}. We conjecture that  
in every cluster coordinate system 
$\{X_i\}$, for every prime $p$  one should have the congruence
\begin{equation} \label{5.7.04.1}
{\Bbb I}_{\cal A}(p\cdot l)(X_i) = {\Bbb I}_{\cal A}(l)(X^p_i) \quad \mbox{\rm modulo $p$}.
\end{equation}
\vskip 3mm
{\bf Example}. Take an element $\delta_i^{\mathbf i} 
\in {\cal X}(\Z^t)$ 
whose coordinates $(x_1,
..., x_n)$ 
in the coordinate system related to a seed ${\mathbf i}$ are $x_j =0$ for $j \not =
    i$, $x_i =1$. 
Setting 
${\Bbb I}_{\cal X}(\delta_i^{\mathbf i} ) :=
    A_i^{\mathbf i}$ we get  a  universally
Laurent 
polynomial by the Laurent phenomenon theorem  of \cite{FZ3}. Its positivity was conjectured in {\it loc. cit.}.  
Therefore the cluster algebra 
sits inside of the algebra of universally Laurent polynomials 
${\Bbb L}({\cal A})$.
However the latter can be strictly 
bigger than the former.
 \vskip 3mm

The set of points  of a positive space ${\cal X}$ with values 
in a semifield is determined by a single 
positive coordinate system on ${\cal X}$. Contrary to this,  
the semiring of  universally positive Laurent polynomials depends on the 
whole collection of positive coordinate systems on ${\cal X}$. So the source of the 
canonical map  is determined by a single coordinate system of a positive atlas on the source space,  while its
image depends on the choice of a positive atlas on the
target space. This shows, for instance, that the set  
of positive coordinate systems on the
target space can not be ''too big'' or ''too small''.

\begin{conjecture} \label{10.10.03.10}
There exist  canonical $\Gamma$-equivariant pairings between the 
real tropical and real positive spaces
\begin{equation} \label{11.23.03.24e}
{\rm I}_{\cal A}: {\cal A}(\R^t)\times {\cal X}^{\vee}(\R_{>0})
\lra {\mathbb R}, \quad 
{\rm I}_{\cal X}: {\cal X}(\R^t) \times {\cal A}^{\vee}(\R_{>0})
\lra {\mathbb R}, \quad 
\end{equation} 
as well as a canonical $\Gamma$-equivariant intersection pairing between the real tropical spaces
\begin{equation} \label{11.23.03.1}
{\cal I}: {\cal A}(\R^t)\times {\cal X}^{\vee}(\R^t) 
\rightarrow {\mathbb R},
\end{equation} 

We expect the following properties of these pairings:
\begin{enumerate}
\item\label{convexity}  All the pairings are  convex
  functions in each of the two variables, and in each of the cluster
  coordinate systems.   

\item\label{homogeneos} All the  pairings are homogeneous with respect
  to the tropical variable(s): for any $\alpha >0 $ one has 
$
{\rm I}_*(\alpha l, m) = \alpha {\rm I}_*(l, m) 
$,  and the same holds for ${\cal I}$.

\item\label{symmetry-add} 
Let 
$l \in {\cal A}(\R^t)$ and $m$ is a point of 
${\cal A}^{\vee}(\R)$. Then 
 $$
{\rm I}_{\cal A}(l, p(m))={\rm I}_{\cal X}(p(l), m). 
$$
Similarly, ${\cal I}(l, p(m)) = 
{\cal I}(p(l), m)$.

\item \label{explicit} If the coordinates $(x_1, ..., x_n)$ of 
a point $l \in {\cal X}(\R^t)$ are positive numbers, 
and a point $m \in {\cal A}^+$ has the logarithmic
coordinates $(a_1, ..., a_n)$, then 
$$
{\rm I}_{\cal X}(l,m) = \sum_{i \in I}x_ia_i, 
$$
and the same holds for ${\cal I}$.  

\item\label{limit-add} Let $m$ be a point of either ${{\cal
    X}^{\vee}}(\R_{>0})$ or 
${\cal A}^{\vee}(\R_{>0})$ 
with logarithmic coordinates $u_1,\ldots, u_n$. Let $C \in \R$. Denote by 
$C\cdot m$ the point with coordinates $C u_1,\ldots, C u_n$.
Let $m_L$ be the point of either ${\cal
    A}(\R^t)$ or 
${\cal X}(\R^t)$  with the coordinates $u_1,\ldots, u_n$. Then, 
for both $* = {\cal A}$ and $* = {\cal X}$, 
\begin{equation} \label{11.26.03.2}
\lim\limits_{C \rightarrow \infty}{\rm I}_*(l, C \cdot m)/C = {\cal I}(l, m_L).
\end{equation}

\item The intersection pairing (\ref{11.23.03.1})
restricts to a pairing
\begin{equation} \label{11.23.03.2}
{\cal I}: {\cal A}(\Q^t)\times {\cal X}^{\vee}(\Q^t) 
\rightarrow {\mathbb Q}.
 \end{equation} \end{enumerate}
\end{conjecture}

{\bf Remark}. The convexity property implies that  the intersection pairing is
continuous. So an intersection pairing (\ref{11.23.03.2}) satisfying the
convexity property \ref{convexity} determines the pairing (\ref{11.23.03.1}). 
Similarly convex pairings 
\begin{equation} \label{11.23.03.24ef}
{\rm I}_{\cal A}: {\cal A}(\Q^t)\times {\cal X}^{\vee}(\R_{>0})
\lra {\mathbb R}, \quad 
{\rm I}_{\cal X}: {\cal X}(\Q^t) \times {\cal A}^{\vee}(\R_{>0})
\lra {\mathbb R}
\end{equation} 
can be uniquely extended to continuous pairings (\ref{11.23.03.24e}).

\vskip 3mm

\subsection{Conjecture \ref{10.10.03.10*} essentially 
implies Conjecture \ref{10.10.03.10}}
 Conjecture \ref{10.10.03.10*} is an algebraic cousin of 
Conjecture \ref{10.10.03.10}. Roughly speaking, it is obtained by replacing 
in the pairings (\ref{11.23.03.24e}) the tropical
semifield 
$\R^t$ by its integral version $\Z^t$,
and the real manifolds ${\cal X}^{\vee}(\R_{>0})$ and 
${\cal A}^{\vee}(\R_{>0})$ by the corresponding 
cluster varieties ${\cal X}^{\vee}$ and 
${\cal A}^{\vee}$. It is handy to 
think about these algebraic pairings as of maps from the sets of integral tropical points  to
positive regular functions on the corresponding scheme. 
Observe that the positive structure of the corresponding schemes has been used  
to define the positive regular functions on them.

One can interpret the map ${\Bbb I}_*$ as a pairing 
${\Bbb I}_*(\ast, \ast)$ between laminations and points 
of the corresponding space: ${\Bbb I}_*(l, z):= {\Bbb I}_*(l)(z)$.  
We are going to show how the canonical pairings ${\cal I}(\ast, \ast)$ 
and ${\rm I}_*(\ast, \ast)$ 
emerge from the one ${\Bbb I}_*(\ast, \ast)$ in the 
tropical and scaling limits. 

{\it The tropical limit and the intersection pairing ${\cal I}(\ast, \ast)$}. 
Let  $F(X_i)$ be a positive integral Laurent polynomial. 
Let $F^t(x_i)$, where $x_i$ belong to a semifield, 
 be the  corresponding tropical polynomial. 

\vskip 2mm
{\bf Example}. If $F(X_1, X_2) = X_1^2 + 3 X_1X_2$, then 
$F^t(x_1, x_2) = {\rm max}\{2x_1,  x_1 + x_2\}$. 
\vskip 3mm

Observe that one has 
\be \la{evtr}
\lim_{C \to \infty}\frac{\log (e^{Cx_1} + ... + e^{Cx_n})}{C} = 
{\rm max}\{x_1, ..., x_n\}.
\ee
Therefore the Laurent polynomial $F$ 
and its tropicalization $F^t$ are related as follows: 
\begin{equation} \label{12.6.04.15}
\lim_{C \to \infty}\frac{\log F(e^{Cx_1}, ..., e^{Cx_n})}{C} = F^t(x_1, ..., x_n), \quad x_i \in \R.
\end{equation}

\vskip 2mm
\begin{definition} \label{12.6.04.10} Let us assume that we have 
the canonical maps ${\Bbb I}_*$ from 
Conjecture \ref{10.10.03.10*}. Then, for an integral tropical point 
$l$ of the $*$-space, where $*$ stands for either ${\cal A}$ or ${\cal X}$,  
the function ${\cal I}_*(l, \bullet)$ is the tropicalization of the 
positive integral 
Laurent polynomial ${\Bbb I}_*(l)$:
\begin{equation} \label{12.6.04.5}
{\cal I}_*(l, \bullet):= {\Bbb I}^t_*(l)(\bullet).
\end{equation} 
\end{definition} 
Therefore (\ref{12.6.04.15}) implies that 
$$
{\cal I}_*(l, m)= \lim_{C \to \infty}\frac{\log {\Bbb I}_*(l)(e^{C\cdot m})}{C}.
$$

\vskip 3mm
{\it The scaling limit}. 
Observe that the restriction of ${\Bbb I}_{\cal A}(l)$ to 
${\cal X}^{\vee}(\R_{>0})$, as well as 
${\Bbb I}_{\cal X}(l)$ to ${\cal A}^{\vee}(\R_{>0})$, 
are positive valued functions, so the logarithm 
$\log {\Bbb I}_{*}(l)$ makes sense. Recall that the group $\Q_+^*$ acts by automorphisms 
of the semiring $\Q^t$, and hence acts on the set of $\Q^t$-points of a positive space 
${\cal X}$. Namely, if $l \in {\cal X}(\Q^t)$ and $C \in \Q^*_+$, we denote by $C\cdot l$ 
the element obtained by multiplying all coordinates of $l$ by $C$, 
in any of the coordinate systems. 
Let us define pairings
\begin{equation} \label{11.23.03.24efd}
{\rm I}_{\cal A}: {\cal A}(\Q^t)\times {\cal X}^{\vee}(\R_{>0})
\lra {\mathbb R}, \quad 
{\rm I}_{\cal X}: {\cal X}(\Q^t) \times {\cal A}^{\vee}(\R_{>0})
\lra {\mathbb R} 
\end{equation} 
by taking the scaling limit, 
where, as usual,  $*$ stands for either ${\cal X}$ or ${\cal A}$. 
\begin{equation} \label{11.26.03.1}
{\rm I}_*(l, u):= \lim_{C \to \infty}\frac{\log{\Bbb I}_*(C\cdot l)(u)}{C}.
\end{equation} 

\begin{conjecture} \label{12.6.04.1}
The scaling limit (\ref{11.26.03.1}) exists for both 
$* = {\cal A}$ and $*= {\cal X}$. 
\end{conjecture} 

\begin{theorem} \label{12.6.04.2}
Let us assume Conjectures \ref{10.10.03.10*} and \ref{12.6.04.1}. 
Then the pairings ${\rm I}_*(l, u)$ and ${\cal I}(l,m)$ satisfy all conditions of 
Conjecture \ref{10.10.03.10}. 
\end{theorem}

{\bf Proof}. Definition \ref{12.6.04.10} and (\ref{11.26.03.1}) 
provide the pairings for integral tropical laminations $l$. 
First of all we have to extend them to rational and real tropical points.

\begin{lemma} \label{12.6.04.2s}
Under the assumptions of Theorem \ref{12.6.04.2},  
${\rm I}_*(l, u)$ enjoys the following properties:

i) It is homogeneous in $l$. 

ii) It is convex in the $l$-variable, and it is convex in the $u$-variable. 
\end{lemma} 

{\bf Proof}. The property i) is clear. To check ii) for the Teichm\"uller, 
$u$-variable, observe that 
$$
{\Bbb I}_*(l, e^{x}) \cdot {\Bbb I}_*(l, e^{y}) \geq {\Bbb I}_*(l, e^{x + y}) 
$$
Indeed, $(\sum a_Ie^{x_I})(\sum a_Ie^{y_I}) \geq (\sum a_Ie^{x_I+y_I})$ since 
$a_I$ are positive integers. So taking the logarithm we get a convex function, 
and a limit of convex functions is convex. 
Similarly the property ii) for the tropical, $l$-variable, follows from 
${\Bbb I}_*(l_1){\Bbb I}_*(l_2) \geq {\Bbb I}_*(l_1+l_2)$, which 
 is an immediate corollary of $c_*(l_1, l_2; l_1+ l_2)=1$ and 
the property (\ref{prod1}) 
in  Conjecture \ref{10.10.03.10*}.
The lemma is proved. 
\vskip 3mm 

Using homogeneity in $l$ we extend the pairings to rational tropical points $l$. 
Then convexity of the pairings (\ref{11.23.03.24efd}) implies that they are 
continuous in the natural topology on the set of rational tropical points, and thus 
can be uniquely extended to real pairings (\ref{11.23.03.24e}). 
So we get the pairings ${\rm I}_*$ and ${\cal I}$ defined for any 
real tropical point $l$ of the corresponding space, which satisfy 
Properties 1 and 2. Property  6 is evidently valid. 
Properties 3 and 4 follow immediately from Properties 3 and 2
in Conjecture \ref{10.10.03.10*}. Finally, Property 5   
follows from the very definition:
$$
\lim_{C_1, C_2 \to \infty}\frac{\log {\Bbb I}_*(C_1\cdot l)(C_2\cdot u)}{C_1C_2}= 
\lim_{C_1 \to \infty}\frac{{i}_*(C_1\cdot l)(u_L)}{C_1} \stackrel{1.}{=}
{\cal I}_*(l,u_L).
$$
The theorem is proved.

\subsection{Quantum canonical map} We define a {\it universally positive 
Laurent polynomial} related to the quantum space ${\cal X}_{q}$ as 
an element of the 
(non commutative) fraction field ${\Bbb T}^q_{\mathbf i}$ which 
for any seed ${\mathbf i}'$ is a Laurent polynomial,  
with positive integral coefficients, 
in $q$ and the corresponding 
quantum ${\cal X}$-coordinates $X_i$. We denote by 
${\Bbb L}_+({\cal X}_{q})$ the {\it semiring of universally positive 
Laurent polynomials} on the quantum space ${\cal X}_{{q}}$. 
 Dropping the condition of positivity of the coefficients we get 
the {\it ring of universally 
Laurent polynomials} ${\Bbb L}({\cal X}_{q})$ 
on the quantum space ${\cal X}_{q}$. 

Let us make a few  remarks preceding the conjecture.

\vskip 3mm
1. Given an order on the set $I$, an element $a = (a_1, ..., a_n) \in \Z^n$
 provides a monomial 
\begin{equation} \label{9.11.04.2}
X_a:= q^{-\sum_{i<j}
\widehat {\varepsilon}_{ij}a_ia_j}\prod_iX_i^{a_i}.
\end{equation} 
It  is $\ast$-invariant and 
does not depend on the choice of ordering of $I$ used 
to define it (Section 3.1).

\vskip 2mm
2. For any positive space ${\cal X}$ and a positive integer $N$
there is a subset ${\cal X}(N\Z^t) \subset {\cal X}(\Z^t)$. It
consists of the points whose coordinates in one, and hence in any positive coordinate
system for ${\cal X}$ are integers divisible by $N$. 
The canonical map of positive spaces $p: {\cal A} \to {\cal
  X}$
induces the map 
\begin{equation} \label{5.7.04.10}
p: {\cal A}(\Z^t) \to {\cal
  X}(\Z^t).
\end{equation}  
Let us take the inverse image of the subset ${\cal X}(N\Z^t)$ under
this map, and consider the positive abelian semigroup it generates. We denote
it by 
$Z_{{\cal A}}(N)$:
$$
Z_{{\cal A}}(N):= \Z_+\Bigl\{p^{-1} \left( 
{\cal X}(N\Z^t)\right)\Bigr\}.
$$
Observe that $\alpha \in {\rm Ker}_L[\ast, \ast]$
 provides an element $l(\alpha)\in {\cal A}(\Z^t)$ which lies in 
$Z_{{\cal
    A}}(N)$ 
for all $N$.

\vskip 2mm
3. Recall that for a split torus $H$ the set $H(\Z^t)$ is 
the group of cocharacters 
$X_*(H)$ of $H$. In particular 
it is an 
   abelian group. The canonical isomorphism 
$$
X_*(H_{\cal A}) = X^*({H_{{\cal X}^{\vee}}})
$$ 
allows us to consider an element 
${\alpha} \in H_{\cal A}(\Z^t)$ as a 
character $\chi_{\alpha}$ of the torus
${H_{{\cal X}^{\vee}}}$. The inverse image of this character under the projection 
$\theta_q^{\vee}$
delivers  a central  element  
\begin{equation} \label{1.9.04.122}
\widetilde \chi^q_{\alpha}:= (\theta_q^{\vee})^*\chi_{\alpha} 
\in {\rm Center}\left(
{\Bbb L}_+({\cal X}_{q}^{\vee})\right).
\end{equation} 

\vskip 2mm
4. The torus ${H}_{\cal A}$ acts on the
${\cal A}$-space. So the abelian group 
${H}_{\cal A}(\Z^t)$ acts on the set ${\cal A}(\Z^t)$. 
We denote this action by $*$. 
In coordinates it looks as follows. Choose a  seed ${\mathbf i}$. 
Then an element $\alpha \in {H}_{\cal A}(\Z^t)$ is given by $\{\alpha_1, 
..., \alpha_n\}$ where 
$\alpha_i \in \Z$ and $\sum_j\varepsilon_{ij} \alpha_j =0$ for all $i \in I$. 
The element  $\{\alpha_1, ..., \alpha_n\} \in 
{H}_{\cal A}(\Z^t)$ acts on 
$\{\beta_1, ..., \beta_n\} \in 
{\cal A}(\Z^t)$  by 
$$
\{\alpha_1, ..., \alpha_n\} * \{\beta_1, ..., \beta_n\}= \{\alpha_1 + \beta_1,
..., 
\alpha_1 + \beta_n\}.
$$

\vskip 2mm
5. Let $p$ be a prime, $q$ a primitive $p$-th root of unity. Then $\Z[q]$ is
  the cyclotomic ring. Let $(1-q)$ be the ideal generated  $1-q$. 
It is the kernel of the surjective map 
$$
\pi_p: \Z[q] \lra \Z/p\Z, \quad q \lms
  1.
$$
For a $\Z[q]$-module $M$ denote by $\pi_p(M)$ the reduction of $M$ modulo the ideal $(1-q)$. 
So $\pi_p(M)$ is an $F_p$-vector
  space. 

\vskip 2mm
The following conjecture is the main conjecture in our paper.

\begin{conjecture} \label{10.10.03.10fgf}
There exists a quantum canonical map, that is, a $\Gamma$-equivariant  isomorphism 
\begin{equation} \label{10.q10.03.10fgf}
\widehat {\Bbb I}^q: \Z_+\{{\cal A}(\Z^t)\} 
\stackrel{\sim}{\lra} {\Bbb L}_+({\cal X}_{q}^{\vee})
\end{equation}
satisfying
the following properties:
\begin{enumerate}
\item $\widehat{\Bbb I}^1(l)={\Bbb I}_{\cal A}(l)$.
\item Let $a = (a_1, ..., a_n)$ be the coordinates of $l \in 
{\cal A}(\Z^t)$ in a cluster coordinate system.  Then the highest   
term of $\widehat{\Bbb I}^q(l)$ is 
$
 X_a = 
q^{-\sum_{i<j}\widehat {\varepsilon}^{ij}a_ia_j}\prod_iX_i^{a_i}.
$ 
\item {\em Self-duality}: 
$\ast \widehat{\Bbb I}^q(l)=\widehat{\Bbb I}^q(l)$.
\item $\widehat{\Bbb I}^q(l_1)\widehat{\Bbb I}^q(l_2)=\sum_l c^q(l_1, l_2;l)
 \widehat{\Bbb I}^q(l)$, where the $c^q(l_1, l_2;l)$
are Laurent polynomials of $q$ with positive integral coefficients. 
Moreover $c^q(l_1, l_2; l_1+ l_2)=1$. 
\item {\em The center at roots of unity}:  
Let $N$ be a positive integer and $q$  a primitive $N$-th root of unity. 
Then the map $\widehat {\Bbb I}^q $ induces an isomorphism
\begin{equation} \label{1.9.04.223}
\widehat {\Bbb I}^q \Bigl(Z_{{\cal A}}(N)\Bigr) 
\stackrel{\sim}{\lra} {\rm Center}\left(
{\Bbb L}_+({\cal X}_{q}^{\vee})\right).
\end{equation}
\item {\em The quantum Frobenius on ${\Bbb L}_+({\cal X})$}: 
Let $N \in \Z_{>0}$ and $q$ be a primitive $N$-th root of unity. 
Then the quantum Frobenius map 
${\Bbb F}^*_N$, (see Theorem \ref{QFROB}), 
is related to the map $\widehat {\Bbb I}^q$ 
as follows: 
$
{\Bbb F}^*_N {\Bbb I}_{\cal A}(l)
 = \widehat {\Bbb I}^q (N\cdot l)
$,  
 i.e. in every cluster coordinate system $\{X_i\}$ one has 
\begin{equation} \label{1.9.04.222}
 {\Bbb F}^*_N {\Bbb I}_{\cal A}(l)(X_i):= {\Bbb I}_{\cal A}(l)(X_i^N)
 = \widehat {\Bbb I}^q (N\cdot l)(X_i).
\end{equation}
\item Let $p$ be a prime, $q$ a primitive $p$-th root of unity. Then there is
 a canonical isomorphism 
\begin{equation} \label{5.9.04.1}
\pi_p\Bigl({\Bbb L}_+({\cal X}_q)\Bigr) =  
{\Bbb L}_+({\cal X})\otimes \Z/p\Z.
\end{equation}
%DIVIDED POWERS NEEDED HERE??? HOW TO FORMULATE?
\item The map $\widehat {\Bbb I}^q $ transforms the action of an element 
$\alpha \in H_{\cal A}(\Z^t)$ on ${\cal A}(\Z^t)$ to the 
multiplication 
by the corresponding central element $\widetilde \chi^q_{\alpha}$, given by 
(\ref{1.9.04.122}):  
\begin{equation} \label{1.9.04.111}
\widehat {\Bbb I}^q (\{\alpha  * \beta\})  = \widetilde \chi^q_{\alpha} 
\cdot \widehat {\Bbb I}^q (\{ \beta\}).
\end{equation}
\item  Restricting the scalars from $\Z_+$ to $\Z$ we get isomorphism
$
\widehat {\Bbb I}^q: \Z\{{\cal A}(\Z^t)\} 
\stackrel{\sim}{\lra} {\Bbb L}({\cal X}_{q}^{\vee})
$.  
\end{enumerate}
\end{conjecture}

{\bf Remarks}. 1. The map ${\Bbb F}_N^*$ was defined only for those $N$ which satisfy
the condition formulated in Theorem \ref{QFROB}. 
Formula (\ref{1.9.04.222}) suggests a definition 
of the quantum Frobenius map ${\Bbb F}_N^*$ on the 
algebra ${\Bbb L}({\cal X})$ for all $N$. 
Indeed, the first equality in (\ref{1.9.04.222}) serves as a definition 
of ${\Bbb F}_N^*$ in a given cluster coordinate system $\{X_i\}$, 
and the second would imply that it is independent of the choice 
of the coordinate system. 

2. The isomorphism (\ref{5.9.04.1}) plus (\ref{1.9.04.222}) obviously imply 
Frobenius Conjecture (\ref{5.7.04.1}).

3. The property (\ref{1.9.04.111}) in the coordinate form looks as follows:
 assuming $\sum_j \varepsilon_{ij}\alpha_j =0$, 
$$
\widehat {\Bbb I}^q (\{\alpha_1 +\beta_1, ..., \alpha_n+\beta_n\}) = q^{-\sum_{i<j}
\widehat {\varepsilon}_{ij}\alpha_i\alpha_j}\prod_iX_i^{\alpha_i} \cdot 
\widehat {\Bbb I}^q (\{\beta_1, ..., \beta_n\}).
$$
Property (\ref{1.9.04.111}) can not even be stated 
if we skip the Langlands dual (defined at the end of Section 1.2) 
in (\ref{10.q10.03.10fgf}), replacing  
${\cal X}_{q}^{\vee}$ by ${\cal X}_{q}$. 
Indeed, there is no canonical isomorphism between 
$X_*(H_{\cal A})$ and $X^*(H_{{\cal X}})$. 
\vskip 3mm

\subsection{An example: the map ${\Bbb I}_{\cal A}$ for the cluster ensemble of type $A_2$} 
Let
$\varepsilon_{ij}=\left(\begin{array}{cc}0&1\\-1&0\end{array}\right)$.
Then the cluster modular group is
${\mathbb Z}/5{\mathbb Z}$. Its generator acts on the ${\cal X}$-space by 
$$
(X,Y)\mapsto
(Y(1+X^{-1})^{-1},X^{-1}),  
$$ and on the tropical ${\cal A}$-space and by $(a,b)\mapsto
(b,\max(b,0)-a)$. The canonical map 
${\Bbb I}_{\cal A}$
is given by:
$$
{\Bbb I}_{\cal A}(a,b)=\left\{\begin{array}{lll} X^{a}Y^b& \mbox{ for } & a\leq 0
\mbox{ and
} b\geq 0\\
\left(\frac{1+X}{XY}\right)^{-b}X^{a}& \mbox{ for } & a\leq 0
\mbox{
and } b\leq 0\\
\left(\frac{1+X+XY}{Y}\right)^a\left(\frac{1+X}{XY}\right)^{-b}&
\mbox{ for } & a\geq 0 \mbox{ and } b\leq 0\\
((1+Y)X)^{b}\left(\frac{1+X+XY}{Y}\right)^{a-b}&
\mbox{ for } & a\geq b\geq 0\\
Y^{b-a}((1+Y)X)^a &\mbox{ for } & b\geq a\geq 0.\\
\end{array}
\right.
$$
Or equivalently
$$
{\Bbb I}_{\cal A}(a,b)=\left\{\begin{array}{lll} X^{a}Y^b& \mbox{ for } & a\leq 0
\mbox{ and } b\geq 0\\
X^a Y^b(1+X^{-1})^{-b} & \mbox{ for } & a\leq 0 \mbox{
and } b\leq 0\\
X^a Y^b(1+X^{-1})^{-b}(1+Y^{-1}+X^{-1}Y^{-1})^a&
\mbox{ for } & a\geq 0 \mbox{ and } b\leq 0\\
X^a Y^b (1+Y^{-1})^b(1+Y^{-1}+X^{-1}Y^{-1})^{a-b}&
\mbox{ for } & a\geq b\geq 0\\
X^aY^b(1+Y^{-1})^a&\mbox{ for } & b\geq a\geq 0.\\
\end{array}
\right.$$ 
One can easily verify, that the formulae agree on the
overlapping domains of values of $a$ and $b$. 
A more elaborate and geometric discussion of this example 
see in \cite{G3}. 

%Product formulae
%$$
%{\Bbb I}_{\cal A}(a,b) {\Bbb I}_{\cal A}(a',b')= {\Bbb I}_{\cal
%  A}(a+a') {\Bbb I}_{\cal A}(b+b') \quad \mbox{for} \quad 
%a,a'\leq 0,\quad b,b'\geq 0
%$$
%$$
%{\Bbb I}_{\cal A}(a,b) {\Bbb I}_{\cal A}(a',b')=\sum {\min(|b|,|b'|)
%\choose i}{\Bbb I}_{\cal A}(a+a'-i,b+b') \quad \mbox{ for
%} a,a',b'\leq 0,\quad b\geq 0
%$$
\vskip 3mm

{\bf Remark}. By freezing some of the variables of a cluster ensemble, 
we prohibit mutations in the direction of  these variables, 
and therefore change the spaces of
universally positive Laurent polynomials for the corresponding ${\cal A}$ and 
${\cal X}$ spaces. However we do not change the sets of the points of these
spaces with values in any semifield. 
So by freezing  some of the variables 
we do not change the source of a canonical map, we do change 
its target.

\section{Canonical pairings in the finite type case} \la{Sec4.5}

\subsection{The canonical pairing between the tropical spaces}

Given a seed ${\bf i}$, let us consider the following function\footnote{
We will see below that it is the 
tropicalization of the more fundamental element (\ref{BCE}).} $P_{\bf i}$ on 
${\cal A}(\R^t) \times {\cal X}(\R^t)$:
\begin{equation} \label{6.01.04.10}
P_{\bf i}:= \sum_{i\in
  I}d_ia_ix_i.
 \end{equation}
For a finite type cluster ensemble $({\cal A}, {\cal X})$ we 
define a version ${\cal I}'$ 
of the canonical pairing between the  tropical spaces by maximizing 
this function over the (finite) set of all seeds:
\be \la{CP}
{\cal I}': {\cal A}(\R^t)\times {\cal X}(\R^t) \lra \R, \qquad {\cal I}'(a,x) := 
{\rm max}_{\{\mbox{\it seeds {\bf i}}\}}\sum_i d_ia_i x_i.
\ee
In the duality conjectures, however, we are looking for a canonical pairing 
involving the tropical points of the Langlands dual cluster ${\cal X}$-variety 
${\cal X}^\vee$: 
\be \la{CPa}
{\cal I}: {\cal A}(\R^t)\times {\cal X}^\vee(\R^t) \lra \R.
\ee
This apparent contradiction is resolved by the following Lemma. 

\bl \la{LL}
Let $x = \{x_i\} \in {\cal X}({\Bbb Z}^t)$. Then the coordinates $\{d_ix_i\}$ 
behave under the mutations in the tropical space ${\cal X}({\Bbb Z}^t)$ 
just like the coordinates of a point $\delta(x)$ of the Langlands dual tropical space ${\cal X}^{\vee}({\Bbb Z}^t)$. So there is a canonical $\Gamma$-equivariant inclusion
$$
\delta: {\cal X}({\Bbb Z}^t) \hra {\cal X}^{\vee}({\Bbb Z}^t),
$$ 
given in any cluster coordinate system by $\{x_i\} \lms \{d_ix_i\}$. 
The map $\delta$ is an isomorphisms for rational or real tropical points. 
\el

{\bf Proof}. The tropical mutation formulas are 
$$
x_i^\sharp = x_i - [\varepsilon_{ik}]_+{\rm max}(0, x_k), \qquad 
x_i' =  \left\{\begin{array}{lll} -x_k& \mbox{ if } & i=k, \\
    x_i + [\varepsilon_{ik}]_+ x_k  & \mbox{ if } &  i\neq k. \\
\end{array}\right.
$$ 
So we get the claim by changing 
$x_j \lms d_j x_j$ and using positivity of $d_k$ as well as the formula (\ref{xc}). 

 \vskip 2mm

Lemma \ref{LL} tells that the pairing ${\cal I}'$  determines
 the 
restriction of the canonical pairing (\ref{CPa})
to the subset $\delta ({\cal X}(\Z^t))$ of 
points of ${\cal X}^\vee(\Z^t)$ with the coordinates $\{d_ix_i\}$. 
Let us check now that the defined this way canonical pairing 
${\cal I}(a,\delta(x))$ satisfies the crucial part 4) of Conjecture \ref{10.10.03.10}.

The latter tells that 
if $y \in {\cal X}^\vee(\Z^t)$ has non-negative coordinates $\{y_i\}$ 
in a cluster coordinate system, then 
one should have  ${\cal I}(a,y) = \sum_i a_iy_i$. 
Theorem \ref{5.31.05.121} guarantees that this is the case 
for the points   $y$ of $\delta({\cal X}(\Z^t))$. Indeed, if 
$y = \delta(x)$ then $y_i = d_i x_i$, so ${\cal I}(a,y) = \sum_i d_ia_ix_i$.

\bt \la{5.31.05.121} Suppose that $({\cal A}, {\cal X})$ a finite type cluster ensemble. 
Then, for any $(a,x) \in {\cal A}(\R^t) \times {\cal X}(\R^t)$, there exists 
a seed ${\bf i}$ at which 
the maximum in (\ref{6.01.04.10})  is attained, 
such that all coordinates $x_i$ 
of $x$ at ${\bf i}$
are non-negative.
\et

{\bf Proof}. Let us investigate how the function $P_{\bf i}$ 
behaves under mutations. Let $\{a_i\}$ be the coordinates of $a \in {\cal A}(\R^t)$. Set 
\begin{equation} \label{6.01.04.15}
\varphi_k(a):= (p^*x_k)(a) = \sum_{j}\varepsilon_{kj}a_j. 
\end{equation}

\begin{proposition} \label{1.6.04.3} Given a mutation $\mu_k: {\bf i} \to {\bf i'}$, 
we have 
$$
({P}_{\bf i'} - {P}_{\bf i})(a,x) = 
\left\{\begin{array}{lll} d_k |x_k| \varphi_k(a) 
 & 
\mbox{ \rm if } x_k \varphi_k(a)< 0\\
0& \mbox{\rm otherwise}.\\
\end{array}
\right. 
$$
This implies the following: 

a) If $x_k = 0$, the 
  mutation at the vertex $k$ does not change the function $P_{\bf i}$.

b) If $x_k > 0$, the 
  mutation at the vertex $k$ strictly decreases the function 
$P_{\bf i}$ if $\varphi_k(a) <0$, 
and does not change it otherwise. 

c) If $x_k < 0$, the 
  mutation at the vertex $k$ strictly increases the function $P_{\bf i}$
 if $\varphi_k(a) >0$, 
and does not change it otherwise. 
\end{proposition}
 
We show in Section 6.1 that this is an easy consequence of Basic Lemma \ref{BL}.
\vskip 2mm

\bl \la{8.3.09.1}
Let ${\cal X}$ be an arbitrary cluster ${\cal X}$-variety. 
Let ${\bf i_0} \to {\bf i_1} \to \ldots \to {\bf i_n} \to {\bf i_0}$ be a sequence of mutations. 
Denote by $k_s$ the direction of the mutation ${\bf i_s} \to {\bf i_{s+1}}$. 

Let us assume that there exists a tropical point $x \in {\cal X}(\R^t)$ 
such that for every $s = 0, ..., n$ the $x_{k_s}$-coordinate of the point $x$ for the seed 
${\bf i}_s$ is 
positive. Then the sequence is trivial, i.e. $n=0$. 
\el

{\bf Proof}. Take a point $a \in {\cal A}(\R^t)$ with $\varphi_{k_0}(a) < 0$. 
Then, by Proposition \ref{1.6.04.3}, after the cluster transformation 
${\bf i_0} \to {\bf i_1} \to \ldots \to {\bf i_n} \to {\bf i_0}$ 
the sum $\sum d_ia_ix_i$ will strictly decrease. 
This contradiction proves the Lemma. 
\vskip 2mm

Take a point $(a,x) \in {\cal A}(\R^t) \times {\cal X}(\R^t)$. 
Take a seed ${\bf i}$ which realizes the maximum  of sum (\ref{6.01.04.10})
 evaluated at this point. 
If all coordinates $x_i$ of $x$ in this seed are non-negative, we are
  done. If not, there exists a vertex $k$ such that $x_k <0$. 
Let us perform a mutation at $k$. If the new
  coordinate system is non-negative for $x$, we are done. If not, we perform a
  mutation at a vertex $k'$ such that $x_{k'} <0$, and so on. Since 
${\cal X}$ is of finite type, after a finite number of mutations we
  get to a certain coordinate system for the second time. 
This contradicts Lemma \ref{8.3.09.1}, and thus proves the Theorem.

\begin{corollary} \label{5.31.05.121c}
Let ${\cal X}$ be a finite
 type cluster ${\cal X}$-variety. Then for every 
$x \in {\cal X}(\R^t)$ there exists a cluster coordinate system 
such that the coordinates $x_i$ of the point $x$ are non-negative. 
\end{corollary}

\paragraph{The positive part of a tropical space} 
Given a positive space ${\cal X}$ 
the set of the points of the tropical space ${\cal X}(\R^t)$ 
which have non-negative coordinates in a given coordinate 
system is a convex cone. We call it the {\it positive cone} assigned to the coordinate system. 
The union of the positive cones forms the {\it positive part} ${\cal X}(\R^t)_+$ 
of ${\cal X}(\R^t)$.

\begin{definition} \label{5.31.05.11} 
A positive space  ${\cal X}$ is of definite (semi-definite, indefinite) type if 
the subset ${\cal X}(\R^t)_+$ is equal to (respectively dense, not dense) 
in ${\cal X}(\R^t)$. 
\end{definition} 
Corollary \ref{5.31.05.121c} tells that a finite type cluster ${\cal X}$-variety is of definite type.

\begin{conjecture} \label{5.31.05.12}
A cluster ${\cal X}$-variety is of finite type if and only if it is of definite type. 
\end{conjecture}

\vskip 2mm
{\bf Examples}. 1. Let $S$ be a surface with holes and marked points on the boundary. 
Denote by $h$ the number of holes without marked points on its boundary. 
There is an action of the group $(\Z/2\Z)^h$ 
 by birational automorphisms of the moduli space 
${\cal X}_{PGL_2, S}$, see \cite{FG1}.  It acts by cluster transformations 
if and only if $h>1$. 
Therefore, by Theorem 12.2
in {\it loc. cit.},   the cluster atlas on ${\cal X}_{PGL_2, S}$  
is of semi-definite type if and only if $h>1$. 

However even if $h=1$, the group $\Z/2\Z$ 
acts by positive transformations which leave invariant the semiring ${\Bbb L}_+({\cal X}_{PGL_2, S})$ 
of positive regular functions. 
Therefore it is natural to extend the cluster atlas on ${\cal X}_{PGL_2, S}$ by 
adding the images of the cluster coordinate systems by the 
action of the group $\Z/2\Z$. We call  the obtained positive atlas the 
{\it extended cluster atlas}. 
The moduli space ${\cal X}_{PGL_2, S}$ 
is of semi-definite type for this atlas. 
%An example of the torus with hole is discussed in Section \ref{Sec6.3}. 

2. The canonical pairing ${\cal I}$ for the cluster 
ensemble $({\cal A}_{SL_2, S}, {\cal X}_{PGL_2, S})$ 
coincides with the intersection pairing between the ${\cal A}$- and ${\cal X}$-laminations defined 
in Section 12 of {\it loc. cit.}.

In Section 5.2 we show that 
for a finite type cluster ${\cal X}$-variety the positive  cones give a finite 
decomposition of the 
space ${\cal X}(\R^t)$.  It is dual to the generalized assaciahedron. 
Combining this with the Laurent Phenomenon Theorem \cite{FZ3} 
we  construct the canonical map 
${\Bbb I}_{\cal X}$ in the finite type case. 
Its tropicalization provides the canonical pairing (\ref{CP}).

\subsection{The canonical map ${\Bbb I}_{\cal X}$}
Any subset of vertices of a seed ${\bf i}$ provides us with a {\it subseed} 
${\bf j} \subset {\bf i}$. Mutating the seed ${\bf i}$ at a vertex  
of ${\bf j}$ we get a new seed with a subseed canonically identified with 
${\bf j}$. So we can mutate at its vertices, and so on. The 
obtained this way collection of seeds is called {\it the set of seeds ${\bf j}$-equivalent 
to a seed ${\bf i}$}. We use the notation 
${\bf i_1} \sim_{\bf j} {\bf i_2}$ for ${\bf j}$-equivalent seeds. 

Let $l \in {\cal X}(\Z^t)$. Let ${\bf i}$ be a non-negative 
seed for $l$, i.e. the $x$-coordinates of $l$ in this seed are non-negative. 
Let ${\bf j}(l)$ be the subseed of ${\bf i}$ determined by the 
zero $x$-coordinates of $l$, i.e. by the 
set of coordinates $x_i$ such that $x_i(l)=0$. 
We call it the {\it zero subseed for $l$}. 

\bt \la{EUNI} Let ${\cal X}$ be a finite type cluster ${\cal X}$-variety, and 
$l \in {\cal X}(\Z^t)$. Let ${\bf i}$ be a non-negative 
seed for $l$, and ${\bf j}(l)\subset {\bf i}$ the zero subseed for $l$. 
Let ${\bf i'}$ be another non-negative 
seed for $l$. Then the seeds ${\bf i}$ and ${\bf i'}$ are ${\bf j}(l)$-equivalent. 
\et

{\bf Proof}. 
Choose a cluster transformation ${\bf c}$ connecting the seeds ${\bf i}$ and ${\bf i'}$. 
Among the sequences of mutations realising 
 $\mathbf c$, choose the 
subset of sequences maximizing the minimal value of $P_{\mathbf s}(l)$, and 
in this subset choose a sequence where this minimum is attained minimal number 
of times. 
Due to the finite type assumption, 
the number of possible values is finite, so such a subset exists. 
Our goal is to show that for any such a sequence 
the $x$-coordinates of all mutating vertices are zero. Let 
$\mathbf i=\mathbf s_0,\mathbf s_1,\ldots,\mathbf s_g=\mathbf i'$ be 
our sequence of seeds

Consider the first seed $\mathbf s_r$ where the minimum of $P_{\mathbf s_r}(l)$ 
is attained for the first time. Suppose that we come to this seed by a mutation $\mu_j$ 
at the vertex $j$ and leave it by a mutaion $\mu_i$ at the vertex $i$. 
Recall that if the pair $(i,j)$ generates the standard $(h+2)$-gon, 
i.e. (\ref{K10}) holds, 
then the cluster transformation 
$\mu_j\mu_i$  equals to the one given by a sequence $\mu_i\mu_j\ldots$ of length $h$,
 where $h+2$ is the period of the sequence $\mu_j\mu_i\ldots$, times $\sigma_{ij}$ 
in the $A_2$ case.   
We shall show that if we replace the subsequence of mutations $\mu_i\mu_j$ 
in our cluster transformation 
by the sequence $\mu_j\mu_i\ldots$ of length $h$, 
then in the new sequence 
either the minimum of $P_{\mathbf s_i}(l)$ will be bigger or it will be attained 
smaller number of times, contradicting to our assumption.

To show this it is sufficient to prove the following

\bl
For the exchange functions 
$\varepsilon_{ij}$ corresponding to Dynkin diagrams 
$A_1\times A_1, A_2, B_2, G_2$,  the 
values $P_{\mathbf s}(l), P_{\mu_j\mathbf s}(l), P_{\mu_i\mu_j\mathbf s}(l), \ldots $ 
change growth to decay only once per period. 
\el

To prove this we will show that in this sequence the $x$-coordinates at the 
mutated vertices change their signs only once per period, and then 
use Proposition \ref{1.6.04.3}.

The claim we have to prove is an immediate corollary of the following observation. 
Recall the mutations $\mu$ (types $A_2$ and $A_1\times A_1$) 
or $\mu_\pm$ (types $B_2$ and $ G_2$). 
Consider the action of the 
sequence of cluster transformations $\mu$, $\mu^2$, $\mu^3$ ... or, respectively,  
$\mu_\pm$, $\mu_\mp\mu_\pm$, , $\mu_\pm\mu_\mp\mu_\pm$, ... on 
the tropical plane $(x_1, x_2)$. 
The sequence of 
$x_1$-coordinates of the points on the orbit of a point 
is the set we were looking for.  
Now Lemmas \ref{OHK1}, \ref{OHK0} and \ref{OHK2} imply that 
the $x_1$-coordinate  changes the sign just once per period. The Theorem is proved.

By the Laurent Phenomenon Theorem \cite{FZ3},   
${\Bbb I}_{\cal X}(l)$ is a universally Laurent polynomial. 

\begin{conjecture}
Theorem \ref{EUNI} is valid for an arbitrary cluster ${\cal X}$-variety. 
\end{conjecture}

%This conjecture would immediately imply a partial construction 
%of the canonical map ${\Bbb I}_{\cal X}$ for an arbitrary cluster ${\cal X}$-variety. 
%Indeed, it tells that there is a canonical bijection 

%\input{clustepd1.tex}
Theorems \ref{5.31.05.121} and \ref{EUNI} tell that the 
space ${\cal X}(\R^t)$ has a canonical decomposition 
into cones. The cones are parametrised by the 
cluster ${\cal X}$-coordinate systems.  Namely, such a cone is given by 
the set of all points with non-negative 
coordinates in a given cluster coordinate system.

\bp \la{PiR} The decomposition of ${\cal X}(\R^t)$ is dual to 
  the generalized associahedron. 
\ep

{\bf Proof}. Follows immediately from the combinatorial description of the 
$n-k$-dimensional faces of the cones 
implied by Theorem \ref{EUNI}: 
the cones are parametrised by the ${\bf j}$-equivalence classes of seeds, where 
 $|{\bf j}|=k$. Furthermore, the subcones of a given  cone are 
parametrised by the subseeds ${\bf j'}$ squeezed between ${\bf j}$ and ${\bf i}$. 

\vskip 2mm

\begin{conjecture}
In any cluster coordinate system the decomposition of 
${\cal X}(\R^t)$ is a decomposition into convex cones.  
\end{conjecture}
%This would provide a convex polytopal realization of the dual 
%polyhedron, which is a generalized associahedron by Proposition \ref{PiR}.   

\paragraph{Construction of the canonical map ${\Bbb I}_{\cal X}$ 
for finite type cluster ensembles.} It is provided by 
Theorems \ref{5.31.05.121} and \ref{EUNI} as follows. 
Given $l \in {\cal X}(\Z^t)$, Theorem \ref{5.31.05.121} tells that  
there exists a seed ${\mathbf i}$
such that all coordinates $(x_1, ..., x_n)$ of $l$ for 
this seed are non-negative. Set 
\be \la{TEX}
{\Bbb I}_{\cal X}(l) := A_1^{x_1} ... A_n^{x_n}
\ee
 where $(A_1, ..., A_n)$ are the 
${\cal A}$-coordinates for the Langlands dual seed ${\mathbf i}^{\vee}$. 
Such a seed ${\bf i}$ may not be unique. However 
Theorem \ref{EUNI} guaranties that expression (\ref{TEX}) 
does not depend on the choice of  
${\mathbf i}$. Indeed, it tells that any other seed ${\bf i'}$ 
is ${\bf j}(l)$-equivalent to  ${\bf i}$, and hence 
the expression (\ref{TEX}) for ${\bf i'}$ is the same as for ${\bf i}$. 
So ${\Bbb I}_{\cal X}(l)$ is well defined.

%\vskip 2mm
%It seems that the proof of the existence of 
%polytopal realizations of the generalised associahedra 
%provided by Theorem \ref{PiR} is simpler then the one in \cite{CFZ}. 

%\input{assocehedron.tex}

\section{Motivic avatar of the form $\Omega$ on the ${\cal A}$-space 
and the dilogarithm} \la{motivic}

Given a seed ${\bf i}$, we introduce a Casimir element 
${\bf P}_{\bf i}$. It does depend on the choice of a seed ${\bf i}$. 
Basic Lemma \ref{BL} provides a transformation formula for the element 
${\bf P}_{\bf i}$. Its applications include a  
proof of Proposition \ref{1.6.04.3} as well as the  
key properties of the motivic dilogarithm class introduced 
 in Section 6.4.

\subsection{The Basic Lemma}

Given a seed ${\bf i}$, the seed cluster tori ${\cal A}_{\bf i}$ 
and ${\cal X}_{\bf i}$ are dual to each other. 
Recall  the cluster 
${\cal A}$- and ${\cal X}$-coordinates 
$\{A_i\}$ and $\{X_i\}$  related to the seed ${\bf i}$. 
The set of of characters 
$\{X_i^{d_i}\}$  of the torus ${\cal X}_{\bf i}$ is dual to the basis of characters 
$\{A_i\}$ of the torus ${\cal A}_{\bf i}$, see (\ref{TT}). So 
there is a natural Casimir element 
\be \la{BCE}
{\bf P}_{\bf i}:= \sum_{i\in I}d_i \cdot A_i \otimes X_i \in \Q({\cal A})^*\otimes
\Q({\cal X})^*.
\ee
Here $\Q({\cal A})^*$ is the multiplicative group of the field 
of rational functions on ${\cal A}$, similarly $\Q({\cal X})^*$. 
We denote by $d \cdot \ast$ multiplication of an element $\ast$ of the 
tensor product by an integer $d$.

Let us investigate how the Casimir element ${\bf P}_{\bf i}$ changes under a mutation 
$\mu_k: {\bf i} \lra {\bf i'}$. 
Denote by $\{A'_i\}$ and $\{X'_i\}$ the cluster coordinates 
related to the seed ${\bf i'}$. 
Recall the notation 
$$
p^* X_k = \prod_{j\in I}A_j^{\varepsilon_{kj}}  = 
\frac{{\Bbb A}_k^+}{{\Bbb A}_k^-}\in \Q({\cal A})^*.
$$

\bl \la{BL} Given a mutation 
$\mu_k: {\bf i} \lra {\bf i'}$ at the vertex $k$, one has 
$$
{\bf P}_{\bf i'} - {\bf P}_{\bf i}  = 
d_k \Bigl( p^*  X_k\otimes (1+X_k)  - (1+ p^*  X_k)\otimes X_k \Bigr).
$$
\el

{\bf Proof}. 
We decompose the mutation $\mu_k$ into the composition $\mu_k = 
\mu_k' \circ \mu_k^\sharp$ (Section 2.1). The map 
$\mu_k'$ clearly preserves the element ${\bf P}_i$. 
Let us calculate the effect of the automorphism $\mu_k^\sharp$. 
We have 
$$
\sum_i d_i \cdot A_i' \otimes X_i' - \sum_i d_i \cdot A_i \otimes X_i 
\stackrel{(\ref{f3**}) -  (\ref{f3**as})}{=} 
$$
$$
\sum_i d_i \cdot A_i\otimes X_i(1+X_k)^{-\varepsilon_{ik}}  + 
(1+p^*X_k)^{-1}\otimes X_k 
 - \sum_i d_i \cdot A_i \otimes X_i = 
$$
$$
\sum_i d_i \cdot A_i\otimes (1+X_k)^{-\varepsilon_{ik}}  - d_k \cdot 
(1+p^*X_k)\otimes X_k 
 = 
$$
$$
\sum_i -d_i \varepsilon_{ik}\cdot  A_i\otimes (1+X_k) - d_k \cdot (1+p^*  X_k) 
\otimes X_k.
$$
Notice that $-d_i \varepsilon_{ik} = -  \varepsilon_{ik} = 
  \varepsilon_{ki} = d_k \varepsilon_{ki}$. So 
the first term equals to 
$$
\sum_i d_k \varepsilon_{ki}\cdot A_i\otimes (1+X_k)  = d_k \cdot  p^*  X_k
\otimes (1+X_k).
$$
The Basic Lemma is proved. 
\vskip 2mm

{\bf Proof of Proposition \ref{1.6.04.3}.} Notice that $P_{\bf i}$ is nothing else but the tropicalization of the element ${\bf P}_{\bf i}$. 
Furthermore, ${\rm max}(0, x) = [x]_+$. So the Basic Lemma implies 
$$
({P}_{\bf i'} - {P}_{\bf i})(a,x)  = 
d_k\Bigl( \varphi^*_k(a) ~[x_k]_+ -  [\varphi^*_k(a)]_+~x_k \Bigr).   
$$
If $x_k>0$, we get
$$
({P}_{\bf i'} - {P}_{\bf i})(a,x)  = 
d_k x_k\Bigl(\varphi^*_k(a) - 
[\varphi^*_k(a)]_+\Bigr) = 
\left\{\begin{array}{lll} 0 & 
\mbox{ if } \varphi^*_k(a)\geq 0\\
d_kx_k\varphi^*_k(a) <0& \mbox{ if } \varphi^*_k(a)< 0.\\
\end{array}
\right. 
$$
If $x_k<0$, we get 
$$
({P}_{\bf i'} - {P}_{\bf i})(a,x)  = - d_k x_k  [\varphi^*_k(a)]_+ = 
\left\{\begin{array}{lll} -d_kx_k\varphi^*_k(a) >0 & 
\mbox{ if } \varphi^*_k(a)> 0\\
0&\mbox{ if } \varphi^*_k(a)\leq 0.\\
\end{array}
\right.   
$$
The Proposition  is proved. 

\subsection{The group $K_2$, the Bloch complex and the dilogarithm} 
\paragraph{The Milnor group $K_2$ of a field.} Let $A$ be an abelian group. Recall that   $\Lambda^2A$ is the quotient 
of $\otimes^2A$ modulo the subgroup generated by the elements 
$a\otimes b + b\otimes a$. 
We denote by $a \wedge b$ the projection of $a\otimes b$ to the quotient. 

Let $F$ be an arbitrary field. The Milnor group $K_2(F)$ is an abelian group 
given as the quotient of the group $F^* \otimes F^*$ by the subgroup generated by the 
{\it Steinberg relations} $(1-x)\otimes x$ where $x \in F^*-\{1\}$. 
The image of the generator $x\otimes y$ in  $F^* \otimes F^*$ is denoted by $\{x,y\}$. 
It is well known \cite{Mi} that the Steinberg relations 
imply that $\{x,y\} = -\{y, x\}$. So one has 
\be \la{kk}
K_2(F) = \frac{\bigwedge^2F^*}{\mbox{Steinberg relations}}.
\ee

\paragraph{A regulator map on $K_2$.} Let $X$ be a smooth algebraic variety. Denote by $\Omega^2_{\rm log}(X)$ 
the space of $2$-forms with logarithmic singularities on $X$. 
Denote by $\Q(X)$ the field of rational functions on $X$. One has $d\log \wedge d\log ((1-f)\wedge f) =0$. 
So there is a group homomorphism
$$
d\log \wedge d\log: K_2(\Q(X)) \lra \Omega^2_{\rm log}(X), \qquad f\wedge g \lms d\log(f) \wedge d\log(g). 
$$

\paragraph{The Bloch complex.} It is clear from (\ref{kk}) that the Milnor group $K_2(F)$ is the cokernel of the map
$$
\delta: \Z[F^*-\{1\}] \lra {\bigwedge}^2F^*, \qquad \{x\} \lms (1-x) \wedge x.
$$
where $\{x\}$ is the generator of $\Z[F^* - \{1\}]$ corresponding to $x \in F^*- \{1\}$.

It turns out that one can describe nicely 
some elements in the kernel of this map parametrised 
by connected varieties of dimension bigger then zero. 
Namely, let $r(*,*,*,*)$ be the cross-ratio of four points on ${\Bbb P}^1$, normalized by 
$r(\infty, 0, 1, x) = x$. Let $R_2(F)$ be the 
subgroup of $\Z[F^*-\{1\}]$ generated by the following elements 
(the ``five term relations''):
$$
\sum_{i=1}^5 (-1)^i\{r(x_1, ..., \widehat x_i, ..., x_5)\}, \qquad x_i \in P^1(F), 
\quad x_i \not = x_j. 
$$
Then one can check that $\delta (R_2(F)) =0$ (see \cite{G1} for a conceptual proof). 
Let us set
$$
B_2(F):= \frac{\Z[F^*-\{1\}]}{R_2(F)}.
$$
Then the map $\delta$ gives rise to a homomorphism 
of $F$: 
$$
\delta: B_2(F) \lra {\bigwedge}^2F^*.
$$
We view it as a complex, called the Bloch complex.  
Let $\{x\}_2$ be  the 
projection of $\{x\}$  to $B_2(F)$.   
It is handy to add the elements $\{0\}_2$, $\{1\}_2$, $\{\infty\}_2$ and put 
  $\delta \{0\}_2 =
\delta \{1\}_2 = \delta \{\infty\}_2 =0$.

Recall that $K_1(F) = F^*$. The product in Quillen's K-theory provides a 
  map
$$
K_1(F) \otimes K_1(F) \otimes K_1(F) \lra K_3(F).  
$$ 
The cokernel of this map is the
{\it indecomposable part} $K_3^{\rm ind}(F)$  of $K_3(F)$. 
By  Suslin's theorem \cite{Sus} there is an isomorphism 
\begin{equation} \label{SUS}
{\rm Ker}\Bigr(B_2(F) \stackrel{\delta}{\lra}  {\bigwedge}^2F^*\Bigl) \otimes \Q
\stackrel{\sim}{= }
K_3^{\rm ind}(F)\otimes \Q. 
\end{equation}
%By Beilinson's rigidity conjecture, an
%embedding $\overline \Q \hra \C$ should give rise to  an isomorphism 
%$$
%K_3^{\rm ind}(\overline \Q)\otimes \Q \stackrel{\sim}{\lra}  
%K_3^{\rm ind}(\C)\otimes \Q. 
%$$

\paragraph{The dilogarithm.} 
Recall the classical dilogarithm function
$$
{\rm Li}_2(z):= 
-\int_0^z \log(1-t) d\log t.
$$
The single-valued cousin of the dilogarithm, the Bloch-Wigner function 
$$
{\cal L}_2(z):= {\rm Im}\Bigl({\rm Li}_2(z) + \arg(1-z)\log|z|\Bigr) 
$$
satisfies the Abel five term functional equation
$$
\sum_{i=0}^4(-1)^i {\cal L}_2(r(x_0, ..., \widehat x_i, ..., x_4)) = 0.
$$ 
Equivalently,  it provides a group homomorphism 
$$
{\cal L}_2: B_2(\C) \to \R, \qquad \sum n_i\{z\}_2 \lms \sum n_i{\cal L}_2(z_i).
$$

\subsection{The $W$-- class in $K_2$} 
Recall that a seed ${\mathbf i}$ provides cluster coordinates 
 $\{A_i\}$ on  
 ${\cal A}$. 
Set 
\begin{equation} \label{5.20.03.4}
W_{\mathbf i}:= \sum_{i,j \in I}\widetilde \varepsilon_{ij} \cdot A_i \wedge A_j  
\in {\bigwedge}^2\Q({\cal A})^*, \quad \widetilde \varepsilon_{ij}:= d_i 
\varepsilon_{ij}.
\end{equation}
Let us consider a map $$
\pi: {\cal A} \lra {\cal X} \times {\cal A}
$$
 given as a composition
$$
{\cal A} \hookrightarrow {\cal A} \times {\cal A} 
\stackrel{p\otimes {\rm Id}}{\lra} {\cal X} \times {\cal A}. 
$$ 
Here the first map is the diagonal map, 
and the second is the map $p$ on the first factor.

\bl \la{IF1} 
The element 
$W_{{\bf i}}$ is the projection of $\pi^*{\bf P}_{{\bf i}}$ to 
$\bigwedge^2\Q({\cal A})^*$. 
\el

{\bf Proof}. One can write the element $W_{\mathbf i}$ as follows:
\be \la{IF}
W_{\mathbf i} = \sum_{i \in I}d_i \cdot A_i \wedge p^*X_i. 
\ee
Indeed, 
$$
\sum_{i,j \in I}\widetilde \varepsilon_{ij} \cdot A_i \wedge A_j = 
\sum_{i,j \in I}d_i  \cdot A_i \wedge A_j^{\varepsilon_{ij}} = 
\sum_{i \in I}d_i \cdot A_i \wedge p^*X_i.
$$ 
Lemma follows immediately from this.

\begin{proposition} \label{5.11.03.7} Let $\mu_k: {\mathbf i} \to {\mathbf i}'$ 
be a mutation.  Then one has 
$$
\mu_k^*W_{{\mathbf i}'} - W_{\mathbf i} = -2d_k \cdot p^*\Bigl ((1+X_k)\wedge X_k  
\Bigr) = -2d_k\cdot \delta\{-p^*X_k\}_2. 
$$ 
\end{proposition}

{\bf Proof}. Follows from Basic Lemma \ref{BL} and Lemma \ref{IF1}.

\begin{corollary}\label{5.11.03.8}
The element
$$
W = \sum_{i,j \in I}\widetilde \varepsilon_{ij}\cdot 
 \{A_i,  A_j\} \in 
K_2(\Q({\cal A})) 
$$ 
does not depend on the choice of the 
cluster coordinate system $\{A_i\}$. 
\end{corollary}

{\bf Proof}. Observe that 
$2 \cdot (1+x) \wedge x = 2 \cdot (1-(-x))\wedge (-x)$ is twice 
a Steinberg relation. 

\begin{corollary}\label{5.11.03.28}
The $2$--form 
$$
\Omega= \sum_{i,j \in I}\widetilde 
\varepsilon_{ij}\cdot  d\log A_i\wedge  d\log A_j
$$ 
does not depend on the choice of 
cluster coordinates $\{A_i\}$. 
\end{corollary}

\begin{lemma} \label{5.11.03.7fd} 
The element $W$ is lifted from ${\cal U}$, i.e. one has 
$
W \in p^*\bigwedge^2\Q({\cal U})^*.
$ 
\end{lemma}

{\bf Proof}. This is equivalent to the following.  
Let $\{b_j\}$ be a set of integers such that $\sum_j 
\varepsilon_{ij}b_j=0$. Then replacing  $A_i$ by
$\lambda^{b_i}A_i$  we do not change the class $W_{\mathbf i}$.  Let us check this
claim. We have, using the skew symmetry of  $\varepsilon_{ij}$:
$$
\sum_{i,j \in I}\widetilde \varepsilon_{ij} \cdot (\lambda^{b_i} A_i) \wedge 
(\lambda^{b_j} A_j) - \sum_{i,j \in I}\widetilde \varepsilon_{ij} \cdot A_i \wedge 
A_j = 2\sum_{i \in I}\sum_{j \in I}  \widetilde \varepsilon_{ij}b_j \cdot    A_i \wedge \lambda = 0. 
$$ 
The last equality follows from the condition $\sum_j \widetilde
\varepsilon_{ij}b_j =0$, which is equivalent to $\sum_j 
\varepsilon_{ij}b_j=0$.  The lemma is proved. 

\vskip2mm

A cluster coordinate $A_k$ provides a valuation $v_{A_k}$
of the field  $\Q({\cal A})$ given by 
$v_{A_k}(A_j) = \delta_{jk}$. 

\begin{lemma}\label{5.11.03.38}
The element $W$ has zero tame symbol with respect to the valuation 
$v_{A_k}$. 
\end{lemma}

{\bf Proof}. The tame symbol for the valuation 
$v_{A_k}$ equals to
$$
\prod_{j\in I} A_j^{\widetilde \varepsilon_{kj}} = 
\Bigl(\prod_{j\in I} A_j^{\varepsilon_{kj}}\Bigr)^{d_k} = p^*X_k^{d_k}.
$$ 
One has 
$
1+ X_k = ({\Bbb A}_k^+ + {\Bbb A}_k^-)/{\Bbb A}_k^- = A_k A_k'/{\Bbb A}_k^-.
$ 
Thus $1+ X_k=0$ if $A_k =0$.
The lemma is proved.

\vskip 2mm
{\bf Remark}. {\it Over local fields}. 
The $K_2$-class $W$ has the following amusing application. 
Let $F$ be a local field. Let $\mu_F$ 
be the group of roots of unity contained in  $F$. 
Let $\alpha: K_2(F) \lra \mu_F$ be the norm residue symbol 
([Mi], ch. 15).
It is a (weakly) continuous  function on $F \times F$ ([Mi], ch. 16). 
Clearly there is a restriction $i_x^*W \in K_2(F)$ of $W$ to any $F$-point $x$ 
of the union of the cluster tori. 
So we get 
a continuous function 
$$
h_F: {\cal A}(F) \lra \mu_F; \quad x \in {\cal A}(F) \lms \alpha(i_x^*W) \in \mu_F.
$$
 If $F = \R$,  we get a continuous function 
$h_{\R}: {\cal A}(\R) \lra  \Z/2\Z$. Its value on ${\cal A}(\R_{>0})$ is $+1$.

\begin{lemma} \label{11.14.03.3b}
The  action of the torus ${H}_{\cal A}$ on the space 
${\cal A}$ preserves the class $W$ in $K_2$. The orbits of this action 
coincide with the null-foliation of the $2$-form $\Omega$. 
\end{lemma}

{\bf Proof}. Follows from the very definitions 
and Lemma \ref{1.18.04.1}.

\begin{corollary} \label{11.14.03.3bc}
The   symplectic structure on the positive space ${\cal U}$ provided
by the form $\Omega$ on the space ${\cal A}$ coincides with the
symplectic structure 
induced by the Poisson structure on the space ${\cal X}$.
\end{corollary}

{\bf Proof}. This is an immediate consequence of  Lemma \ref{1.18.04.1}, 
 Proposition \ref{1.8.04.12}, and Lemma  \ref{11.14.03.3b}.

\vskip3mm

{\it A degenerate symplectic structure on the space of real tropical ${\cal
    A}$-space}. Let $\{a_i\}$ be the coordinates corresponding to a
    seed ${\mathbf i}$. Consider the $2$-form 
$
\omega:= \sum_{i,j \in I}\widetilde \varepsilon_{ij}da_i \wedge da_j
$ 
on the space ${\cal A}(\R^t)$. 
Since mutations are given by piece-wise linear transformations, it makes sense 
to ask whether this form is invariant under mutations. 
The following easy Lemma follows from Proposition \ref{5.11.03.7}. 
\begin{lemma} \label{7.8.04.15}
The form $\omega$ does not depend on the choice of a cluster coordinate
system.  
\end{lemma} 
So 
the $2$-form $\omega$ provides the real tropical space ${\cal A}(\R^t)$ 
with a degenerate symplectic
structure invariant under the modular group.

\subsection{The motivic dilogarithm class} 

\paragraph{The weight two motivic complexes.} 
Let $X^{(k)}$ be the set of all codimension $k$ irreducible subvarieties of $X$. 
Recall the tame symbol  map 
\begin{equation} \label{12.12.02.11we}
{\rm Res}: {\bigwedge}^2\Q(X)^* \stackrel{}{\lra}  \prod_{Y \in X^{(1)}}\Q(Y)^*; \quad 
f\wedge g \lms {\rm Rest}_Y\Bigl(f^{v_Y(g)}/g^{v_Y(f)}\Bigr)
\end{equation}
where $v_Y(f)$ is the order of zero of a rational function $f$ at the generic point of an 
irreducible divisor $Y$, and ${\rm Rest}_Y$ denotes restriction to 
$Y$. 
The {\it weight two motivic complex} $\Gamma(X;2)$ of a regular irreducible 
  variety $X$ with the field of functions $\Q(X)$ is  the following 
complex of abelian groups: 
\begin{equation} \label{12.12.02.11}
\Gamma(X;2):= \qquad B_2(\Q(X)) \stackrel{\delta}{\lra} {\bigwedge}^2 \Q(X)^* \stackrel{\rm Res}{\lra}  \prod_{Y \in X^{(1)}}\Q(Y)^*
\stackrel{{\rm div}}{\lra}  \prod_{Y \in X^{(2)}}\Z.
\end{equation}
where the first group is in degree $1$ and ${\rm Res}$ is the tame symbol map 
(\ref{12.12.02.11we}). If $Y \in X^{(1)}$ is normal, 
the last map 
is given by the divisor ${\rm div}(f)$ of $f$. 
If $Y$ is not normal,  we take its normalization $\widetilde Y$, compute ${\rm div}(f)$ 
on $\widetilde Y$, 
and then push it down to $Y$.  
The rational cohomology of this complex of groups is 
the weight two 
motivic cohomology of the scheme $X$:
$$
H^i(X, \Q_{\cal M}(2)) := H^i\Gamma(X;2)\otimes \Q.
$$ 
The complex (\ref{12.12.02.11}) is a complex of global sections of a complex 
of acyclic 
sheaves in the Zariski topology. We denote this complex of sheaves by 
$\underline \Gamma(2)$.

\paragraph{Background on equivariant cohomology.}  
Let $P$ be an oriented polyhedron, possibly infinite. 
We denote by $V(P)$ the set of its vertices. 
Suppose that we have a covering ${\cal U} = \{U_{v}\}$  of a scheme $X$ 
by Zariski open subsets, parametrized by the set $V(P)$ of vertices $v$ 
of $P$. Then for every face $F$ of $P$ there is 
a Zariski open subset 
$
U_F:= \cap_{v \in V(F)} U_{v}.
$ 
For every inclusion $F_1 \hra F_2$ of faces of $P$ there is an 
embedding $j_{F_1, F_2}: U_{F_1} \hra U_{F_2}$. 
In particular there is an embedding $j_{F}: U_F \hra X$. 
So we get a diagram of schemes $\{{\cal U}_F\}$, 
whose objects are parametrized by the faces of $P$, and arrows 
correspond to codimension one inclusions of the faces. 

Given  a complex of sheaves ${\cal F}^{\bullet}$ on $X$, the 
diagram $\{{\cal U}_F\}$ provides a bicomplex 
${\cal F}_{{\cal U}, P}^{\bullet,\bullet}$ defined as follows. For any  integer $k\geq 0$, let 
${\cal F}_{{\cal U}, P}^{\bullet,k}$ be the following direct sum of the complexes of sheaves:
$$
{\cal F}_{{\cal U}, P}^{\bullet,k}:= \oplus_{F: 
{\rm codim}F =k} ~{j_{F}}_*{j_{F}}^*{\cal F}^{\bullet}.
$$
The second differential 
$d_2: {\cal F}_{{\cal U}, P}^{\bullet,k}\lra {\cal F}_{{\cal U}, P}^{\bullet,k+1}$ is a sum 
of  
restriction morphisms $j_{F_1, F_2}^*$ where ${\rm codim}~F_1 =k$ and 
${\rm codim}F_2 = k+1$ with the signs reflecting the 
orientations of the faces of $F$. 

Assume that a group 
$\Gamma$ acts on  $X$, and it acts 
freely on a contractible polyhedron $P$, preserving its 
polyhedral structure. Moreover, we assume that 
for any $\gamma \in \Gamma$ and any face $F$ of $P$ one has 
$\gamma({\cal U}_F) = {\cal U}_{\gamma(F)}$. Finally, 
let us assume that ${\cal F}^{\bullet}$ is a $\Gamma$-equivariant complex of sheaves on $X$. 
Then the group $\Gamma$ acts freely on the bicomplex 
${\cal F}_{{\cal U}, P}^{\bullet,\bullet}$.  If 
the sheaves $j^*_{U_F}{\cal F}^i$ are acyclic -- 
 this will be the case below -- then 
the $\Gamma$-equivariant hypercohomology  of $X$ 
with coefficients in ${\cal F}^{\bullet}$ are computed as follows: 
\begin{equation} \label{8.25.04.1}
H^*_{\Gamma}(X, {\cal F}^{\bullet})= 
H^*\left({\rm Tot}({\cal F}_{{\cal U}, P}^{\bullet,\bullet} )^{\Gamma}
\right).
\end{equation}
Here the right hand side has the following meaning: 
we take the total complex of the bicomplex ${\cal F}_{{\cal U}, P}^{\bullet,\bullet}$, 
take its subcomplex of $\Gamma$-invariants and compute its cohomology. 
Formula (\ref{8.25.04.1}) holds modulo $N$-torsion if 
the group $\Gamma$ acts on $P$ with finite stabilisers, whose orders 
divide $N$.

\paragraph{The motivic dilogarithm class.} 
The  special cluster modular group $\widehat \Gamma$ acts on 
 ${\cal  U}$. 
Let us apply the above construction  
in the case when $X = {\cal  U}$ and the cover ${\cal U}$ is given by 
(the images) of the cluster tori:

\begin{definition} \label{} The second integral weight two 
$\widehat \Gamma$-equivariant motivic cohomology group of ${\cal  U}$, 
denoted 
$%\begin{equation} \label{12.9.03.100}
H_{\widehat \Gamma}^2({\cal U}, \Z_{\cal M}(2)), 
$ %\end{equation} 
is obtained by  the construction of Section  5.4
in the following situation: 

1. The scheme $X$ is ${\cal U}$.

2. The group $\Gamma$ is the special cluster modular group $\widehat \Gamma$.

3. The  polyhedral complex 
$P$ is the modular complex $\widehat M$ 
to which we glue cells of dimension $\geq 3$ to make it contractible. 

4. The complex of sheaves ${\cal F}^{\bullet}$ is the complex of sheaves 
$\underline \Gamma(2)$ 
on ${\cal U}$. 
\end{definition}

{\bf Remarks}. i) One easily sees from the construction that the group 
$H_{\widehat \Gamma}^2({\cal U}, \Z_{\cal M}(2))$ does not 
depend on the choice of the cells of dimension $\geq 3$ glued to 
the universal cover on $M$.

ii) One can not define the 
$H_{\widehat \Gamma}^i({\cal U}, \Z_{\cal M}(2))$ for $i<2$ in a 
similar way using the complex (\ref{12.12.02.11}), since it  computes only the rational, 
but not  the integral 
motivic cohomology in the degree $1$.

\begin{theorem} \label{5.20.03.1}
There is   a class 
$
{\Bbb W} \in H_{\widehat \Gamma}^2({\cal U}, \Z_{\cal M}(2)).
$ 
\end{theorem}

{\bf Proof}. A construction of a cocycle representing a cohomology class 
on the right boils down to the following procedure.

i) For each seed  ${\mathbf i}$ exhibit a class $W \in 
\bigwedge^2\Q({\cal U})^*$ 
such that for every irreducible divisor $D$ 
in ${\cal U}$ the tame symbol 
of $W_{\mathbf i}$ at $D$ vanishes. 

ii) To any mutation  ${\mathbf i}\to {\mathbf i}'$ in $\widehat M$ find 
an element 
\begin{equation} \label{5.20.03.2}
\beta_{{\mathbf i} \to  {\mathbf i}'}\in B_2(\Q({\cal U}))\quad \mbox{ such that} \quad  
\delta \Bigl( \beta_{{\mathbf i} \to  {\mathbf i}'}\Bigr)  =   W_{{\mathbf i}} - W_{{\mathbf i}'}
\end{equation}

iii) Prove that for any $2$-dimensional cell of  
$\widehat M$ the sum of the elements  (\ref{5.20.03.2}) assigned to the 
oriented 
edges of 
the boundary 
of this cell is zero.

Recall the element 
$W_{\mathbf i} \in \bigwedge^2\Q({\cal U})^*$ assigned to a seed ${\mathbf i}$, 
see (\ref{5.20.03.4}). 
For a mutation $\mu_k: {\mathbf i} \to {\mathbf i}'$  set 
$$
\beta_{{\mathbf i} \to  {\mathbf i}'}:= -2d_k\cdot \{-X_k\}_2 \in B_2(\Q({\cal U})).
$$ 
Then Proposition \ref{5.11.03.7} is
 equivalent to
(\ref{5.20.03.2}). The $2$-cells in $\widehat M$ are the standard   
$(h+2)$-gons. Therefore we have to check the statement only for the 
cluster transformations (\ref{K10}), which are
 given by $(h+2)$-fold composition of the recursion (\ref{5.22.03.1a}). 
Therefore iii) reduces to the following

\begin{lemma}\label{5.24.03.2}
Let us define elements $x_i$ of a field $F$ by
recursion (\ref{5.22.03.1a}). Then if the sequence $\{x_i\}$ 
is periodic with the period $h+2$, 
one has 
\begin{equation} \label{5.24.03.4}
\sum_{i=1}^{h+2}d_i\{-x_i\}_2 =0 
 \quad \mbox{in $B_2(F)\otimes \Q$}, \qquad\qquad d_i =  \left\{ \begin{array}{ll}b & 
\mbox{ $i$: even }\\
c
& \mbox{ $i$: odd. }\end{array}\right.
\end{equation}
\end{lemma}

{\bf Proof}. Suppose that we have rational
functions $f_i$ such that 
$\sum (1-f_i) \wedge f_i =0$. Then (see \cite{G1}) if there is a point 
$a\in F$ such that $\sum \{f_i(a)\}_2 =0$, then $\sum \{f_i(x)\}_2 =0$ for any
$x \in F$.  

Formula (\ref{5.24.03.4}) in the case $A_2$ 
is equivalent to the famous five term  relation for the dilogarithm. 
In the case $A_1\times A_1$ it is  the well known 
inversion relation $2(\{x\}_2 + \{x^{-1}\}_2) =0$ in the group $B_2(F)$. 
It is easy to check using the recursion (\ref{5.22.03.1a}) 
that $\delta \sum_{i=1}^{h+2}d_i\{-x_i\}_2 =0$ in 
$\bigwedge^2F^*$. (Modulo $2$-torsion this follows from 
Proposition \ref{5.11.03.7}). Specializing $x_1=-1$ in the $G_2$ case we get 

$$
3 \{1\}_2 + \{-x_2\}_2 + 3 \{1+x_2\}_2 + \{x^2_2\}_2 + 3 \{1-x_2\}_2 +
\{x_2\}_2 +3 \{1\}_2 + \{0\}_2.  
$$
Specializing further $x_2 =0$ we get $12\{1\}_2 + 4\{0\}_2 =0$. 
 The $B_2$ case is similar. The lemma and hence 
the theorem are proved.

\vskip2mm

{\bf Exercise}. The lemma asserts that the element on the left of (\ref{5.24.03.4})
can be presented as linear combinations of the 
five term relations. Write them down in the $B_2$ and $G_2$ cases.

 \begin{corollary}\label{5.24.0}
Let us define elements $x_i \in \C$ by
recursion (\ref{5.22.03.1a}). 
Then one has 
\begin{equation} \label{5.24.03.4a}
\sum_{i=1}^{h+2}d_i{\cal L}_2(-x_i) =0. 
\end{equation}
\end{corollary}

{\bf Proof}. Follows from Lemma \ref{5.24.03.2} and and the five-term
functional equation for the dilogarithm. 

\vskip3mm

\begin{corollary}\label{5.24.0fr}   
A path $\alpha: {\mathbf i} \to {\mathbf j}$ in 
$\widehat {\cal G}$ provides an element 
$
\beta_{\alpha} \in B_2(\C)$  such that  
$\delta ( \beta_{\alpha})  =   W_{{\mathbf i}} - W_{{\mathbf j}}.
$ 
\end{corollary}

{\bf Proof}. Decomposing the path $\alpha$ as a composition of mutations 
${\mathbf i}= {\mathbf i}_1 \to {\mathbf i}_2 \to ... \to {\mathbf i}_n= {\mathbf j}$, set 
$$
\beta_{\alpha}:= \sum_{i=1}^{n-1}\beta_{{\mathbf i}_i \to  {\mathbf i}_{i+1}}\in B_2(\C).
$$
The relations in the groupoid 
$\widehat {\cal G}$ are generated by the ones corresponding to the standard   
$(h+2)$-gons. Thus Lemma \ref{5.24.03.2} implies that this element 
does not depend on the 
choice of a decomposition.  It evidently satisfies formula (\ref{5.20.03.2}). 
The corollary is proved.

\subsection{Invariant points of the modular group
  and $K_3^{\rm ind}(\overline \Q)$} 
Let $g$ be an element of the group $\widehat \Gamma$. It acts 
  by an automorphism of the scheme ${\cal U}$. 
It follows that any stable point of $g$ is defined over $\overline \Q$: 
it is determined by a set of equations with rational coefficients. 
Let $p\in {\cal U}(\overline \Q)$ be a stable point of 
  $g$. The element  $g$ can be presented by 
a loop $\alpha(g)$ based at $v$. 

\begin{proposition} Let $p\in {\cal U}(\overline \Q)$ be a stable point of 
  an 
  element $g \in \widehat \Gamma$. 
Then there is an invariant
$$
\beta_{g, p}:= \beta_{\alpha(g)} \in K_3^{\rm ind}(\overline \Q)\otimes \Q.
$$
\end{proposition}

{\bf Proof}. Since $p \in {\cal U}$, the element $W_{\bf i}$ can be evaluated 
at  $p$.  Corollary \ref{5.24.0fr}, applied to the  
loop $\alpha(g)$ based at $v$, implies 
$
\delta \beta_{g, p} = W_{\mathbf i}(p) - W_{{\mathbf i}}(p) = 0.
$ 
So $\beta_{\alpha(g)} \in H^1(B(\overline \Q; 2))$. 
Using (\ref{SUS}) get an element in 
$K_3^{\rm ind}(\overline \Q)\otimes \Q$.  

\vskip2mm

Here is how we compute the invariant $\beta_{g, p}$. 
Let us present $g$ as a composition of mutations: $g = \gamma_1
  ... \gamma_n$. Each mutation $\gamma$ determines the corresponding
  rational function $X_{\gamma}$ on  ${\cal U}$. 
Then 
$$
\beta_{g,p}= \sum_{i=1}^n 2d_\gamma \{X_{\gamma_i}(p)\}_2.
$$
 Here $d_\gamma$ is the multiplier assigned to the 
cluster coordinate $X_\gamma$. A similar procedure can be applied to any stable point 
$p \in {\cal U}$ of $g$, assuming that the functions 
$X_{\gamma_i}$ can be evaluated at $p$. 

\vskip 2mm
{\bf Remark}. According to Borel's theorem, the rank of 
$K_3^{\rm ind}(F)$ for a number field $F$ equals to the number $r_2(F)$ of 
embeddings $F\subset \C$ up to complex conjugation.

\end{document}